\documentclass[11 pt, reqno]{article}

\usepackage{amsmath,amsfonts,amssymb,amsthm}
\usepackage[colorlinks=true,hyperindex=true]{hyperref}
\usepackage{cancel,bbm}
\usepackage{mathabx} 
\usepackage{comment}
\usepackage{enumitem}
\usepackage{cite}

\usepackage{chngcntr}
\counterwithin*{equation}{section}

\usepackage[a4paper, left=2cm, right=2cm, top=2 cm, bottom= 2 cm]{geometry}
\usepackage{color}
\usepackage{fancyhdr}
\usepackage{latexsym}

\newtheorem{ccounter}{ccounter}[section]
\newtheorem{thm}[ccounter]{Theorem}
\newtheorem{lem}[ccounter]{Lemma}
\newtheorem{cor}[ccounter]{Corollary}

\newtheorem{prop}[ccounter]{Proposition}
\newtheorem{ass}[ccounter]{Assumption}

\theoremstyle{definition}
\newtheorem{defn}[ccounter]{Definition}
\newtheorem{ex}[ccounter]{Example}

\def\bet{\begin{thm}}
\def\eet{\end{thm}}
\def\bel{\begin{lem}}
\def\eel{\end{lem}}
\def\bas{\begin{ass}}
\def\eas{\end{ass}}
\def\bec{\begin{cor}}
\def\eec{\end{cor}}
\def\bed{\begin{defn}}
\def\eed{\end{defn}}
\def\bep{\begin{prop}}
\def\eep{\end{prop}}
\def\beq{\begin{equation}}
\def\eeq{\end{equation}}
\def\proof{\noindent {\bf Proof.}\ \ }
\def\bea{\begin{equation*}}
\def\eea{\end{equation*}}
\def\tr{\mathrm{tr}}
\def\bex{\begin{ex}}
\def\eex{\end{ex}}
\def\remark{\noindent{\bf Remark. }}

\def\rr{\mathbb{R}}

\def\cc{\mathbb{C}}
\def\1{\boldsymbol{1}}
\def\Im{\mathrm{Im}}
\def\Re{\mathrm{Re}}
\def\e{\mathrm{e}}
\def\i{\mathrm{i}}
\def\del{\partial}
\def\d{\mathrm{d}}
\def\eps{\varepsilon}
\renewcommand\leq\varleq
\renewcommand\geq\vargeq
\def\ee{\mathrm{E}}

\def\F{\mathcal{F}}
\def\O{\mathcal{O}}

\def\ee{\mathbb{E}}

\def\pp{\mathbb{P}}
\def\cb{c_\beta}
\def\tilm{\tilde{m}}
\def\hatgam{\hat{\gamma}}

\def\rhosc{\rho_{\mathrm{sc}}}

\def\dist{\stackrel{\mathrm{d}}{=}}
\def\mv{m_v}
\def\tilmv{\tilde{m}_v}
\def\kapgam{\kappa}

\def\F{\mathcal{F}}
\def\I{\mathcal{I}}

\def\G{\mathcal{G}}
\def\msc{m_{\mathrm{sc}}}
\def\hata{\hat{a}}
\def\hatgam{\hat{\gamma}}

\def\rhosc{\rho_{\mathrm{sc}}}
\def\S{\mathcal{S}}

\def\ulj{\underline{j}}
\def\P{\mathcal{P}}
\def\hats{\hat{s}}
\def\ss{\mathbb{S}}
\def\J{\mathcal{J}}
\def\Dk{\mathcal{D}_{\mathfrak{m}}}
\def\Gamk{\Gamma_\mathfrak{m}}
\def\tilGamk{\tilde{\Gamma}_{\mathfrak{m}}}
\def\tilGam{\tilde{\Gamma}}
\def\tilf{\tilde{f}}
\def\tillam{\tilde{\lambda}}
\def\Feps{\mathcal{F}_{\varepsilon_1}}
\def\xb{x_b}
\def\kapp{\mathfrak{m}}

\def\gamu{\tilde{\gamma}_u}
\def\gam1{\tilde{\xi}_1}
\def\gamt{\tilde{\gamma}}
\def\etat{\tilde{\eta}}
\def\gamh{\tilde{\xi}}
\def\gammu{\xi_u}

\begin{document}

\begin{table}
\centering

\begin{tabular}{c}
\multicolumn{1}{c}{\Large{\bf Fluctuations of the 2-spin SSK model with magnetic field}}\\
\\
\\
\end{tabular}
\begin{tabular}{c c c}
Benjamin Landon & \phantom{blblbl} & Philippe Sosoe\\
\\
\small{Department of Mathematics} & & \small{Department of Mathematics} \\
\small{University of Toronto} & & \small{Cornell University} \\
\small{blandon@math.toronto.edu} &  & \small{ps934@cornell.edu}  \\
\\
\end{tabular}
\\
\begin{tabular}{c}
\multicolumn{1}{c}{\today}\\
\\
\end{tabular}

\begin{tabular}{p{15 cm}}
\small{{\bf Abstract:} We analyze the fluctuations of the free energy, replica overlaps, and overlap with the external field in the quadratic spherical SK model with a magnetic field. We identify several different behaviors for these quantities depending on the size of the magnetic field, confirming predictions by Fyodorov-Le Doussal and recent  work of Baik, Collins-Wildman, Le Doussal and Wu.}
\end{tabular}
\end{table}


\section{Introduction}
In this paper, we study the 2-spin spherical Sherrington-Kirkpatrick (SSK) model with external magnetic field. This is the random Gibbs measure on the $N-1$ dimensional sphere of radius $\sqrt{N}$ (denoted by $\ss^{N-1} := \{ \sigma \in \rr^N : ||\sigma||_2 = \sqrt{N} \}$), 
 with Hamiltonian given by
\beq
H_N ( \sigma ) := \frac{1}{ \sqrt{2N}} \sum^N_{i, j=1, i \neq j} g_{ij} \sigma_i \sigma_j + h v^T \sigma,
\eeq
where $h \in \rr$, and $v \in \rr^N$ is a deterministic vector such that
$
\|v\|_2^2 = N,
$
and the $\{ g_{ij}\}_{i, j=1}^N$ are iid standard normal random variables. The partition function is given by an integral over $\ss^{N-1}$ and is defined by,
\beq
Z_{N,\beta,h}=\int_{\mathbb{S}^{N-1}} \exp(\beta H_N(\sigma))\, \frac{ \d \omega_{N-1}(\sigma)}{ | \ss^{N-1} | }.
\eeq
Here, $\omega_{N-1}$ is the uniform measure on $\ss^{N-1}$.

\paragraph{Model with zero magnetic field.} The SSK (with zero magnetic field $h=0$) was introduced by Kosterlitz, Thouless and Jones \cite{KTJ} by analogy with the Sherrington-Kirkpatrick model with Ising spins. The spherical model with quadratic spins which we study here turns out to have a rather different behavior from the SK model with Ising spins or the $p$-spin models with $p\geq 3$ (including such higher spin spherical models). This is due to the fact that the quadratic nature of the Hamiltonian and the continuous state space reduce the complexity of the \emph{energy landscape}, which is not exponential as in other ``true'' spin glasses. Indeed, spin glasses are commonly characterized by the existence of an exponential (in $N$) number of critical points of the Hamiltonian for generic realizations of the disorder, while in the quadratic case on the sphere, this number is at most linear and indeed is $O(1)$ as the size of the magnetic field increases to order 1. See \cite{auffinger-cpx} Theorem 2.8 and Remark 2.11, as well as \cite{FD}, Section 3. These differences with other spin glasses are also reflected in our methods, which come largely from random matrix theory.

Nevertheless, the model exhibits an interesting phase transition at the critical temperature $\beta=1$, which was already identified by Kosterlitz-Thouless-Jones. The SSK without magnetic field (including higher spin models) was studied by Panchenko and Talagrand in \cite{panchenko-talagrand}. Baik and Lee \cite{baiklee} used  results in random matrix theory including the local semicircle law and the method of steepest descent to obtain second order asymptotics for the free energy. They later applied the same methods to analyze a number of variants of the spherical model, including the bipartite SSK model as well as a model including a deterministic Curie-Weiss term in the Hamiltonian \cite{baiklee1,baiklee2,BLW}. The latter model, although its definition is relatively simple, has an intricate phase diagram comprising three different regimes, which the authors of \cite{BLW} term ferromagnetic, paramagnetic and spin glass.

Baik and Lee \cite{baiklee} show that in the high temperature phase $\beta< 1$, the fluctuations of the free energy for the SSK with no magnetic field  are Gaussian. This is the spherical analog of the central limit theorem by Aizenman-Lebowitz-Fröhlich \cite{alr} for the  free energy in the SK model with Ising spins. In the low temperature regime, however, Baik and Lee found that the free energy fluctuations asymptotically follow the Tracy-Widom GOE distribution, the limiting distribution of the largest eigenvalue of a real symmetric random matrix consisting of Gaussian entries.

Following Baik and Lee's work, the second author and  Nguyen analyzed the overlap between two independent samples (known as \emph{replicas}) from the Gibbs measure in the case $h=0$ \cite{vl-ps}. They showed that in the high temperature region, the overlap,  properly rescaled, is asymptotically normally distributed.  The authors of the present paper then extended this result to the low temperature phase, where the asymptotic distribution is no longer Gaussian, but instead given by an explicit distribution which is a function of the Airy$_1$ random point field \cite{SSK-ls} (i.e., the joint limit of the largest eigenvalues of the GOE).

\paragraph{Model with magnetic field.} In the current work we investigate the effect of a non-zero magnetic field on the fluctuations of the partition function and overlap. In addition to the works in physics of Fyodorov-Le Doussal \cite{FD} and Cugliandolo-Dean-Yoshino \cite{cdy}, which we discuss in more detail below, the closest precedent to the current work in the mathematical spin glass literature is the works of Chen-Dey-Panchenko \cite{csp} and Chen-Sen \cite{chensen}.  The work \cite{csp} considers general $p$-spin models with Ising spins and $h =c > 0$, showing that the free energy fluctuations are Gaussian at any temperature; their methods apply also to spherical models.  The work \cite{chensen} shows for general even mixed $p$-spin models, that the fluctuations of the ground state are Gaussian if $h=c>0$. In particular for $p=2$, the Tracy-Widom fluctuations found by Baik and Lee for $h=0$ disappear.

In addition to the case of constant order magnetic field, we study also the case $h\sim N^{-\alpha}$ for $\alpha > 0$, that is, when the magnetic field  vanishes as the system grows. As $\alpha$ varies, the nature of the fluctuations of the free energy and overlaps also changes.  As mentioned, this situation was studied in some detail in the physics literature by Fyodorov-Le Doussal \cite{FD} (see also  Dembo-Zeitouni \cite{dembo-zeitouni}). Their work is formulated in terms of the constrained optimization problem
\begin{equation}\label{eqn: FD}
E_{\mathrm{min}}(h)=\min_{\sigma \in \mathbb{S}^{N-1}} \big(-H_N(\sigma)\big) ,
\end{equation}
which corresponds to the ``zero temperature'' regime $\beta =\infty$.

Fyodorov and Le Doussal identify two distinguished scaling regimes for the magnetic field $h$ in the low temperature phase $\beta> 1$:
\begin{enumerate}

\item \textbf{Microscopic magnetic field}, $h\sim N^{-1/2}$:  this regime is the boundary of the region of ``small'' magnetic field, in the sense that the fluctuations of the free energy about its first-order asymptotic limit are not affected by the presence of the magnetic term. These fluctations coincide with the case $h=0$ treated by Baik and Lee, and are given by the Tracy-Widom GOE distribution. 

\item \textbf{Moderate magnetic field}, $h\sim N^{-1/6}$: in this regime the magnetic field modifies the fluctuations asymptotics of the free energy, which are no longer given by the Tracy-Widom distribution.  Fyodorov and Le Doussal analyze the tail of the distribution of the minimum in \eqref{eqn: FD}. 
\end{enumerate}

In this paper we obtain several results confirming and extending the predictions of Fyodorov and Le Doussal \cite{FD}.  In addition to results about the fluctuations of the free energy, we also describe the fluctuations of the overlap (normalized inner product) between two independent samples from the Gibbs measure (``replicas'') as well as the overlap between a sample and the external field. This allows for a more refined description of each of the regimes found by Fyodorov and Le Doussal.   For example, our results show following: 
\begin{enumerate}
    \item For microscopic magnetic field, we confirm that the limiting fluctuations remain the same as in the case $h\neq 0$. Thus the Tracy-Widom fluctuations of the free energy proved by Baik and Lee in the $h=0$ case persist for values of $h \sim N^{-1/2}$. 
    Nevertheless, we show that the nonzero magnetic field can still be detected at the level of overlaps between replicas, whose distribution does differ from the case $h=0$.
    \item  Previous works \cite{dembo-zeitouni}, \cite{FD} only consider the far tail of the limiting distribution of the free energy density.
    We give various expressions for the limiting distribution in terms of quantities from random matrix theory. This includes an intermediate regime $h \sim N^{-1/6}$ and $\beta >1$ where the free energy fluctuations are neither Tracy-Widom or Gaussian. 
     We also show that for values of the magnetic field $h\gg N^{-1/6}$, the fluctuations of the free energy are Gaussian even in the low temperature regime $\beta > 1$.
    \item We also analyze the Gaussian in regime in great detail, finding that as long as the magnetic field is not too small (depending on how close $\beta$ is to the critical temperature), the fluctuations of quantities involving the overlaps are Gaussian. 
\end{enumerate}

The susceptibility of spherical spin glass models has been considered by  Cugliandolo-Dean-Yoshino \cite{cdy}.  Here, they consider two different cases of vanishing magnetic field, corresponding to whether one takes the limits $h \to 0$ or $N \to \infty$ first.  They find that these two limits agree in the high temperature regime, but differ at low temperature. 

While completing this work we learned that Baik, Collins-Wildman, Le Doussal and Wu \cite{BCDW} have completed a paper discussing various properties of a spherical model with random external field similar to the one we study here. The focus of this work is different from ours; it includes some formal computations at the physics level of rigor, but on the other hand it treats aspects of the model not considered in our work, including various transitions between scaling regimes and the geometry of the Gibbs measure. Some of the results in \cite{BCDW} had previously appeared in Wu's thesis \cite{wu}.

We also note that Kivimae  studied the fluctuations of the free energy at zero temperature $\beta=\infty$ in the regime $h=\O(N^{-1/6})$, identifying a family of distributions interpolating between Tracy-Widom and Gaussian \cite{pkiv}.

\subsection{Organization of paper}

The remainder of the paper is organized as follows.  We will first introduce much of the notation used throughout the paper in Section \ref{sec:notation}, organizing it in a single section for convenient reference.  The remainder of Section \ref{sec:mr} then contains our main results on the fluctuations of the SSK model.  We have split this into three subsections, one for each of the three scaling regimes we consider:  Section \ref{sec:mr-gauss} contains the high temperature and/or large magnetic field regime of Gaussian fluctuations; Section \ref{sec:mr-int} the regime of intermediate magnetic field and low temperature; Section \ref{sec:mr-micro} the regime of microscopic magnetic field and low temperature.  In Section \ref{sec:rmt}, we collect the results from random matrix theory that we will use in our paper.  Section \ref{sec:rep} collects the various contour integral representations we use to prove our results.   The remainder of the paper then consists of the proofs of the results stated in Section \ref{sec:mr}.  The high temperature or large magnetic field regime is analyzed in Section \ref{sec:gauss}.  The regimes of intermediate and microscopic magnetic fields at low temperature are considered in Sections \ref{sec:int} and \ref{sec:micro}, respectively.  Finally, the appendices contain some proofs of auxiliary results or those straightforward or repetitive enough to be omitted from the main body of the paper.

\paragraph{Acknowledgements.} The work of B.L. is supported by an NSERC Discovery grant.  The work of P.S. is partially supported by NSF grants DMS-1811093 and DMS-2154090.

\section{Main results and notation} \label{sec:mr}

Our results concern fluctuations of the free energy, the overlap between two independent copies of $\sigma \in \ss^{N-1}$ distributed according to the Gibbs measure associated to the SSK Hamiltonian (``replicas'') and the overlap between a replica and the external field $v$.  Our results are different depending on the scaling behavior of the magnetic field and inverse temperature.  We therefore organize our results into the three different scaling regimes.  The first regime is characterized by Gaussian fluctuations. In this regime, either the temperature is high or the magnetic field does not tend to $0$ too quickly -- this is quantified by the assumption \eqref{eqn:gaussianassumption} below. These Gaussian results are collected in Section \ref{sec:mr-gauss}. The other regimes are the cases of intermediate ($h \sim N^{-1/6}$) and microscopic ($h \sim N^{-1/2}$) magnetic fields and low temperature.  The results in this regimes are collected in Sections \ref{sec:mr-int} and \ref{sec:mr-micro}, respectively.  Before stating our results we organize the notation of the paper in the following section for convenient reference.

\subsection{Notation}
\label{sec:notation}

The Hamiltonian of the $2$-spin spherical Sherrington-Kirkpatrick model with magnetic field is 
\beq \label{eqn:hamiltonian-def}
H_N ( \sigma ) := \frac{1}{ \sqrt{2N}} \sum_{i \neq j } g_{ij} \sigma_i \sigma_j + h v^T \sigma,
\eeq
where 
$
\|v\|_2^2 = N
$
is a fixed, deterministic vector and the $\{ g_{ij}\}_{i < j }$ are a family of  iid standard normal random variables.  Above, the strength of the magnetic field is given by $h \in \rr$, which we may assume to be non-negative.  The phase space of our system is the $N-1$ dimensional sphere of radius $\sqrt{N}$ which is denoted by,
$
\ss^{N-1} := \{ \sigma \in \rr^N : \|\sigma\|_2^2 = N \}.
$
The partition function is then,
\beq
Z_{N, \beta, h} = \int \exp \left[ \beta H_N ( \sigma ) \right]  \frac{ \d \omega_{N-1} ( \sigma ) }{ | \ss^{N-1} | },
\eeq
where $\d \omega_{N-1} $ is the uniform surface measure on $\ss^{N-1}$ under which $\ss^{N-1}$ has volume,
$
\omega_{N-1} ( \ss^{N-1} ) = \frac{2  \pi^{ \frac{N}{2}}}{  \Gamma \left( \frac{N}{2} \right)} N^{ \frac{N-1}{2}}.
$
Given a Hermitian matrix $A$ we will denote its eigenvalues in decreasing order by $\lambda_i (A)$ and the associated eigenvectors by $u_i (A)$. 

 We will denote by $M$ the random matrix formed from the disorder variables in the SSK Hamiltonian via
\beq \label{eqn:M-def-f}
M_{ij} = - \frac{1}{ \sqrt{2N}} ( g_{ij} + g_{ji} ), \quad i \neq j
\eeq
and $M_{ii} = 0$.  In terms of the SSK Hamiltonian, we have
\beq
H_N ( \sigma ) = - \sigma^T M \sigma + h v^T \sigma.
\eeq

Up to the diagonal being $0$, the matrix $M$ is a matrix from the Gaussian Orthogonal Ensemble.  In fact, we will often compare the eigenvalues and eigenvectors of $M$ to those of a certain GOE matrix, which we denote by $H$.  The off-diagonal elements of $H$ are the same as $M$, and the diagonal entries are given by $H_{ii} = \sqrt{2} g_{ii} N^{-1}$ where $\{g_{ii} \}_{i=1}^N$ are independent standard normal random variables.  

Associated with $M$ are the following spectral quantities which describe the limiting density of states of the eigenvalues.  The semicircle distribution (being the limit of the empirical distribution of the eigenvalues of $H$ and $M$) and its Stieltjes transform will be denoted by,
\beq \label{eqn:rhosc-def-f}
\rhosc(E) := \frac{1}{ 2 \pi } \sqrt{ (4 - E^2 )_+ }, \qquad \msc (z) =\int \frac{ \rhosc (E) \d E }{ E - z }.
\eeq
Note that,
\beq \label{eqn:quadratic-msc}
\msc^2(z) + z \msc (z) +1 = 0, \qquad \msc' (z) = \frac{\msc(z)^2}{1 - \msc(z)^2} .
\eeq
The $N$-quantiles (also called classical locations) $\gamma_i$ of $\rhosc(E)$ are the points defined by
\beq \label{eqn:gamma-i-def}
\int_{\gamma_i}^{2} \rhosc (E) \d E = \frac{i}{N}.
\eeq
For the matrix $M$ we will denote,
\beq
m (z) := \frac{1}{N} \tr \frac{1}{M-z} = \frac{1}{N} \sum_{i=1}^N \frac{1}{ \lambda_i (M) - z }.
\eeq
We will also denote by $\mv (z)$ the following quadratic form in the resolvent,
\beq
\mv (z) := \frac{1}{N} v^T \frac{1}{M-z} v = \frac{1}{N} \sum_{i=1}^N \frac{ (v^T u_i (M) )^2}{  \lambda_i (M) - z },
\eeq
where $v$ will always be the external field as above. 
We will denote by $v_i$ the projections of the eigenvectors $u_i(M)$ onto $v$,
\beq
v_i := v^T u_i(M).
\eeq
Note that we will not need to refer to the components of $v$ in the coordinate basis, and so $v_i$ should not be mistaken for such quantities.

At one point we will need to separate the contribution of the largest eigenvalue of $M$ from the quantities $m(z)$ and $m_v (z)$ and so we denote,
\beq \label{eqn:m-tilde-def-f}
\tilm (z) := \frac{1}{N} \sum_{j=2}^N \frac{1}{ \lambda_j (M) - z } , \qquad \tilmv (z) := \frac{1}{N} \sum_{j=2}^N \frac{ ( v^T u_j (M))^2}{ \lambda_j (M) - z }.
\eeq

Our analysis of the various integrals over the sphere $\ss^{N-1}$ will proceed via the method of steepest descent.  In our application of this method we will have 
 use of the function $G(z)$, which is defined by
\beq \label{eqn:not-G}
G(z) := \beta z - \frac{1}{N} \sum_{i=1}^N \log (z - \lambda_i (M) ) - \frac{ h^2 \beta }{N} v^T \frac{1}{M-z} v,
\eeq
as well as its deterministic approximation, the function $g(z)$ defined by
\beq
g(z):= \beta z - \int \log(z - x ) \rhosc (x) \d x - h^2 \beta \msc ( z).
\eeq
We define here as well the saddles $\gamma$ and $\hatgam$ which are the unique solutions $\gamma >\lambda_1 (M)$  and $\hatgam >2$ satisfying
\beq \label{eqn:gam-def-f}
G' (\gamma) = 0, \qquad g' ( \hatgam ) = 0.
\eeq
The parameter
\beq \label{eqn:kap-def-f}
\kappa = \hatgam -2 >0
\eeq
will appear in many error estimates.

We use angular brackets to denote the Gibbs expectation,
\beq
\langle f ( \sigma ) \rangle = \frac{1}{Z_{N, \beta, h} } \int f ( \sigma ) \exp [ \beta H_N ( \sigma ) ] \frac{ \d \omega_{N-1} ( \sigma)}{ | \ss^{N-1} | }.
\eeq
The notation $\langle f ( \sigma_1, \sigma_2  ) \rangle$ indicates the expectation with respect to independent copies $ \sigma_1, \sigma_2$ of the Gibbs expectation.

Given two positive, possibly $N$-dependent quantities $a_N$ and $b_N$, the notation $a_N \asymp b_N$ means there is a constant $c>0$ so that
\beq
c a_N \leq b_N \leq c^{-1} a_N.
\eeq
The notation $a_N \ll b_N$ means $a_N / b_N \to 0$ as $N \to \infty$. We will only use the $\ll$ notation informally to motivate or explain sections of our proofs; in general we always have an explicit rate of convergence.

We will say that an event $\F$ holds with \emph{overwhelming probability} if for any $D>0$ we have $\pp [ \F] \geq 1- N^{-D}$ for all $N$ large enough.  Given a set of events $\F_i$ depending on a parameter $i \in \I$ in some possibly infinite index set $\I$, we say that the family $\{ \F \}_{i \in \I}$ holds with overwhelming probability if for any $D>0$ we have for all $N$ large enough that $\pp[ \F_i ] \geq 1 - N^{-D}$ for all $i \in \I$.

\subsection{High temperature or slowly decaying magnetic field (Gaussian regime)}
\label{sec:mr-gauss}

If either the inverse temperature $\beta$ satisfies $\beta <1$ or if the magnetic field $h$ is not too small, the fluctuations of the SSK model are in general Gaussian.  All of the results in this subsection are proven in Section \ref{sec:gauss}. 

\paragraph{Free energy.}

The leading order fluctuations to the free energy will be seen to be contributed by the quantity $G ( \hatgam)$. The first theorem below is our main result on the fluctuations of the free energy in the Gaussian regime. To state it, we define the following constant,  
\beq
C^{(g)}_N := \frac{1}{N} \log \Gamma(N/2) + \frac{1}{N} (N/2-1) \log(2 / (N \beta )) - \frac{1}{2N} \log ( N g'' ( \hatgam) \pi ) + \frac{ g ( \hatgam)}{2} .
\eeq
The following theorem is proven in Section \ref{sec:gauss-fe}.  As indicated in the introduction, when the magnetic field is of constant order, Gaussian fluctuations for spin glass models were proven by Chen-Dey-Panchenko \cite{csp}.  This case is covered by our result, which also allows for the magnetic field to tend to $0$ as $N \to \infty$. Essentially, there are two reasons that fluctuations might be Gaussian; first, if the temperature is above the critical temperature, $\beta=1$, and second, if the magnetic field is large (at any temperature). Both of these regimes are covered by our theorem. Note that in our results, $\beta$ and $h$ are allowed to depend on $N$. 
\bet[Free energy fluctuations] \label{thm:mr-gauss-fe}
We have the following in various parameter regimes:
\begin{enumerate}[label=\normalfont(\arabic*)]
\item Suppose that $c \leq \beta \leq 1-c$ and $h^2 \leq c^{-1}$ for some $c>0$. Then,
\beq
N \frac{ F_N - \hat{E}_N}{\hat{V}_N^{1/2}}
\eeq
converges to a standard random variable where $ \hat{E}_N$ and $\hat{V}_N$ are explicit functions of $v$, $\beta$ and $h$ defined later in \eqref{eqn:gaussian-variance-f}. We have $\hat{V}_N \asymp 1 + h^2 N$. 
\item Suppose that $c \leq \beta \leq c^{-1}$ and that 
\beq
c^{-1} \geq h^2 \geq N^{c} \left( N^{-1/2} (1-\beta)_+^{1/2}  + N^{-1/3} (\beta-1)_+ + N^{-2/3} \right)
\eeq
for some $c>0$. Then, 
\beq
2 \frac{N^{1/2} \kappa^{1/4}}{h^2 \beta V_N^{1/2} } (F_{N, \beta, h} -C^{(g)}_N )
\eeq
converges to a standard random variable where
\beq
V_N := \frac{ \hatgam + \sqrt{ \hatgam^2-4}}{ \sqrt{ \hatgam+2}} \msc(\hatgam)^4 \left( \msc^2 ( \hatgam) + (1 - N^{-2} \|v\|_4^4 ) (1 - \msc^2 ( \hatgam) ) \right).
\eeq
Moreover, $V_N \asymp 1$. 
\end{enumerate}
\eet

There is one high temperature regime not covered by our result, when $\beta \to 1^{-}$ as $N \to \infty$ at a sufficiently slow rate. In fact, we expect Gaussian fluctuations to hold whenever
\beq \label{eqn:gaussianassumption}
c^{-1} \geq (1 - \beta)_+  + \frac{ h^2}{ |1- \beta | + |h| } \geq N^{-1/3+\tau}, \qquad c \leq \beta \leq c^{-1}, \quad |h| \leq c^{-1}.
\eeq 
for some $c, \tau >0$. 
The regime $\beta \to 1^{-}$ (and $h \to 0$ sufficiently fast so that this is not covered by the second regime in the above theorem) requires a mesoscopic central limit theorem for non-compactly supported linear spectral statistics for random matrices near the spectral edge. Such a result is not available in the literature but would not be hard to prove with existing techniques; we refrain from doing so in the interest of space in this paper. The first part of the proof of the above theorem, Proposition \ref{prop:gaussianexpansion} that expands the free energy in terms of random matrix quantities in fact holds in this regime; it is only that the fluctuation results are not available in the random matrix literature. We expect that the first statement in the above theorem would hold without substantial change in this regime.

\paragraph{Overlap with external field.} Our next theorem concerns the overlap between a sample $\sigma$ drawn from the Gibbs measure and the external field $v$.  We obtain the following result for the Laplace transform of $\sigma \cdot v$, exhibiting quenched Gaussian fluctuations wrt the fixed Gibbs measure. 
\bet[Quenched external field overlap fluctuations] \label{thm:mr-gauss-ext}
Assume that \eqref{eqn:gaussianassumption} hold for a $\tau >0$ and $c>0$.  The following holds for some $c_1 >0$.   For any  $C>0$ we have with overwhelming probability that  uniformly for $|\lambda| \leq C$,
\begin{align}
 &\log \left(  \left\langle \exp \left[ \lambda N^{-1/2} v \cdot \sigma \right] \right\rangle\right) \nonumber \\
=& \lambda^2 \frac{ \msc' ( \hatgam)}{ 2 \beta g'' ( \hatgam ) } \left( - \msc ( \hatgam) + 2 h^2 \beta ( \msc' ( \hatgam) )^2 \right)  - \lambda \left(N^{-1/2}  h v^T (M - \gamma )^{-1} v \right) + \O(N^{-c_1} )
\end{align}
where $\gamma$ is the unique solution of $G' ( \gamma ) =0$ satisfying $\gamma > \lambda_1 (M)$.  The coefficient of $\lambda^2$ satisfies,
\beq
 \frac{ \msc' ( \hatgam)}{ 2 \beta g'' ( \hatgam ) } \left( - \msc ( \hatgam) + 2 h^2 \beta ( \msc' ( \hatgam) )^2 \right) \asymp 1.
\eeq
\eet  This theorem is proven in Section \ref{sec:gauss-ext}. 
The term linear in $\lambda$ corresponding to the quenched expectation of $N^{-1/2} v \cdot \sigma$ depends on the disorder.  Its fluctuations with respect to the disorder random variables are larger than the order $1$ fluctuations of $N^{-1/2} v \cdot \sigma$ with respect to the Gibbs measure, and so cannot be replaced by a deterministic quantity in the above statement.  Note also that it involves the random saddle $\gamma$ instead of its deterministic approximation $\hatgam$; the random saddle cannot in general be replaced by $\hatgam$ as the error is too large to still obtain a statement with an $o(1)$ error in the above estimate. 

We remark that the term linear in $\lambda$ is indeed, up to an $o(1)$ error term, the quenched expectation of the rescaled overlap of $\sigma$ with the magnetic field $v$, as can be seen from Lemma \ref{thm:mr-gauss-over-b1} that appears later in the paper. Therefore, the above result implies the asymptotic quenched Gaussian fluctuations for $N^{-1/2} v \cdot \sigma$ after recentering by its Gibbs expectation.

\paragraph{Overlap between replicas.} Our last major result for the Gaussian regime is the following theorem for the fluctuations of the overlap between two replicas. It again exhibits quenched Gaussian fluctuations wrt the fixed Gibbs measure. 
\bet[Quenched replica overlap fluctuations] \label{thm:gauss-overlap}
Assume that \eqref{eqn:gaussianassumption} holds for some $\tau >0$ and $c>0$.  There is a $c_1 >0$ so that the following holds.  For any $C>0$ we have with overwhelming probability that uniformly for $|\lambda| \leq C$,
\begin{align}
&\log \left( \left\langle \exp\left[ \lambda N^{-1/2} \kappa^{1/4} \sigma^{(1)} \cdot \sigma^{(2)} \right] \right\rangle \right) \nonumber \\
=& \frac{\lambda^2}{2 \beta } ( \msc' ( \hatgam) \kappa^{1/2} ) \frac{ \msc' ( \hatgam) - 2 h^2 \beta \msc'' ( \hatgam)}{ \msc' ( \hatgam) - h^2 \beta \msc'' ( \hatgam ) } + \lambda N^{-1/2} \kappa^{1/4} h v^T (M - \gamma)^{-2} v + \O (N^{-c_1} ).
\end{align}
\eet
This theorem is proven in Section \ref{sec:gauss-overlap}.  Again, the term linear in $\lambda$ is, up to an $o(1)$ error term, the quenched expectation of the rescaled overlap as can seen from Lemma \ref{thm:mr-gauss-b2} appearing later in the paper.  The coefficient of $\lambda^2$ is again $\asymp 1$. 

\subsection{Low temperature and intermediate magnetic field} \label{sec:mr-int}

The second scaling regime we find is when the temperature is low $(\beta >1)$ and the magnetic field scales like $h \sim N^{-1/6}$.  To be more precise, we assume
\beq \label{eqn:mr-int-ass}
1+c \leq \beta \leq c^{-1}, \qquad h^2 \beta = \theta N^{-1/3}
\eeq
for a fixed $c>0$ and fixed $\theta > 0$.  Note also that in the low temperature regime $ \beta > 1 + c$, the assumption \eqref{eqn:gaussianassumption} is that $h \geq N^{-1/6+\eps}$, i.e., the Gaussian fluctuations hold down to just above the current intermediate regime \eqref{eqn:mr-int-ass}. 

\paragraph{Free energy.} Our first results in this regime concern the free energy fluctuations, which are no longer Gaussian.  The following theorem provides an expansion of the free energy in terms of spectral quantities associated to the GOE; convergence of these quantities is stated below in Theorem \ref{thm:conv1}. The following is proven in Section \ref{sec:int-fe}. 
\bet[Free energy expansion] \label{thm:mr-int-fe-1} Assume that \eqref{eqn:mr-int-ass} holds.  For sufficiently small $\eps >0$ we have that the following estimate holds with probability at least $1-N^{-\eps}$ for $N$ large enough.
We have, 
\beq
N^{2/3} \left( F_{N, \beta, h}- \frac{1}{N} \log ( \Gamma (N/2)) + \frac{1}{N} ( N/2-1) \log (2 / (N \beta )) -\mathcal{C} \right) = \frac{1}{2} X_N  + \O (N^{-\eps} ),
\eeq
where 
\beq
\mathcal{C} := - \int \log(2-x) \rhosc (x) \d x,
\eeq
and 
$X_N$ is a random variable that is equal in distribution to a random variable $Y_N$ that satisfies the following.  With probability at least $1-N^{-\eps}$,
\beq
Y_N =  N^{2/3} ( \beta-1) ( x_b -2 ) - \frac{\theta}{N^{2/3}} \sum_{i=1}^N \frac{g_i^2}{ \mu_i - x_b}  + \O(N^{-\eps} )
\eeq
where the $\{ \mu_i\}_{i=1}^N$ are GOE eigenvalues, $\{g_i\}_i$ are independent standard normals and $x_b$ is the smallest solution larger than $\mu_1$ to the equation
\beq
( \beta-1) = \frac{ \theta}{N^{4/3}} \sum_i \frac{g_i^2}{ ( \mu_i -x_b)^2}.
\eeq
\eet
The random variable $X_N$ is constructed from the eigenvalues and eigenvectors of $H$, the GOE matrix associated to $M$ via $H= M+V$ where $V$ is a diagonal matrix of independent Gaussians with variance $2/N$ (see Section \ref{sec:notation}).  The proof of the above theorem explicitly constructs $X_N$ (it is of course very similar to $Y_N$ but not stated here for concerns of brevity). 

The second result concerns the convergence of the leading order contribution to the free energy.   For this statement we let $\{ g_i \}_{i=1}^\infty$ be an infinite sequence of iid standard normal random variables, and $\{ \chi_i \}_{i=1}^\infty$ be the Airy$_1$ random point field whose first $n$ particles are the limits of the top $n$ rescaled eigenvalues of the GOE, $\{ N^{2/3} ( \mu_i - 2) \}_{i=1}^n$.  We note that the particles of the Airy$_1$ random point field are almost surely distinct (see, e.g., Proposition 3.5 of \cite{rrv}).

We let $a>0$ be the unique solution in the positive half-line to 
\beq \label{eqn:mr-a-def}
\beta -1 = \theta \sum_{i=1}^\infty \frac{ g_i^2}{ ( \chi_i - \chi_1 - a)^2 }
\eeq
and $\xi$ be the random variable,
\beq \label{eqn:mr-xi-def}
\xi := \lim_{n \to \infty} \left(  (\beta-1)(\chi_1 + a) - \theta \left(   \sum_{i=1}^n  \frac{ g_i^2}{ \chi_i - \chi_1 -a }  + \frac{1}{ \pi} \int_{0}^{ \left( 3 \pi n /2  \right)^{2/3}}  \frac{1}{ \sqrt{x}} \d x \right) \right).
\eeq
The following theorem is proven in Section \ref{sec:int-conv}.  In this section we also prove that $a$ and  $\xi$ are well-defined and almost surely finite. 

\bet[Free energy fluctuations] \label{thm:conv1}
Let $\{ \mu_i \}_{i=1}^N$ be the eigenvalues of a GOE matrix.  Let $\{g_i \}_{i=1}^N$ be iid Gaussian random variables independent from $\{ \mu_i \}_{i=1}^N$.  Let $\xb$ be the unique solution to the following equation satisfying $\xb > \mu_1$,
\beq
\beta -1 = \frac{\theta}{N^{4/3}} \sum_{i=1}^N \frac {g_i^2}{ ( \mu_i - \xb )^2 }.
\eeq
Then $\xb$ and the random variable,
\beq
N^{2/3} (\beta-1) ( \xb- 2 ) + \frac{\theta}{N^{2/3} } \sum_{i=1}^N \left( \frac{ g_i^2}{ ( \mu_i - \xb) } - \frac{1}{ ( \gamma_i - 2 ) } \right)
\eeq
converge in distribution to the corresponding quantities defined in \eqref{eqn:mr-a-def} and \eqref{eqn:mr-xi-def} for the Airy$_1$ random point field.  As a consequence, the  free energy $F_{N, \beta, h}$, after an appropriate rescaling, converges to $\xi$ in the regime \eqref{eqn:mr-int-ass}.
\eet

%
%
%

\paragraph{Overlap with external field.} For the overlap with the external field, we have the following theorem, which is proven in Section \ref{sec:int-ext}. It exhibits quenched Gaussian fluctuations wrt the fixed Gibbs measure.
\bet[Quenched external field overlap fluctuations] \label{thm:mr-int-ext}
For $|t| \leq C$ we have with probability at least $1 - N^{-\eps}$ and $\eps >0$ small enough,
\beq
\log \langle \exp \left[ \beta^{1/2} t N^{-1/2} v^T \sigma \right] \rangle = \frac{t^2}{2} - \theta^{1/2} (N^{1/3} t) ( N^{-1} v^T(M- \gamma)^{-1} v ) + \O(N^{-\eps} )
\eeq
some $\eps>0$.
\eet
For the quenched expectation we have the following, proven at the end of Section \ref{sec:int-ext}.
\bet[Quenched external field overlap expectation] \label{thm:mr-int-ext-exp}
There is an $\eps >0$ so that with probability at least $1 - N^{-\eps}$ we have,
\beq \label{eqn:exp-exp-yy1}
\langle N^{-1/2} v^T \sigma \rangle = - ( \theta / \beta )^{1/2} N^{-2/3} v^T (M - \gamma )^{-1} v + \O ( N^{-\eps} )
\eeq
and that the random variable
\beq
N^{1/3} \left( \frac{1}{N} v^T (M- \gamma)^{-1} v  - \msc(2) \right)
\eeq
converges in distribution to the random variable,
\beq
\lim_{n \to \infty} \left(  \sum_{i=1}^n \frac{g_i^2}{ \chi_i - a} + \frac{1}{ \pi} \int_{0}^{ (3 \pi n /2 )^{2/3} } \frac{ 1}{ \sqrt{x}} \d x \right)
\eeq
where $a$ is as in \eqref{eqn:mr-a-def}, and $\{ \chi_i \}_{i=1}^\infty$ is the Airy$_1$ point field as defined above.
\eet

\paragraph{Overlap between two replicas.} For the overlap between two replicas we have the following.  It is proven in Section \ref{sec:int-over}. 
\bet[Quenched replica overlap fluctuations] \label{thm:mr-int-over}
There is a small $c_1 >0$ so that the following holds.  For every $\alpha >0$ and positive integer $k$, we have for all sufficiently small $\eps >0$ that there is an $\eps_1 >0$ such that the following holds on an event of probability at least $1 - N^{-\eps_1}$.  We have for all $|t| \leq N^{-\alpha}$
\begin{align}
\log \left\langle \exp \left( t N^{-2/3} \sigma^{(1)} \cdot \sigma^{(2)} \right) \right\rangle &= \sum_{j=1}^k t^j Z_j + \O \left( N^{-c_1} + |t|^{k+1} N^{\eps} \right)
\end{align}
where the $Z_j$ are specific functions of the spectral quantities of $M$ (described below).  On the above event we have $|Z_j| \leq N^{\eps}$.  Additionally,  the $Z_j$ jointly converge in distribution to functions of the Airy$_1$ random point field.
\eet
We were not able to find a simple form for the $Z_j$; they come from a certain Taylor expansion and involve general sums of products of the Stieltjes transforms of $M$ and the quadratic form $v^T (M-\gamma)^{-k} v$ evaluated at the point $\gamma$ defined above which comes from the steepest descent analysis.  However, we are able to show that they converge in distribution to similar quantities of the Airy$_1$ random point field.  The latter quantities are constructed by replacing the eigenvalues of $M$ by the particles of the Airy$_1$ random point field and the inner products $v \cdot u_i (M)$ by iid standard normal random variables.  Note that the expansion and estimate given above are strong enough to interpret the $Z_j$ as the leading order contribution to the quenched moments of the overlap between two replicas (which then converge in distribution).

\subsection{Low temperature and microscopic magnetic field} \label{sec:mr-micro}
The third and final scaling regime we consider is when the temperature is low, $\beta >1$ and the magnetic field is of order $h \sim N^{-1/2}$.  We assume,
\beq \label{eqn:mr-micro-ass}
1 + c \leq \beta \leq c^{-1}, \qquad h^2 \beta = \theta N^{-1},
\eeq
for fixed $c>0$ and $\theta >0$.  The results in this section are proven in Section \ref{sec:micro}.

\paragraph{Free energy.} Our first results show that for the free energy, the fluctuations coincide with the low temperature result of \cite{baiklee}  in that they are governed by the largest eigenvalue of $M$ and are asymptotically Tracy-Widom.  This theorem is proven in Section \ref{sec:micro-fe}. 
\bet[Free energy fluctuations] \label{thm:mr-micro-fe}
Under the assumption \eqref{eqn:mr-micro-ass} we have for all sufficiently small $\eps$ that with probability at least $1- N^{-\eps_1}$ for some $\eps_1 >0$ and $N$ large enough that,
\begin{align}
\frac{1}{N} \log Z_{N, \beta, h} &= \frac{\beta-1}{2} ( \lambda_1(M) -2) + C^{(m)}_N+ \O(N^{-1+\eps} ).
\end{align}
where,
\beq
C^{(m)}_N :=2 \beta - \int \log(2-x) \rhosc (x) \d x  +\frac{1}{N} \left( (1-N/2) \log (\beta) +\frac{N}{2} \log(2 \pi) + \frac{1}{2} \log(N) \right).
\eeq
Hence, $N^{2/3} 2(\beta-1)^{-1} \left( F_{N, \beta, h} - C^{(m)}_N \right)$ converges to a TW$_1$ random variable.
\eet

\paragraph{Overlap between two replicas.} The following results concern the overlap between two replicas.  They all are proven in Section \ref{sec:micro-over}.   The first concerns the fluctuations of the weights in the Parisi measure.  Recall the notation $v_1 := v^T u_1 (M).$
\bet[Overlaps measure] \label{thm:mr-micro-over-1}
Suppose that \eqref{eqn:mr-micro-ass} holds.    For all sufficiently small $\eps >0$ there is an $\eps_1>0$ so that with probability at least $1-N^{-\eps_1}$  and all $1>t>0$,
\beq
\left\langle \1_{ \{  | N^{-1}  \sigma^{(1)} \cdot \sigma^{(2)} \mp  (1-\beta^{-1} ) |  \leq t \} } \right\rangle = \frac{1}{2} \pm \frac{1}{2} \tanh^2 \left( \sqrt{ v_1^2 \theta(\beta-1)} \right) + N^{\eps} \O \left( t+ N^{-2/3+\eps} t^{-2} + N^{-1/3} \right).
\eeq
We have also that $v_1$ converges to a standard normal random variable as $N \to \infty$, and so we get a convergence in distribution result for the random variable on the LHS for any $t$ satisfying $N^{-1/3+3\eps} \leq t \leq N^{-2\eps}$.
\eet

Recall the definition of $\tilm (z)$ in \eqref{eqn:m-tilde-def-f}. 
We will often consider $\tilm ( \lambda_1(M))$ which we will later see satisfies
\beq
\left| \tilm( \lambda_1 (M) ) +1 \right| \leq N^{-1/3+\eps}
\eeq
with probability at least $1- N^{-\eps/10}$ for any sufficiently small $\eps >0$ and large enough $N$.  We have the following result for the quenched expectation of the overlap.

\bet[Quenched overlaps expectation] \label{thm:mr-micro-over-2}
For all sufficiently small $\eps >0$ there is an $\eps_1>0$ so that with probability at least $1-N^{-\eps_1}$,
\beq \label{eqn:mr-micro-1}
\frac{1}{N} \langle  \sigma^{(1)} \cdot \sigma^{(2)} \rangle = \frac{ \beta + \tilm ( \lambda_1 (M) ) }{ \beta} \left( \tanh \left( \sqrt{ v_1^2 \theta ( \beta  + \tilm (\lambda_1 (M) ) ) } \right) \right)^2 + \O(N^{-2/3+\eps} ).
\eeq
In particular, $N^{-1} \langle \sigma^{(1)} \cdot \sigma^{(2)} \rangle $ converges in distribution to the random variable $(1-\beta^{-1} )\tanh^2 \sqrt{Z^2 \theta(\beta-1) } $ where $Z$ is a standard normal random variable.  
\eet
In the work \cite{SSK-ls} we considered the random variable,
\beq
\Xi_N := N^{1/3} (\tilm ( \lambda_1 (M) )+1)
\eeq
and showed that it has a limit $\Xi$ given by,
\beq
\Xi := \lim_{n \to \infty} \left( \sum_{i=2}^n \frac{1}{ \chi_i - \chi_1}  + \int_0^{ \left( \frac{ 3 \pi n}{2}\right)^{2/3} } \frac{1}{ \pi \sqrt{x} } \d x \right)
\eeq
where $\{ \chi_{i} \}_{i=1}^\infty$ are the particles of the Airy$_1$ random point field.  

The random variable $\Xi_N$ appears in Theorem \ref{thm:mr-micro-over-2} and it would be of interest to obtain a kind of second order fluctuation or conditional result about the quenched expectation of the overlap.

If we Taylor expand the leading order term in \eqref{eqn:mr-micro-1} in $1 + \tilm ( \lambda_1 (M))$ we obtain a convergence statement for a randomly rescaled version of $N^{-1}\langle  \sigma^{(1)} \cdot \sigma^{(2)}\rangle$ to $\Xi$ after subtracting and dividing by explicit functions of $v_1^2$.  However, in our view this is not quite satisfactory as it falls short of a statement of the fluctuations of the overlap \emph{conditional} on $v_1^2$.  This is due to the fact that, for the SSK as defined above, the eigenvector projection $v_1$ is not stochastically independent from the eigenvalues of $M$.  

However, we are able to obtain a finer result under either of two possible modifications to the model.  

The first modification is if we replace the SSK Hamiltonian by a simpler model,
\beq \label{eqn:SSK-simp}
\tilde{H}_N ( \sigma) = \frac{1}{ \sqrt{2N} } \sum_{i, j} g_{ij} \sigma_i \sigma_j + h v^T \sigma.
\eeq
Denote now by $H$ the matrix given by $H_{ij} = (2N)^{-1} (g_{ij} + g_{ji})$.  In this case, the same estimate as \eqref{eqn:mr-micro-1} holds with the spectral quantities of $M$ now replaced by those of $H$.  In the case of $H$, the eigenvectors and eigenvalues are independent and so we can obtain a statement of the conditional fluctuations of the quenched overlap.  The second modification is to simply let $v$ be a vector uniformly distributed on the sphere. In either case we obtain the following.
\bet \label{thm:mr-micro-over-3}
Assume either that the SSK model is replaced by the simpler GOE-associated model \eqref{eqn:SSK-simp}, or that the vector $v$ is uniformly distributed on the $N-1$ sphere of radius $\sqrt{N}$.  Let $\eps >0$ be sufficiently small, and $F$ a Lipschitz function.  Denote by $z$ the random variable $v^T \cdot u_1(H)$ in the case that $M$ is replaced by $H$ or $v^T u_1(M)$ in the case that it is uniformly distributed.  There is an event of probability at least $1-N^{-\eps}$ so that the following holds.
\beq
\ee\left[ F ( N^{1/3} b_z^{-1} \left( N^{-1} \langle \sigma^{(1)} \cdot \sigma^{(2)}\rangle-a_z \right) )\vert z \right] = \ee[ F( \Xi_N ) ] + \O (N^{-\eps} ),
\eeq
where,
\begin{align}
a_z&:= \frac{ \beta-1}{ \beta} \tanh^2(Z), \qquad b_z:= \tanh(Z) \frac{ \sinh(Z) \cosh(Z) + Z }{ \beta \cosh^2 (Z) }, \nonumber\\
Z &:= \sqrt{z^2 \theta ( \beta-1 ) }.
\end{align}
\eet
We have a final result on the fluctuations of the quenched variance.
\bet[Quenched overlap variance] \label{thm:mr-micro-over-4}
For all sufficiently small $\eps >0$ there is an $\eps_1 >0$ so that we have the following estimate with probability at least $ 1- N^{-\eps_1}$,
\begin{align}
&\frac{1}{N^2} \left\langle \left( \sigma^{(1)} \cdot \sigma^{(2)} \right)^2 \right \rangle - \left( \frac{1}{N} \left\langle \sigma^{(1)} \cdot \sigma^{(2)} \right\rangle \right)^2 \nonumber \\
=& \frac{( \beta + \tilm ( \lambda_1 (M) ))^2 }{\beta^2} (1 - \tanh^4 \left( \sqrt{ v_1^2 \theta ( \beta + \tilm ( \lambda_1 (M) ) ) } \right) ) + \O(N^{-2/3+\eps} ),
\end{align}
and so the random variable on the LHS converges in distribution to the random variable $(1-\beta^{-1})^2 (1-\tanh^4 ( \sqrt{ W^2 \theta(\beta-1) } ) )$ where $W$ is a standard normal random variable.  
\eet

\section{Results from random matrix theory} \label{sec:rmt}
In this section we collect some results from random matrix theory.  Recall the matrix $M$ defined by \eqref{eqn:M-def-f}. We have the following local semicircle law and rigidity result, quantifying the convergence of the eigenvalues of $M$ to the semicircle distribution.  The statements and proofs can be found in Theorem 2.6 and Theorem 10.3 of \cite{local-lectures}. Throughout this section we denote $z  = E + \i \eta$ for $E \in \rr$ and $\eta > 0$. 
\bet \label{thm:locallaw}
Let $M$ be as above.  Let $\omega >0$.  Consider the spectral domain,
\beq
\S_1 := \{ z  = E + \i \eta : |E| \leq \omega^{-1}, N^{-1+\omega} \leq \eta \leq \omega^{-1} \}.
\eeq
For any $\eps >0$ and $D>0$, the following estimate holds for $N$ large enough.  
\beq
\pp\left[ \bigcap_{z \in \S_1} \left\{ |m(z) - \msc (z) | \leq \frac{N^{\eps}}{N \eta } \right\} \right] \geq 1 - N^{-D}.
\eeq
Consider the spectral domain,
\beq
\S_2 :=\{ z = E + \i \eta :2+ N^{-2/3+ \omega} \leq |E| \leq \omega^{-1}, 0 \leq \eta \leq \omega^{-1} \}.
\eeq
For any $\eps >0$ and $D>0$, the following estimate holds for $N$ large enough.   
\beq
\pp \left[ \bigcap_{ z \in \S_2 } \left\{ | m(z) - \msc (z) | \leq \frac{N^{\eps}}{N ( | |E| - 2 | + \eta )} \right\} \right] \geq 1 - N^{-D}.
\eeq
The same estimates hold for $H$, where $H$ is defined as in Section \ref{sec:notation}. 
\eet
We recall the definition of the quantiles $\gamma_i$ in \eqref{eqn:gamma-i-def}. We have the following rigidity theorem, a statement of which can be found in Theorem 2.9 of \cite{local-lectures}.

\bet \label{thm:rigi}
For any $\eps >0$ and $D>0$ we have for $N$ large enough,
\beq
\pp \left[ \bigcap_{ j } \left\{ | \lambda_j(M) - \gamma_j | \leq N^{\eps} N^{-2/3} \min\{ j^{-1/3}, (N+1-j)^{-1/3} \} \right\} \right] \geq 1 - N^{-D}.
\eeq 
The same estimate holds for $H$ as defined in Section \ref{sec:notation}. The quantiles satisfy,
\beq
2 - \gamma_i \asymp \frac{i^{2/3}}{N^{2/3}} .
\eeq
\eet
We will also need the so-called ``isotropic'' local law and delocalization results for the matrix $M$. To state them, we introduce the resolvent of $M$ by
\beq
R(z) = \frac{1}{ M - z}.
\eeq
The following statements are from Theorems 2.12 and 2.15 of \cite{alex-iso}. 
\bet \label{thm:iso}
Let $\omega >0$ and $\S_1$ and $\S_2$ be as above.  Let $\eps >0$ and $D>0$.  Let $S^{N-1}$ denote the unit sphere in $\rr^N$.  For $N$ large enough, we have the following estimates,
\beq
\inf_{u, w \in S^{N-1} } \pp \left[ \bigcap_{z \in \S_1} \left\{ \left| u^T R (z) w - (u^T w) \msc (z) \right| \leq N^{\eps} \left( \sqrt{ \frac{ \Im [ \msc (z) ] }{ N \eta } } + \frac{1}{ N \eta } \right) \right\} \right] \geq 1 - N^{-D}.
\eeq
and
\beq
\inf_{u, w \in S^{N-1} } \pp \left[ \bigcap_{z \in \S_2} \left\{ \left| u^T R (z) w - (u^T w) \msc (z) \right| \leq N^{\eps} \sqrt{ \frac{ \Im [ \msc (z) ] }{ N \eta } }  \right\} \right] \geq 1 - N^{-D}.
\eeq
If $u_i$ are the eigenvectors of $M$, then for any $\eps >0$ and $D>0$ we have for $N$ large enough,
\beq
\inf_{w \in S^{N-1} } \pp\left[ \bigcap_{j} \{ N^{1/2} |w^T u_i| \leq N^{ \eps} \} \right] \geq 1 - N^{-D}.
\eeq
The same estimate holds for $H$ as defined in Section \ref{sec:notation}. 
\eet
We note the following behavior for $\Im [ \msc (E+ \i \eta)]$ as it appears in the estimates above (see Lemma 3.3 of \cite{local-lectures}),
\begin{align}
\Im[ \msc ( E + \i \eta ) ] \asymp \begin{cases} \sqrt{ | |E| -2 | + \eta } & |E| < 2 \\ \frac{ \eta}{ \sqrt{ | |E| -2 |  + \eta } }& |E| \geq 2 \end{cases}.
\end{align}

Due to the Cauchy integral formula, Theorems \ref{thm:locallaw}, \ref{thm:rigi} and \ref{thm:iso} imply the following.
\bet \label{thm:ll-derivs}
With $\S_2$ as above we have for any $k$ and $\eps >0$ that  for $N$ large enough, 
\beq
\pp \left[ \bigcap_{ z \in \S_2 } \left\{ | m^{(k)}(z) - \msc^{(k)} (z) | \leq \frac{N^{\eps}}{N ( | |E| - 2 | + \eta )^{k+1}} \right\} \right] \geq 1 - N^{-D},
\eeq
and
\beq
\inf_{u, w \in S^{N-1} } \pp \left[ \bigcap_{z \in \S_2} \left\{ \left| u^T \del^k_z R (z) w - (u^T w) \msc^{(k)} (z) \right| \leq N^{\eps}\frac{1}{( | |E| - 2 | + \eta )^{k}} \sqrt{ \frac{ \Im \msc (z) }{ N \eta } }  \right\} \right] \geq 1 - N^{-D}.
\eeq
\eet

The following level repulsion estimate is proven in Lemma 3.4 of \cite{SSK-ls}. 
\bet \label{thm:lr} Let $\eps >0$.  For $N$ large enough it holds for all $N^{-1/3} < s < 1$ that
\beq
\pp[ | \lambda_1 (M) - \lambda_2 (M) | \leq s N^{-2/3}] \leq N^{\eps} s, \qquad \pp [ | \lambda_1 (H) - \lambda_2 (H) | \leq s N^{-2/3} ] \leq N^{\eps} s.
\eeq
\eet

\subsection{Comparison of $M$ to the GOE}

Recall that $H$ was introduced in Section \ref{sec:notation} via $H = M + V$ where $V$ is a diagonal matrix whose entries are iid Gaussians with variance $2 /N$.  Then $H$ is a matrix from the GOE.  The first estimate in the following was proven in Proposition A.1 of \cite{SSK-ls}.  The second estimate is proven in Appendix \ref{a:compare}.
\bet \label{thm:compare}
Let $\eps >0$ and $D>0$.  For $N$ large enough it holds that,
\beq \label{eqn:compare1}
\pp \left[ \bigcap_{ 1 \leq j \leq N^{1/20} } \left\{ | \lambda_j (M) - \lambda_j (H) | \leq \frac{ N^{\eps}}{N} \right\} \right] \geq 1- N^{-D}.
\eeq
Let $u_1 (M)$ and $u_1 (H)$ be the first eigenvectors of $M$ and $H$.  There is a $c >0$ so that the following estimate holds.
\beq \label{eqn:compare2}
\sup_{v \in S^{N-1} } \pp \left[ | (v^T u_1 (M) )^2 - (v^T u_1 (H) )^2 | \geq N^{-1-c} \right] \leq N^{-c}.
\eeq
\eet
We also have the following lemma, proven in Appendix \ref{a:compare}. 
\bel \label{lem:HViso}
Let $ \delta >0$ and $\omega >0$.   For any $\eps$ and $D>0$,  the following estimates hold for $N$ large enough.  First,
\begin{align} 
\inf_{ v \in \S^{N-1} } &  \pp\left[  \bigcap_{ z \in \S_1 \cap \{ \eta  \leq N^{-\delta} \} } \left\{ \left| v^T (M - z )^{-1} v - v^T (H-z )^{-1} v \right| \leq N^{\eps} \left( \frac{1}{ \sqrt{N}} + \frac{1}{ (N \eta)^2}  + \frac{ \Im [ \msc (z) ]}{ N \eta } \right) \right\} \right] \notag \\
\geq& 1 - N^{-D}. \label{eqn:isocompare1}
\end{align}
Second,
\beq
\inf_{ v \in \S^{N-1} } \pp\left[  \bigcap_{ z \in \S_2 \cap \{ \eta  \leq N^{-\delta} \} } \left\{ \left| v^T (M - z )^{-1} v - v^T (H-z )^{-1} v \right| \leq N^{\eps} \left( \frac{1}{ \sqrt{N}}  + \frac{ \Im [ \msc (z) ]}{ N \eta } \right) \right\} \right] \geq 1 - N^{-D}. \label{eqn:compare3}
\eeq
\eel

\section{Representation formulas} \label{sec:rep}

This section contains  the various representation formulas for the free energy and various Gibbs expectations that we will use.  They reduce various high-dimensional integrals over the sphere $(N-1)$-sphere to  low dimensional contour integrals.  Such representations were first used in \cite{KTJ} and \cite{baiklee} to study the free energy.  These representations were extended in \cite{vl-ps,SSK-ls} to study the overlap in the model without magnetic field.  Here we extend these representations to allow for a magnetic field in the Hamiltonian.

Recall that $\d \omega_{N-1}$ is the uniform measure on the $N-1$ sphere of radius $\sqrt{N}$, which we denote by $\ss^{N-1}$.  As the proofs of the following results are similar as to previous results in the literature we defer their proofs to Appendix \ref{a:rep}.  
\bep \label{prop:rep}
Let $M$ be a real symmetric matrix, $v$ a vector and $\beta, \lambda >0$.  Then,
\begin{align} \label{eqn:representation1}
\int \exp \left[ \frac{ \beta}{2} \sigma^T M \sigma + \lambda v^T \sigma \right] \d \omega_{N-1} ( \sigma ) 
= \frac{ \beta N^{1/2}}{2 \pi \i }  \left( \frac{ 2 \pi}{\beta} \right)^{ \frac{N}{2} }  \int_{\gamma - \i \infty}^{ \gamma  +\i \infty} \exp \left[ \frac{N}{2} G_o (z, v, \lambda, \beta ) \right] \d z
\end{align}
where,
\beq
G_o(z, v, \lambda, \beta) := \beta z - \frac{1}{N} \sum_{i=1}^N \log (z  - \lambda_i (M) ) - \frac{ \lambda^2}{N \beta } v^T \frac{1}{M-z} v
\eeq
and $\gamma > \lambda_1 (M)$.  Furthermore,
\begin{align} 
&\int \int \exp \left[ \frac{\beta}{2} \sigma_1^T M \sigma_1 + \frac{ \beta}{2} \sigma_2^T M \sigma_2 + t \sigma_1^T \sigma_2 + \lambda v^T ( \sigma_1 + \sigma_2 ) \right] \d \omega_{N-1} ( \sigma_1 ) \d \omega_{N-1} ( \sigma_2 ) \nonumber \\
= &\frac{ \beta^2 N}{ ( 2 \pi \i)^2} ( 2 \pi)^N \int \int \exp\left[ \frac{N}{2} G_d (z, w, v, \lambda, t, \beta ) \right] \d z \d w, \label{eqn:representation2}
\end{align}
where
\begin{align}
G_d = \beta (z + w ) - \frac{1}{N} \sum_{i=1}^N \log (\beta^2 (z - \lambda_i )(w - \lambda_i ) - t^2 ) -\frac{ \lambda^2}{N } \sum_{i=1}^N (v^T u_i(M))^2\frac{ \beta(2\lambda_i - z - w) - 2 t }{ \beta^2 ( \lambda_i - w ) ( \lambda_i - z) - t^2 },
\end{align}
and the double contour integral is over the two lines $\gamma \pm \i \infty$ where $\gamma > \lambda_1 (M) +|t|/\beta$. 
\eep

\bep \label{prop:moment-rep}
Let $M$ be a real symmetric matrix, $v$ a vector, and $\beta, \lambda >0$.  Then,
\begin{align}
& \int ( \sigma^{(1)} \cdot \sigma^{(2)}) \exp \left[ \frac{\beta}{2} \left( (\sigma^{(1)})^T M \sigma^{(1)} +(\sigma^{(2)})^T M \sigma^{(2)}\right) + \lambda v^T ( \sigma^{(1)} + \sigma^{(2)}) \right] \d \sigma^{(1)} \d \sigma^{(2)}  \nonumber  \\
 = & \frac{ \beta^2N}{ (2 \pi \i)^2} \left( \frac{2 \pi}{ \beta} \right)^N  \int_{\Gamma^2}  \frac{ \lambda^2}{ \beta^2} v^T \frac{1}{ (M-z)(M-w)} v \exp\left[ \frac{N}{2} (G_o(z, v, \lambda, \beta) + G_o(w,v, \lambda, \beta) ) \right] \d z \d w \label{eqn:2repa}
\end{align}
where $\Gamma$ is a vertical line in the complex plane lying to the right of $\lambda_1 (M)$ and the function $G_o$ is as in the previous proposition.  Additionally,
\begin{align}
& \int ( \sigma^{(1)} \cdot \sigma^{(2)})^2 \exp \left[ \frac{\beta}{2} \left( (\sigma^{(1)})^T M \sigma^{(1)} +\sigma^{(2)})^T M \sigma^{(2)}\right) + \lambda v^T ( \sigma^{(1)} + \sigma^{(2)}) \right] \d \sigma^{(1)} \d \sigma^{(2)} \nonumber \\
= &  \frac{ \beta^2N}{ (2 \pi \i)^2} \left( \frac{2 \pi}{ \beta} \right)^N  \int_{\Gamma^2}  \d z \d w  \exp\left[ \frac{N}{2} (G_o(z, v, \lambda, \beta) + G_o(w,v, \lambda, \beta) ) \right] \nonumber \\
\times & \bigg\{  \frac{1}{ \beta^2} \sum_i \frac{1}{ ( \lambda_i - z ) ( \lambda_i  - w ) }  - \frac{ \lambda^2}{ \beta^3} v^T \left( \frac{1}{(M-z)^2(M-w) } + \frac{1}{ (M-w)^2 (M-z)} \right) v \nonumber \\
+& \frac{ \lambda^4}{ \beta^4} \left( v^T \frac{1}{ (M-w)(M-z) } v \right)^2 \bigg\}
\end{align}
\eep

\section{Gaussian regime} \label{sec:gauss}
In this section we will investigate the fluctuations of the SSK model in the regime of Gaussian fluctuations, i.e., in which either the temperature is high or the magnetic field is not decaying too quickly.  This section contains all the proofs of the results of Section \ref{sec:mr-gauss}.
\subsection{Free energy fluctuations} \label{sec:gauss-fe}
In this section we consider the free energy fluctuations.  
We define $\theta = h^2 \beta$ for notational simplicity.  We will assume that \eqref{eqn:gaussianassumption} holds for some $c>0$ and $\tau >0$. In applying this assumption it can be helpful to note that
\beq
(1 - \beta)_+  + \frac{ \theta}{ |1- \beta | + \sqrt{ \theta} } \asymp \begin{cases} 1-\beta + \sqrt{\theta} , & 1 \geq \beta \\ \frac{\theta}{\beta-1+\sqrt{\theta}} , & \beta \geq 1 \end{cases}.
\eeq
Our first task is to derive the following expansion of the free energy in terms of random matrix quantities.

\bep[Free energy expansion] \label{prop:gaussianexpansion}
Assume that \eqref{eqn:gaussianassumption} holds for some $\tau, c >0$. 
The following holds with overwhelming probability:
\begin{align}
F_{N, \beta, h} =&  \frac{G ( \gamma)}{2} + \frac{1}{N}\log \Gamma (N/2) +\frac{1}{N} (N/2-1) \log(2 / (N\beta) ) - \frac{1}{2N} \log(N G'' ( \gamma) \pi ) + \O (N^{-1-\tau/10}) \nonumber \\
&=   \frac{ G ( \hatgam)}{2} + \frac{1}{N} \log \Gamma (N/2) + \frac{1}{N} (N/2 - 1) \log (2 / (N \beta ) ) - \frac{1}{2N} \log (N g'' ( \hatgam) \pi ) \nonumber \\
&+ N^{-\tau/10}   \O \left( \frac{1}{N} + \frac{\theta}{ N^{1/2}  \kapgam ^{1/4} } \right)
\end{align}
where $G$ is as in \eqref{eqn:not-G}. 
The quantity $\kappa$ has the asymptotics,
\beq \label{eqn:kap-ass-f}
\kappa^{1/2} \asymp \frac{h^2 \beta}{ |1- \beta | + \sqrt{ h^2 \beta} } + (1- \beta)_+ .
\eeq
\eep

Throughout Section \ref{sec:gauss} we will abbreviate the eigenvalues of $M$ as
\beq
\lambda_i := \lambda_i (M).
\eeq
Recall they are arranged in decreasing order.  
From Theorem \ref{thm:ll-derivs} we have for any $\eps >0$ that for any $E > 2 + N^{-2/3+\eps}$,
\begin{align} \label{eqn:G-der-gauss-f}
G' (E) &= \beta + m (E) - \theta N^{-1} v^T (M -z)^{-2} v \nonumber\\
 =& \beta + \msc (E) - \theta \msc' (E) + N^{\eps} \O \left( \frac{1}{ N (E-2)} + \frac{ \theta}{ N^{1/2} (E-2)^{5/4} } \right)
\end{align}
with overwhelming probability. 
   Since $G' (E)$ is monotonic, is positive for large $E$ and tends to $-\infty $ as $E$ tends to $\lambda_1$ from the right, we see that the equation
$
G' ( \gamma ) = 0
$
has a unique solution for $E > \lambda_1$.  We first seek to establish preliminary estimates on the location of the saddle $\gamma$.   
Recall our definition of the function $g(z)$ in Section \ref{sec:notation} so that,
\beq
g' (z) = \beta + \msc ( z) - \theta \msc' ( z). 
\eeq
Note $g(z)$ is an approximation to $G(z)$.  Recall that $\hatgam$, the deterministic saddle, is defined by $g'(\hatgam) = 0$. 
Using the fact that,
\beq
1 + \msc ( \hatgam) \asymp \sqrt{ \hatgam - 2 } , \qquad \msc' ( \hatgam) \asymp ( \hatgam - 2 )^{-1/2},
\eeq
a routine calculation using the quadratic formula shows that
\beq
\sqrt{ \hatgam - 2 } \asymp ( 1 - \beta)_+  + \frac{ \theta}{ |1- \beta | + \sqrt{ \theta} } \asymp \begin{cases} (1 - \beta) + \sqrt{ \theta} , \quad & 1 \geq \beta \\ \frac{ \theta}{ \beta-1 + \sqrt{ \theta} } , \quad & \beta \geq 1 \end{cases}
\eeq
proving \eqref{eqn:kap-ass-f}. 
Hence, under the assumption \eqref{eqn:gaussianassumption} we have that $
\kappa \geq N^{-2/3+ \tau}
$
for $N$ large enough.  
\bel  \label{lem:sad-est}
Assume that \eqref{eqn:gaussianassumption} holds.  Then, for $N$ large enough, $\hatgam -2 \geq N^{-2/3+\tau}$ and with overwhelming probability, the random saddle $\gamma$ satisfies $\gamma - 2 \geq c N^{-2/3+\tau}$ for some $c>0$ and
\beq
| \hatgam - \gamma | \leq N^{\eps'} \frac{1}{1 +  \theta (\hatgam-2)^{-1} } \left( \frac{1}{ N \sqrt{  \hatgam - 2 } } + \frac{ \theta} { N^{1/2} ( \hatgam -2 )^{3/4} } \right) \leq N^{\eps'} N^{-2/3-\tau/10},
\eeq
for any $\eps' >0$. 
\eel
\proof Observe for $E-2 \geq N^{-2/3+\eps}$, we have by Theorem \ref{thm:ll-derivs} that with overwhelming probability,
\begin{align} \label{eqn:zz-1}
G''(E) &= m' (E) -2 \theta N^{-1} v^T (M -z )^{-3} v \nonumber \\
=& \msc' (E) - 2 \theta \msc''(E) + N^{\eps} \O \left( \frac{1}{ N (E-2)^2} + \frac{ \theta}{ N^{1/2} (E-2)^{9/4} }  \right) \nonumber \\
 \asymp & \frac{1}{ \sqrt{ E-2}} + \frac{ \theta}{ (E-2)^{3/2}} \asymp g''(E)
\end{align}
Note that the condition that $E-2 \geq N^{-2/3+\eps}$ ensures that the error term in the second to last expression is smaller than the final line appearing above, justifying the inequalities implicit in the last line (our definition of $\asymp$ appears in Section \ref{sec:notation}).
For $E-2 \geq N^{-2/3+\eps}$ we have, using \eqref{eqn:G-der-gauss-f} and Theorem \ref{thm:ll-derivs},
\beq \label{eqn:aa-1-f}
 \frac{ |g' (E) - G' (E) | }{ | g'' (E) | } \leq N^{\eps'}  \frac{1}{ 1 + \theta / (E-2)} \left( \frac{1}{N \sqrt{ E-2 }} + \frac{\theta}{ N^{1/2}(E-2)^{3/4} } \right) \leq C (E-2) \frac{N^{\eps'}}{N^{1/2} (E-2)^{3/4} }
\eeq
with overwhelming probability for any $\eps, \eps' >0$.   Recall that $\hatgam -2 > N^{-2/3+\tau}$, and $g'(E) $ and $G' (E)$ are monotonic. Using the second inequality in \eqref{eqn:aa-1-f} for the quantity on the LHS together with \eqref{eqn:zz-1} we conclude first that $\gamma \geq c N^{-2/3+\tau}$ for some $c>0$ with overwhelming probability, by using the definitions $g'(\hatgam) = G'(\gamma) = 0$. Then, using the first inequality in \eqref{eqn:aa-1-f}  we can conclude the estimate of the lemma by Taylor expanding the LHS of the equation $g'( \hatgam) - g' ( \gamma) = G'(\gamma) - g'(\gamma)$. \qed

With these preparations, we implement the method of steepest descent to prove Proposition \ref{prop:gaussianexpansion} in the next subsection.

\subsubsection{Proof of Proposition \ref{prop:gaussianexpansion}} 
Using Proposition \ref{prop:rep} we are led to investigate the contour integral,
\beq
\int_{\gamma - \i \infty}^{ \gamma + \infty} \exp \left[ \frac{N}{2} (G(z) - G ( \gamma) ) \right] \d z.
\eeq
We decompose the contour into two components, $\Gamma_1 \cup \Gamma_2$ where,
\beq
\Gamma_1 := \{ \gamma \pm \i t, |t| \leq N^{\eps} \}, \qquad \Gamma_2 := \{ \gamma \pm \i t, |t| > N^{\eps} \},
\eeq  for some small $\eps >0$.  The contribution from $\Gamma_2$ is exponentially small, and is handled via calculations almost identical to those appearing in Lemma 5.4 of \cite{baiklee}, and so we omit this part of the proof.  

Via Taylor expansion and the estimates of Theorem \ref{thm:ll-derivs}, we have for any $|t| \leq 10$, with overwhelming probability,
\begin{align}\label{eqn: taylor-infant}
G ( \gamma+ \i t ) - G ( \gamma) = - \frac{t^2}{2} G'' ( \gamma) +  \O \left( |t^3|  \kapgam^{-3/2}(1+\theta  \kapgam^{-1} )\right).
\end{align}
Fixing $\alpha >0$ sufficiently small, we will use this expansion when
\beq
|t| \leq t_\alpha := \frac{N^{\alpha}}{N^{1/2} }|G''( \gamma ) |^{-1/2} \asymp \frac{N^{\alpha}}{N^{1/2}}  \kapgam^{1/4} (1 + \theta /  \kapgam )^{-1/2},
\eeq
where the final inequalities hold with overwhelming probability using \eqref{eqn:zz-1}. 
Then, for $|t| \leq t_\alpha$, the error term in \eqref{eqn: taylor-infant} is of order
\beq
|t^3| \kappa^{-3/2} ( 1 + \theta \kapgam^{-1} ) \leq \frac{C}{N} \frac{N^{3 \alpha}}{N^{1/2} \kapgam^{3/4} } \frac{1}{ (1 + \theta /\kapgam)^{1/2} }  \leq N^{-1-\tau/10}
\eeq 
for $\alpha < \tau/100$ and $N$ large enough.  Hence,
\begin{align}
\int_{\gamma - \i t_\alpha}^{\gamma + \i t_\alpha }  \exp \left[ \frac{N}{2} (G(z) - G ( \gamma) ) \right] \d z &= \int_{- t_\alpha}^{t_\alpha } \exp \left[ - \frac{N}{4} G''( \gamma ) t^2 \right] \left( 1 + \O (N^{-\tau/10} ) \right) \i \d t \nonumber \\
&=\i \sqrt{ \frac{4\pi}{N G'' ( \gamma) }} \left(1 + \O (N^{-\tau/10} ) \right)
\end{align}
Using also that $\del_{ \eta} \Re [ G ( E + \i \eta ) ] = - \Im [ G' (E + \i \eta ) ] <0$ for $E > \lambda_1$ and the fact that 
\beq
\Re[ G (E + \i t_\alpha ) - G( \gamma ) ] \asymp -t_\alpha^2 G'' ( \gamma) \asymp - N^{ \alpha-1}
\eeq
we see that the contribution of the contour integral from $t_\alpha < |t| < N^{\eps}$ is  $\O (\e^{ -c N^{\alpha}})$.  Summarizing, this proves that
\beq
|\ss^{N-1} | Z_{N, \beta, h} = \frac{ \beta N^{1/2}}{2 \pi} \left( \frac{ 2 \pi }{ \beta} \right)^{ \frac{N}{2} } \exp \left[ \frac{N}{2} G ( \gamma ) \right]  \sqrt{ \frac{4\pi}{N G'' ( \gamma) }} \left(1 + \O (N^{-\tau/10} ) \right),
\eeq
with overwhelming probability.  Hence, with overwhelming probability,
\beq
\frac{1}{N} \log Z_{N, \beta, h} = \frac{G ( \gamma)}{2} + \frac{1}{N}\log \Gamma (N/2) +\frac{1}{N} (N/2-1) \log(2 / (N\beta) ) - \frac{1}{2N} \log(N G'' ( \gamma) \pi ) + \O (N^{-1-\tau/10}),
\eeq
which proves the first claim of the proposition.  
Now, with overwhelming probability, for any $\eps >0$,
\begin{align}
G ( \gamma) - G ( \hatgam ) = \O ( G'' ( \gamma ) ( \gamma - \hatgam )^2 ) &= \frac{N^{\eps}}{ \kapgam^{1/2}(1+ \theta / \kapgam) }  \O\left( \frac{1}{ N^2 \kapgam } + \frac{\theta^2}{ N  \kapgam^{3/2} } \right) = N^{-\tau/10}  \O \left( \frac{1}{N} + \frac{\theta}{ N^{1/2} \kapgam^{1/4}} \right),
\end{align}
and 
\beq
\log(G'' ( \gamma) ) = \log (g'' ( \hatgam ) ) + \O (N^{-\tau/10} ).
\eeq
Above, we used \eqref{eqn:zz-1} and Lemma \ref{lem:sad-est}. 
Hence, with overwhelming probability,
\begin{align}
\frac{1}{N} \log Z_{N, \beta, h} &= \frac{ G ( \hatgam)}{2} + \frac{1}{N} \log \Gamma (N/2) + \frac{1}{N} (N/2 - 1) \log (2 / (N \beta ) ) - \frac{1}{2N} \log (N g'' ( \hatgam) \pi ) \nonumber \\
&+ N^{-\tau/10}   \O \left( \frac{1}{N} + \frac{\theta}{ N^{1/2} \kapgam^{1/4} } \right),
\end{align}
which is the final claim of the proposition. \qed

\subsubsection{Proof of Theorem \ref{thm:mr-gauss-fe}}

We start by considering the second parameter regime of the theorem statement, where the magnetic field is of at least moderate size. Under the assumptions of the theorem statement, it is straightforward to see that \eqref{eqn:gaussianassumption} holds and moreover that
\beq \label{eqn:order-compare-f}
\frac{N^{\eps}}{N} \leq \frac{ \theta}{ N^{1/2} \kappa^{1/4} }
\eeq
for some $\eps >0$ using \eqref{eqn:kap-ass-f}. Recall now that $G( \hatgam )$ has two fluctuating components, the linear statistic $N^{-1} \sum_i \log ( \hatgam - \lambda_i )$ and the resolvent component $\frac{\theta}{N} v^T (M- \hatgam)^{-1} v$. By the rigidity estimates of Theorem \ref{thm:rigi} it is straightforward to check that the former has fluctuations of size at most $N^{\eps/2-1}$. By \eqref{eqn:order-compare-f} these are of lower order than $\frac{\theta}{N^{1/2} \kappa^{1/4} }$ which by Theorem \ref{thm:iso-clt} are the size of the fluctuations of the resolvent component. The statement now follows from a direct application of Theorem \ref{thm:iso-clt}.

We now consider the first regime of the theorem statement. In this regime we have that $\kappa \asymp 1$. For the linear statistic component of $G( \hatgam)$ we have, see, e.g., Section 4 of \cite{landon2020applications} that the quantity
\beq
\frac{1}{\hat{B}_N^{1/2}} \left( \sum_i \log ( \hatgam - \lambda_i ) - N \int \log ( \hatgam - x) \rhosc(x) \d x - \hat{A}_N \right)
\eeq
converges to a standard normal random variable, where
\begin{align}
\hat{B}_N &:= \frac{1}{2 \pi^2} \int_{-2}^2 \int_{-2}^2 \left( \frac{F(x) - F(y) }{x-y} \right)^2 \frac{4-xy}{\sqrt{4-x^2} \sqrt{ 4 - y^2} } \d x \d y - \frac{1}{ 2 \pi^2} \left( \int_{-2}^2 \frac{F(x) x}{ \sqrt{4-x^2}} \right)^2 \\
\hat{A}_N &:= - \frac{1}{2 \pi} \int_{-2}^2 \frac{F(x)(1-4+ 2x^2)}{\sqrt{4 -x^2}} \d x + \frac{ \log (\hatgam^2 -4)}{4} \\
F(x) &:= \log | \hatgam -x | .
\end{align}
But by Corollary 2.5 of \cite{functional-CLT}, the re-scaled fluctuations of the linear spectral statistic and the resolvent are asymptotically independent. The theorem statement now follows from this as well as Theorem \ref{thm:iso-clt}, with
\begin{align} \label{eqn:gaussian-variance-f}
\hat{V}_N := \frac{1}{4} \left( \hat{B}_N + V_N N \theta^2 \kappa^{-1/2} \right) , \qquad \hat{E}_N := C_N + \frac{ \hat{B}_N^{1/2} \hat{A}_N}{2N}
\end{align}
Here $V_N$ is the quantity from the other part of the theorem, controlling the size of resolvent fluctuations. 
This completes the proof of Theorem \ref{thm:mr-gauss-fe}. \qed

\subsection{Overlap with external field} \label{sec:gauss-ext}

This section will prove Theorem \ref{thm:mr-gauss-ext}, which will involve the calculation of the Laplace transform,
\beq \label{eqn:gauss-lt-a1}
 \langle e^{ \lambda v^T \sigma } \rangle.
\eeq
We will assume that $|\lambda | \leq C N^{-1/2}$ as well as \eqref{eqn:gaussianassumption}. 
 From Proposition \ref{prop:rep} we see that the Laplace transform \eqref{eqn:gauss-lt-a1} of the overlap with the external field can be represented as a ratio of two contour integrals, one involving the function $G(z)$ as above, and the second involving the function,
\begin{align}
G_u (z) := G (z) - u N^{-1} v^T (M-z )^{-1} v.
\end{align}
where the parameter $u$ is defined as
\beq
u = 2 h \lambda + \lambda^2 \beta^{-1} = 2  \lambda ( \theta / \beta)^{1/2}  + \lambda^2 \beta^{-1}.
\eeq
Note that under \eqref{eqn:gaussianassumption} and $|\lambda| \leq CN^{-1/2}$ we have for $N$ large enough,
\beq \label{eqn:gaussianoverlapu}
|u| \leq N^{-\tau/100} (1 + \theta /\kapgam ) \kapgam.
\eeq
We extend our expansion Proposition \ref{prop:gaussianexpansion} to the contour integral with $G_u (z)$.  The estimate \eqref{eqn:gaussianoverlapu} implies that in the low temperature regime we have $|u| \ll \theta$ (because $\kappa \leq C \theta$ here) and in the high temperature regime that $|u| \ll \kappa$ (as $\theta \leq C \kappa$ here) and so \eqref{eqn:gaussianassumption} will be satisfied for $\theta$ replaced by $\theta + u$.  The point of this discussion is then that we can apply the expansion obtained in the previous subsection in Proposition \ref{prop:gaussianexpansion} to calculate the Laplace transform.  In what follows we seek to expand the quantities coming from the saddle point analysis applied to $G_u$ around those coming from $G$.

 We now define $g_u(z)$ and $\hatgam_u$ by
\beq
g_u(z) = g(z)  - u \msc ( z), \qquad 
g_u' ( \hatgam_u ) = 0.
\eeq
We have,
\beq \label{eqn: raubtier}
g' ( \hatgam_u ) - g' (\hatgam) =   u  \msc' ( \hatgam_u ).
\eeq
Taylor expanding this relation, we find,
\beq
( \hatgam_u - \hatgam ) \asymp \frac{ u}{ 1  + \theta /\kapgam }.
\eeq
Substituting this estimate back into the Taylor expansion of \eqref{eqn: raubtier} we find,
\begin{align}
( \hatgam_u - \hatgam ) =&  \frac{ u \msc' ( \hatgam ) }{g'' ( \hatgam ) } \nonumber- \frac{1}{2 g''( \hatgam )} g''' ( \hatgam ) ( \hatgam_u - \hatgam )^2 + u \frac{ \msc'' ( \hatgam)}{ g'' ( \hatgam ) }  ( \hatgam_u - \hatgam ) 
+  \O \left( (1 + \theta / \kapgam )^{-3} u^3\kapgam^{-2} \right) \nonumber \\
= &  \frac{ u \msc' ( \hatgam ) }{ g'' ( \hatgam ) } - \frac{1}{ 2 g'' ( \hatgam ) } g''' ( \hatgam) \frac{ u^2 \msc' ( \hatgam)^2 }{ g'' ( \hatgam )^2 } + u^2 \frac{ \msc' ( \hatgam ) \msc'' ( \hatgam ) }{ (g'' ( \hatgam ))^2}
+  \O \left( (1 + \theta/ \kapgam )^{-3} u^3 \kapgam^{-2} \right).
\end{align}
Defining then $\gammu$ by
\beq
G_u ( \gammu  ) = 0
\eeq
a similar calculation to that above (using Theorem \ref{thm:ll-derivs}) gives the following estimates with overwhelming probability.  First,
\beq \label{eqn:gamexp-a1}
(\gammu - \gamma) \asymp \frac{u}{1 + \theta / \kapgam}.
\eeq
Second:
\begin{align}
 ( \gammu - \gamma ) &=  \frac{ u N^{-1} v^T (M-\gamma)^{-2} v }{ G'' ( \gamma ) } - \frac{G''' ( \gamma ) }{ 2 G'' ( \gamma ) } \frac{ u^2 (N^{-1} v^T (M- \gamma)^{-2} v )^2 }{ G'' ( \gamma)^2} \nonumber \\
+& u^2 \frac{ N^{-1} v^T (M- \gamma)^{-2} v N^{-1} 2 v^T (M - \gamma)^{-3} v }{ (G'' ( \gamma ))^2}   \label{eqn:gamexp-a2}
+ N^{\eps} \O \left( (1 + \theta / \kapgam )^{-3} u^3 \kapgam^{-2} \right),
\end{align}
 for any $\eps >0$.    With these preparations, we can prove the following.
\bet \label{thm:gauss-ext}
Suppose that $|t| \leq C  $ and that \eqref{eqn:gaussianassumption} holds.  Then there is an $\eps >0$ so that with overwhelming probability,
\begin{align}
\log \left( \langle \e^{ t N^{-1/2} v^T \sigma } \rangle \right) &= t^2 A_N + t B_N + \O (N^{-\eps} )
\end{align}
where,
\beq
A_N :=\frac{1}{ \beta} \left( - \frac{ \theta \msc' ( \hatgam)^2}{g'' ( \hatgam)} - \frac{1}{2} \msc ( \hatgam ) \right),
\eeq
and
\beq
B_N := - N^{1/2} \theta^{1/2} (N^{-1} v^T(M- \gamma)^{-1} v ) \beta^{-1/2}.
\eeq
\eet
\proof 
Proposition \ref{prop:gaussianexpansion} implies that with overwhelming probability, with $\lambda = t N^{-1/2}$,
\begin{align}
\log \langle \e^{ \lambda v^T \sigma } \rangle &= \frac{N}{2} \left( G_u ( \gammu) - G ( \gamma ) \right) - \frac{1}{2} \log( G'' ( \gamma) / G_u'' ( \gammu ) ) + \O(N^{-\tau/10} ).
\end{align}
Using \eqref{eqn:gamexp-a1} we find that with overwhelming probability,
\begin{align}
\frac{N}{2} (G_u ( \gammu ) - G ( \gamma ) ) &= N \frac{G'' ( \gamma )}{4} ( \gammu - \gamma )^2 - \frac{ u}{2} v^T(M-\gammu)^{-1} v + N^{\eps} \O\left( N |u|^3 \kapgam^{-3/2} (1 + \theta /\kapgam)^{-2} \right)
\end{align}
Applying \eqref{eqn:gamexp-a2} we find that with overwhelming probability,
\begin{align}
N \frac{ G'' ( \gamma ) }{ 4} ( \gammu - \gamma )^2 = & N u^2 \frac{1}{4} \frac{  (N^{-1} v^T ( M - \gamma)^{-2} v )^2}{ G'' ( \gamma ) } +  N^{\eps} \O\left( |u|^3 N\kapgam^{-3/2} (1 + \theta / \kapgam )^{-2}  \right).
\end{align}
We also calculate the linear contribution,
\begin{align}
\frac{N}{2} u N^{-1} v^T (M  - \gammu )^{-1} v &= \frac{N}{2} u N^{-1} v^T (M  - \gamma )^{-1} v + \frac{N}{2} u ( \gammu - \gamma ) N^{-1} v^T ( M - \gamma )^{-2} v \nonumber \\
&+ N^{\eps}  \O\left( |u|^3 N \kapgam^{-3/2} (1 + \theta / \kapgam )^{-2}  \right)\nonumber \\
&=  \frac{N}{2} u N^{-1} v^T (M  - \gamma )^{-1} v + \frac{N}{2} u^2 \frac{  (N^{-1} v^T (M - \gamma )^{-2} v)^2}{G'' ( \gamma ) } \nonumber \\
&+ N^{\eps}  \O\left( |u|^3 N \kapgam^{-3/2} (1 + \theta / \kapgam )^{-2}  \right)
\end{align}
We put these two calculations together and substitute back $u = 2 \lambda ( \theta / \beta)^{1/2} + \lambda^2 \beta^{-1}$, to find that with overwhelming probability,
\begin{align}
&\frac{N}{2} ( G ( \gammu ) - G ( \gamma ) - u N^{-1} v^T (M - \gammu)^{-1} v ) \nonumber\\
=& \lambda^2 \frac{N}{ \beta} \left( - \theta \frac{  (N^{-1} v^T (M - \gamma )^{-2} v )^2}{ G'' ( \gamma ) } - \frac{1}{2} ( N^{-1} v^T (M - \gamma )^{-1} v ) \right) - N \lambda ( \theta \beta^{-1} )^{1/2}  (N^{-1} v^T ( M - \gamma )^{-1} v )  \nonumber \\
+&  \O\left( |u|^3 N \kapgam^{-3/2} (1 + \theta / \kapgam )^{-2}  \right)+ \O \left( N ( | \lambda|^3 \theta^{1/2} + \lambda^4 )\kapgam^{-1/2} ( 1 + \theta / \kapgam)^{-1} \right). \label{eqn:ext-1}
\end{align}
Note that with overwhelming probability,
\beq
\log ( G'' ( \gamma) ) - \log (G_u'' (\gammu ) ) = \O \left( \frac{ |u|}{  \kapgam (1+ \theta /\kapgam ) } \right). \label{eqn:ext-1a}
\eeq
We now let $\lambda = t N^{-1/2}$ for $|t| \leq C$ to see that the error terms in \eqref{eqn:ext-1} and \eqref{eqn:ext-1a} are $\O (N^{-c})$ for some $c>0$.  Finally, using Lemma \ref{lem:sad-est}, Theorems \ref{thm:locallaw} and \ref{thm:iso} we see,
\beq
\left|  \left( - \theta \frac{  (N^{-1} v^T (M - \gamma )^{-2} v )^2}{ G'' ( \gamma ) } - \frac{1}{2} ( N^{-1} v^T (M - \gamma )^{-1} v ) \right) -\left( - \frac{ \theta ( \msc' ( \hatgam) )^2}{ g'' ( \hatgam) } - \frac{1}{2} \msc ( \hatgam) \right) \right| \leq N^{-c}
\eeq
with overwhelming probability for some $c>0$.  This yields the claim. \qed

Theorem \ref{thm:mr-gauss-ext} follows from Theorem \ref{thm:gauss-ext} and the lemma below which computes the coefficient $A_N$ and also finds its order of magnitude. The proof is deferred to Section \ref{a:gauss-1}

\bel \label{thm:mr-gauss-over-b1} We have that $A_N = \frac{\msc' ( \hatgam)}{ 2 g'' ( \hatgam)}\left( - \msc ( \hatgam) + 2 \theta ( \msc' ( \hatgam))^2 \right) \asymp 1$. Moreover, there is a $c_1 >0$ so that the following holds.  For any $D>0$ with probability at least $1-N^{-D}$ and $N$ large enough,
\beq
\langle N^{-1/2} v \cdot \sigma \rangle = N^{-1/2} h v^T (M- \gamma)^{-1} v + \O(N^{-c_1} ).
\eeq
\eel

\subsection{Overlap between two replicas} \label{sec:gauss-overlap}
In this section, we consider the asymptotic fluctuations of the overlap,
$
R_{12} := \frac{1}{N} \sigma^{(1)} \cdot \sigma^{(2)},
$
and prove Theorem \ref{thm:gauss-overlap}. 
From Proposition \ref{prop:rep} we have,
\beq \label{eqn:over-rep-1-f}
\left\langle \exp \left[ \beta t \sigma^{(1)} \cdot \sigma^{(2)} \right] \right\rangle =  \frac{ \int_{\Gamma^2} \exp \frac{N}{2}G_1 (z, w) \d z \d w }{ \left( \int_\Gamma \exp \frac{N}{2} G (z) \d z \right)^2}
\eeq
where $\Gamma$ is a vertical line in the complex plane lying to the right of the spectrum of $M$.  Here the function $G(z)$ is as above, and
\beq
G_1 (z, w) := \beta(z + w ) - \frac{1}{N} \sum_{i=1}^N \log ( (z- \lambda_i )( w - \lambda_i ) - t^2 ) - \frac{ \theta}{N} \sum_{i=1}^N v_i^2 \frac{ ( 2 \lambda_i - z - w ) - 2 t }{( \lambda_i - w  )( \lambda_i - z ) - t^2 }.
\eeq
Recall the notation $
v_i = v^T u_i (M)$. 
 The following lemma provides an expansion of $G_1$ around a point $\gamh$. The proof is deferred to Appendix \ref{a:gauss-1-f}. The point $\gamh$ will eventually be chosen to be an approximate saddle for $G_1 (z, w)$.
\bel \label{lem:gauss-1-f}
Assume,
\beq
t^2 \leq \log(N) \frac{ \kappa^{1/2}}{N}
\eeq
and let $\gamh \in \rr$ satisfy $| \gamh - \gamma | \leq \kappa / \log(N)$.  Let $\delta >0$ be sufficiently small.  For $z$ and $w$ satisfying
\beq \label{eqn:zw-ass}
|z - \gamh | + |w - \gamh | \leq N^{ \delta}  \frac{ \kappa^{1/4}}{ N^{1/2} (1 + \theta / \kappa)^{1/2} },
\eeq
we have with overwhelming probability,
\begin{align} \label{eqn:gauss-1-f}
G_1 (z, w) &= \beta(z+w) - N^{-1} \sum_i \log( ( \lambda_i -z)( \lambda_i -w )) - \theta N^{-1} [ v^T (M- z +t)^{-1} v + v^T (M - w+t)^{-1} v] \nonumber \\
&+ t^2 m' ( \gamh ) \nonumber \\
&+N^{\eps} \O ( t^4 \kappa^{-5/2} + t^2\kappa^{-5/4} N^{\delta-1/2} ) \nonumber\\
&+ N^{\eps}\O ( \theta \kappa^{-2}N^{-1}N^{2 \delta}(1+ \theta / \kappa)^{-1} (|t| + t^2 \kappa^{-1} + N^{\delta} \kappa^{1/4} N^{-1/2} ).
\end{align}
for any $\eps >0$.
\eel

  We are now ready to prove the following.
\bep \label{prop:gauss-over-f}
With overwhelming probability, we have for a $c>0$, and $t$ satisfying,
\beq
t^2 \leq \log(N) \sqrt{ \kapgam} N^{-1}
\eeq
the expansion,
\begin{align}
\log \langle \e^{t \beta \sigma^{(1)} \cdot \sigma^{(2)} } \rangle  &= \frac{N t^2}{2} m' ( \gamma) \left( 1 - \theta (2 N^{-1} v^T (M- \gamma)^{-3} v) (G'' ( \gamma)^{-1} \right)  + t \theta v^T (M- \gamma)^{-2} v \nonumber\\
+& \O (N^{-c} ).
\end{align}
\eep
\proof  
Define the function, (it appears naturally on the RHS of \eqref{eqn:gauss-1-f})
\beq
G_2 (z) := \beta z - \frac{1}{N} \sum_i \log( z - \lambda_i ) - \theta N^{-1} v^T (M- z + t )^{-1} v.
\eeq
and define $\gamh$ by
\beq
\beta + m ( \gamh ) = N^{-1} \theta v^T (M - \gamh + t )^{-2} v.
\eeq
From the equation,
\beq
G' ( \gamh ) - G' (\gamma) = N^{-1} \theta ( v^T (M- \gamh + t)^{-2} v - v^T (M- \gamh)^{-2} v ),
\eeq
we see that with overwhelming probability we have, for any $\eps >0$,
\beq \label{eqn:gam-exp-b1}
| \gamh - \gamma| \leq C t \theta \kappa^{-1} ( 1+ \theta / \kappa)^{-1} \leq C t, \qquad \gamh - \gamma = \frac{ -2t \theta N^{-1} v^T(M- \gamma)^{-3} v}{ G'' ( \gamma) } + N^{\eps} \O \left( t^2/ \kappa \right).
\eeq
From the equation
\beq
G' ( \gamh - t) - G' ( \gamma ) = m ( \gamh -t) - m ( \gamh),
\eeq
we see that
\beq
| \gamh - t - \gamma| \leq C t (1+ \theta/ \kappa)^{-1}
\eeq
and
\begin{align}
 \gamh - t - \gamma &= \frac{ - t m' ( \gamh - t)}{ G'' ( \gamma) }  +\O( t^2 (1 + \theta/\kappa)^{-1} \kappa^{-1} ) \nonumber \\
&= \frac{ - t m' ( \gamma)}{ G'' ( \gamma)} + \O( t^2 (1 + \theta/\kappa)^{-1} \kappa^{-1} ) , \label{eqn:gam-exp-b2}
\end{align}
with overwhelming probability. 
Note that by expanding $\gamh -t -\gamma$ instead of $\gamh - \gamma$ we have picked up some extra factors of $(1+\theta / \kapgam)^{-1}$ which will prove useful later.

With $z = \gamh + \i s_1$ and $w = \gamh + \i s_2$ and $|s_i| \leq N^{\delta} \kapgam^{1/4} N^{-1/2} (1 + \theta / \kapgam)^{1/2}$ we see that by Lemma \ref{lem:gauss-1-f} and a Taylor expansion that
\begin{align} \label{eqn:G1-expand-f}
\frac{N}{2} (G_1 ( z, w) -2 G_2 ( \gamh ) -t^2 m'_N ( \gamh) ) &= - N \frac{G_2'' ( \gamh) }{4} (s_1^2 + s_2^2 ) + \O (N^{-c})
\end{align}
for some $c>0$, taking $\delta >0$ sufficiently small, with overwhelming probability.  Note that since by the above estimates and assumptions on $t$ we see that $| \gamh - \gamma| + | \gamh - \gamma - t | \leq N^{-\eps} \kapgam$. Therefore,  we have that,
\beq
G_2'' ( \gamh)  \asymp \frac{1}{ \kapgam^{1/2}} (1+ \theta / \kapgam ).
\eeq
Therefore, a straightforward modification of the saddle point analysis presented in Proposition \ref{prop:gaussianexpansion} gives (i.e., using \eqref{eqn:over-rep-1-f} and \eqref{eqn:G1-expand-f} to control the contour integral in the numerator),
\begin{align}
\log \langle \e^{t \beta \sigma^{(1)} \cdot \sigma^{(2)} } \rangle &= N( G_2 ( \gamh  ) - G ( \gamh ) ) + \frac{N t^2 m' ( \gamh)}{2}  \nonumber\\
+& \log(G_2'' ( \gamh) )  - \log ( G'' ( \gamma) ) +  \O ( N^{-c})
\end{align}
for some $c>0$, with overwhelming probability.  From our expansions of $\gamh$ in terms of $\gamma$ and $t$ it is easy to see that
\beq
\left| t^2 m' ( \gamh) )  - t^2 m' ( \gamma) )  \right| + \left| \log(G_2'' ( \gamh) )  - \log ( G'' ( \gamma) ) \right| \leq N^{-c}.
\eeq
Therefore we have the following estimate with overwhelming probability,
\begin{align}
\log \langle e^{t \beta \sigma^{(1)} \cdot \sigma^{(2)} } \rangle &= N( G_2 ( \gamh  ) - G ( \gamma ) ) + \frac{N t^2 m' ( \gamma)}{2} + \O(N^{-c}).
\end{align}
Now, Taylor expanding the two components of $G_2 (z)$ around their value at $\gamma$ we find,
\begin{align}
G_2 ( \gamh ) - G ( \gamma) &= (\beta+ m ( \gamma) )( \gamh - \gamma) - \theta N^{-1} v^T(M - \gamma)^{-2} v ( \gamh - \gamma - t) \nonumber \\
&+ \frac{m' ( \gamma)}{2}  ( \gamh - \gamma)^2  - \theta N^{-1} v^T (M - \gamma)^{-3} v ( \gamh - \gamma - t)^2 \nonumber \\
&+\O (|t|^3 \kappa^{-3/2} ) .
\end{align}
For the first line of the above, using $G' ( \gamma)=0$ we see,
\beq
 (\beta+ m ( \gamma) ) ( \gamh - \gamma) - \theta N^{-1} v^T(M - \gamma)^{-2} v ( \gamh - \gamma - t) = t \theta N^{-1} v^T (M- \gamma)^{-2} v.
\eeq
Applying now \eqref{eqn:gam-exp-b1} and \eqref{eqn:gam-exp-b2} we have with overwhelming probability,
\begin{align}
 & \frac{m' ( \gamma)}{2}  ( \gamh - \gamma)^2  -\theta N^{-1} v^T (M - \gamma)^{-3} v ( \gamh - \gamma - t)^2 \nonumber \\
= &- \frac{m' ( \gamma) \theta 2 N^{-1} v^T (M - \gamma)^{-3}v }{2 G'' ( \gamma) } t^2 + N^{\eps} \O ( |t|^3 \kappa^{-3/2} ).
\end{align}
Hence,
\begin{align}
\log \langle \e^{t \beta \sigma^{(1)} \cdot \sigma^{(2)} } \rangle  &= \frac{N t^2}{2} m' ( \gamma) \left( 1 - \theta (2 N^{-1} v^T (M- \gamma)^{-3} v) G'' ( \gamma)^{-1} \right)  + t \theta v^T (M- \gamma)^{-2} v \nonumber\\
+& \O (N^{-c} ).
\end{align}
with overwhelming probability, for some $c>0$. \qed

From Theorem \ref{thm:ll-derivs} as well as Lemma \ref{lem:sad-est} we see that,
\begin{align}
\sqrt{ \kappa} \left|  m' ( \gamma) \left( 1 - \theta (2 N^{-1} v^T (M- \gamma)^{-3} v) (G'' ( \gamma)^{-1} \right) - \msc' ( \hatgam) \frac{ \msc' ( \hatgam) - 2 \theta \msc'' ( \hatgam) }{\msc' ( \hatgam) - \theta \msc'' ( \hatgam)} \right| \leq N^{-c}
\end{align}
for some $c>0$ with overwhelming probability.  Theorem \ref{thm:gauss-overlap} follows from this and Proposition \ref{prop:gauss-over-f}.

The following is proved in a similar fashion to the second statement of Lemma \ref{thm:mr-gauss-over-b1} and the proof is therefore omitted. It again simply verifies that the linear term in the expansion of the Laplace transform in Theorem \ref{thm:gauss-overlap} is indeed to leading order the Gibbs expectation. 
\bel  \label{thm:mr-gauss-b2}
 There is a $c_1 >0$ so that the following holds.  For any $D>0$ with probability at least $1-N^{-D}$ and $N$ large enough,
\beq
\frac{ \kappa^{1/4}}{N^{1/2}} \langle \sigma^{(1)} \cdot \sigma^{(2)} \rangle = N^{-1/2} \kappa^{1/4} h v^T (M  - \gamma)^{-2} v + \O(N^{-c_1} ).
\eeq
\eel

\section{Intermediate magnetic field $h = \O (N^{-1/6} )$ and $ \beta >1 $} \label{sec:int}

In this section we will prove the results in Section \ref{sec:mr-int}.  We assume that the parameters $\beta$ and $h$ satisfy,
\beq \label{eqn:int-ass}
1+c \leq \beta \leq c^{-1}, \qquad h^2 \beta = \theta N^{-1/3}
\eeq
for a fixed $c>0$ and fixed $\theta > 0$.

\subsection{Free energy}
\label{sec:int-fe}

In this section we investigate the free energy and prove Theorem \ref{thm:mr-int-fe-1}.  
 From our choice of $h^2 \beta = \theta N^{-1/3}$ and Proposition \ref{prop:rep} we are led to investigate the function,
\beq \label{eqn:int-fe-Gz-def}
G(z) =\beta z - \frac{1}{N} \sum_{i=1}^N \log ( z - \lambda_i  ) - \frac{\theta}{N^{4/3}} v^T \frac{1}{M-z} v.
\eeq
We recall the definition of $H = M+V$ with $V$ a diagonal matrix of independent centered Gaussians of variance $2/N$, so that $H$ is a GOE matrix.  We define the quantities $\gamma$, $\hatgam_M$ and $\hatgam_H$ as the largest solutions to,
\begin{align}
\beta + m ( \gamma ) = \frac{ \theta}{N^{4/3}} v^T \frac{1}{ (M - \gamma)^2} v, \quad 
\beta - 1 = \frac{ \theta}{N^{4/3}}v^T \frac{1}{ (M - \hatgam_M)^2} v,  \quad 
\beta - 1 = \frac{ \theta}{N^{4/3}} v^T \frac{1}{ (H - \hatgam_H )^2} v . \label{eqn:int-gam-def}
\end{align}
From Theorem \ref{thm:compare} and the fact that $u_1 (H)$ is uniformly distributed on the $N-1$ sphere, we have the following lemma.
\bel \label{lem:Mulb}
For all sufficiently small $\eps >0$ we have,
\beq 
\pp[ (v^T u_1(M))^2  \leq N^{-\eps}  ] \leq N^{-\eps/2} ,
\eeq
for $N$ large enough.
\eel

We now define an event $\Feps$ on which a number of estimates for the eigenvalues and eigenvectors of $M$ and $H$ hold. 
\begin{defn} \label{def:int-F}
Let $\eps_1 >0$ be sufficiently small.  We define $\Feps$ to be the event that all of the following estimates hold.  First, we assume that the following  \emph{level repulsion} estimates hold,
\beq
|\lambda_1(M) - \lambda_2(M)| \geq N^{-2/3-\eps_1/1000}, \quad |\lambda_1(H) - \lambda_2(H)| \geq N^{-2/3-\eps_1/1000}.
\eeq
We assume the following lower bound for the projection of the largest eigenvectors of $M$ and $H$ onto $v$,
\beq \label{eqn:F-v-lb}
 (v^T u_1(M))^2 \geq N^{-\eps_1/1000}, \quad (v^T u_1(H) )^2 \geq N^{-\eps_1/1000}.
\eeq
We assume that the \emph{rigidity} estimates of Theorem \ref{thm:rigi} hold with $\eps = \eps_1 / (10^6)$, for the eigenvalues of $M$ and $H$.  We also assume the \emph{delocalization} bounds,
\beq
\sup_i (v^T u_i(M))^2 + (v^T u_i (H) )^2 \leq N^{\eps_1/10^6}.
\eeq
We assume that the events of Lemma \ref{lem:HViso} hold with $\omega = \eps_1  / 10^6$ in the definitions of $\mathcal{S}_1$ and $\mathcal{S}_2$ and $\delta = 1/1000$ and $\eps = \eps_1/10^6$.  With the same choice of parameters we also assume both \emph{isotropic} estimates of Theorem \ref{thm:iso} hold with $u, w = N^{-1/2} v$.  We assume also that the events of Theorem \ref{thm:compare} hold with $\eps = \eps_1 / 10^6$ in the first estimate and $c = \eps_1 /1000$ in the second. \qed
\end{defn}
From all of the results in Section \ref{sec:rmt} we have that 
there is a small $c_e >0$ so that for sufficiently small $\eps_1 >0$ that $\Feps$ holds with probability at least $1 - N^{-c_e \eps_1}$, for $N$ large enough.

Let us introduce the parameters $s, \hats_M$ and $\hats_H$ via,
\beq
\gamma = \lambda_1 (M) + \frac{s}{N^{2/3}}, \qquad \hatgam_M = \lambda_1 (M) + \frac{ \hats_M}{N^{2/3}}, \qquad \hatgam_H = \lambda_1 (M) + \frac{ \hats_H}{N^{2/3}}.
\eeq
Note that a-priori it is possible for $\hats_H <0$ but the lemma below shows that this does not occur on the event $\Feps$. 
For these quantities we have the following preliminary estimates. The proof is deferred to Appendix \ref{a:int-gams-comp}. 
\bel \label{lem:int-gams-comp}
For all sufficiently small $\eps_1 >0$ we have on the event $\Feps$ that,
\beq \label{eqn:int-gams-comp-1}
N^{-\eps_1 /100} \leq \hats_M, \hats_H, s \leq N^{\eps_1 / 100}, \qquad |s - \hats_M| \leq N^{ \eps_1/5-1/3}.
\eeq
There is moreover a small $c_1 >0$ so that if $\eps_1 >0$ is sufficiently small then on $\Feps$,
\beq \label{eqn:int-gams-comp-2}
| \hats_M - \hats_H | \leq N^{-c_1}, \qquad | \hatgam_M - \hatgam_H | \leq N^{-2/3-c_1}.
\eeq
\eel

With these preliminary estimates on the saddle locations, we turn to the remainder of the proof.  The method of steepest descent will allow us to derive an expansion for the free energy in terms of the saddle $\gamma$; our estimates above then allow us to replace this by quantities involving only the spectral quantities of $H$.  The convergence of this latter system will be carried out in the next section.

For  all sufficiently small $\eps_1 >0$, we have on the event $\Feps$ that, due to  the level repulsion, rigidity and the estimates of Lemma \ref{lem:Mulb} that,
\beq
N^{2k/3 - \eps_1/5 } \leq | G^{ (k+1)} ( \gamma ) | \leq N^{ 2k/3 + \eps_1/5 }, \qquad k=1, 2
\eeq
and the expansion,
\beq
G ( \gamma + \i t ) - G ( \gamma ) = -\frac{ G'' ( \gamma ) }{2} t^2 + \O ( N^{-1-1/10} )
\eeq
for $|t| \leq N^{-5/6+1/100}$.  Using the above estimates, a straightforward modification of the steepest descent argument in Proposition \ref{prop:gaussianexpansion} gives the following. We omit the proof for brevity.
\bep  \label{prop:int-steepest} Assume that \eqref{eqn:int-ass} holds.  
There is a small $c_1>0$ so that for sufficiently small $\eps_1 >0$ we have that on the event $\Feps$, the following expansion for the free energy holds
\begin{align}
\frac{1}{N} \log Z_{N ,\beta, h} &= \frac{1}{2} G ( \gamma) + \frac{1}{N} \log ( \Gamma (N/2) ) \nonumber \\
&+ \frac{1}{N} (N/2-1) \log(2/ (N \beta )) - \frac{1}{2N} \log (N G'' ( \gamma) \pi) + \O(N^{-1-c_1} ),
\end{align}
where $\gamma$ is the saddle defined by \eqref{eqn:int-gam-def}. 
\eep

We expect the fluctuations of the first term, $G( \gamma)$ to be $\O (N^{-2/3})$, so we will ignore the $N^{-1}$ log term in examining fluctuations.  We now replace the quantities that we obtained from the saddle point analysis with spectral quantities involving only the GOE matrix $H$; the convergence of the latter will be easier to analyze later.
\bel \label{lem:int-fe-1}  There is a small $c_1 >0$ so that for sufficiently small $\eps_1>0$,  on the event $\Feps$ it holds that,
\begin{align}
&N^{2/3} \left( F_{N, \beta, h}- \frac{1}{N} \log ( \Gamma (N/2)) + \frac{1}{N} ( N/2-1) \log (2 / (N \beta )) -\mathcal{C} \right) \nonumber\\
=& \frac{1}{2} \left( (\beta-1) N^{2/3}( \hatgam_H - 2) - \frac{ \theta}{N^{2/3}} v^T (H - \hatgam_H)^{-1} v  \right) + \O(N^{-c_1} ), \label{eqn:int-fe-a1}
\end{align}
where 
$
\mathcal{C} := - \int_{-2}^2 \log(2-x) \d \rhosc (x). 
$
\eel
\proof We begin by investigating the term $G(\gamma)$ appearing in the main estimate of Proposition \ref{prop:int-steepest}. Using the rigidity estimates of the event $\Feps$ and the estimates of Lemma \ref{lem:int-gams-comp} we first have,
\begin{align}
\beta ( \gamma - 2)  - \frac{1}{N} \sum_{i=1}^N \log ( \gamma - \lambda_i ) = ( \beta - 1 ) ( \hatgam_H - 2 ) +\mathcal{C}+ \O (N^{-2/3-c} )
\end{align}
for some $c>0$.   
We estimate
\begin{align}
N^{-4/3} \left| v^T (H - \hatgam_H )^{-1} v - v^T (M - \gamma )^{-1} v \right| &\leq N^{-4/3} \left| v^T (H - \hatgam_H )^{-1} v - v^T (M - \hatgam_H )^{-1} v \right| \label{eqn:int-est-a1}   \\
&+N^{-4/3} \left| v^T (H - \hatgam_H )^{-1} v - v^T (H - \gamma )^{-1} v \right|  \label{eqn:int-est-a2}
\end{align}
The term \eqref{eqn:int-est-a1} may be estimated in a similar fashion to \eqref{eqn:aa} with $E = \hatgam_H$; indeed, by the same proof of our estimate of \eqref{eqn:aa} we find that \eqref{eqn:int-est-a1} is $\O(N^{-2/3-c})$ for some $c>0$.  The same estimate holds for \eqref{eqn:int-est-a2} given the second estimate of \eqref{eqn:int-gams-comp-2}.  
Therefore,  on $\Feps$
\beq
N^{-4/3} \left| v^T (H - \hatgam_H )^{-1} v - v^T (M - \gamma )^{-1} v \right| \leq N^{-2/3-c}.
\eeq
 This concludes the proof. \qed

This completes our expansion of the free energy in terms of GOE quantities.  We now prove Theorem \ref{thm:mr-int-fe-1}.

\noindent{\bf Proof of Theorem \ref{thm:mr-int-fe-1}.}  The random variable $X_N$ is the first term on the RHS of \eqref{eqn:int-fe-a1}.  The eigenvectors and eigenvalues of the GOE matrix $H$ are independent, and the eigenvector matrix is uniformly distributed on the orthogonal group.  Hence, the quantities $\{ (v^T u_i (H))^2\}_{i=1}^N$ have the same joint distribution as $\{N  g_i^2 / \sum_j g_j^2 \}_{i=1}^N$ where the $\{ g_i \}_{i=1}^N$ are iid standard normal random variables.  Let now $\{ \mu_i \}_{i=1}^N$ be a vector of GOE eigenvalues in decreasing order independent of the $\{g_i\}_i$.

Let $x_a$ be the largest solution to
\beq
\beta -1 = \frac{\theta}{N^{4/3} } \frac{1}{ N^{-1} \sum_{i=1}^N g_i^2 } \sum_{i=1}^N \frac{ g_i^2}{ ( \mu_i -x_a )^2 }.
\eeq
Then the main term $X_N$ in the expansion of the free energy has the same distribution as,
\beq
(\beta-1) N^{2/3} ( \hatgam_H -2 ) - \frac{ \theta}{N^{2/3}} v^T (H - \hatgam_H)^{-1} v \dist ( \beta-1)N^{2/3} ( x_a - 2) - \frac{1}{ N^{-1} \sum_{i=1}^N g_i^2} \frac{ \theta}{N^{2/3}} \sum_{i=1}^N \frac{ g_i^2}{ ( \mu_i - x_a )}.
\eeq
The quantity on the RHS is the random variable $Y_N$ in the theorem statement.  
Consider now the alternative system where $x_b$ is the solution to
\beq
\beta -1 = \frac{ \theta}{N^{4/3}} \sum_{i=1}^N \frac{ g_i^2}{ ( \mu_i - x_b )^2}
\eeq
and the quantity,
\beq
( \beta - 1) N^{2/3} ( x_b - 2) - \frac{ \theta}{N^{2/3}} \sum_{i=1}^N \frac{g_i^2}{ ( \mu_i - x_b )}.
\eeq
Now we have that 
\beq
N^{-1} \sum_{i=1}^N g_i^2 = 1 + \O(N^{-1/2+\eps} ), \qquad \sup_i |g_i| \leq N^{\eps}
\eeq
 with overwhelming probability for any $\eps >0$.    On the event that the above estimates as well as the rigidity estimates hold for $\mu_i$ with sufficiently small $\eps >0$ it is easy to see that
\beq
\left| \frac{ \theta}{N^{4/3}} \sum_{i=1}^N \frac{ g_i^2}{ ( \mu_i - E )^2} - \frac{\theta}{N^{4/3} } \frac{1}{ N^{-1} \sum_{i=1}^N g_i^2 } \sum_{i=1}^N \frac{ g_i^2}{ ( \mu_i - E )^2 } \right| \leq N^{3 \eps_2+10\eps} N^{-1/2}
\eeq
for $E$ satisfying $N^{-2/3-\eps_2 } \leq E- \mu_1$.   Hence, similar arguments as those leading our estimates \eqref{eqn:int-gams-comp-2} yield,
\beq
|x_a - x_b | \leq N^{-5/6}
\eeq
with probability at least $1-N^{-\eps}$ for sufficiently small $\eps >0$.  Furthermore, we write,
\begin{align}
 \frac{1}{ N^{-1} \sum_{i=1} g_i^2 } \frac{ \theta}{ N^{2/3}} \sum_{i=1} \frac{ g_i^2}{ \mu_i - x_a} -\frac{ \theta}{ N^{2/3}} \sum_{i=1} \frac{ g_i^2}{ \mu_i  - x_b } &= \frac{ \theta}{ N^{2/3}} \sum_{i=1} \frac{ g_i^2}{ \mu_i  - x_a }-\frac{ \theta}{ N^{2/3}} \sum_{i=1} \frac{ g_i^2}{ \mu_i  - x_b } \notag \\
&+ \left( 1 - \frac{1}{ N^{-1} \sum_i g_i^2} \right) \frac{ \theta}{N^{2/3}} \sum_{i=1}^N \frac{ g_i^2}{ \mu_i - x_a }.
\end{align}
With probability at least $1 - N^{-\eps}$ the RHS of the first line is $\O \left( | x_a - x_b | N^{2/3+10\eps} \right)$ and the second line is $\O \left(N^{1/3+\eps} N^{-1/2} \right)$.  This concludes the proof. \qed

\subsection{Convergence of the saddle system} \label{sec:int-conv}

In this section we prove Theorem \ref{thm:conv1}. 
For this we require the following lemma, a consequence of the arguments in Section 6 of \cite{SSK-ls} (specifically, (6.33), (6.34) of \cite{SSK-ls}, the rigidity estimates of Theorem \ref{thm:rigi} and a union bound). 
\bel \label{lem:RoughParticleEst}
There is a $c_0 >0$ so that the following statement holds.  Let $\eps >0$.  Then there is a $K >0$ depending on $\eps$ so that for $N$ large enough,
\beq
\pp \left[ \bigcap_{j=K}^N  \{ c_0 j^{2/3} \leq N^{2/3} ( \mu_j - 2 ) \leq c_0^{-1} j^{2/3} \} \right] \geq 1 - \eps.
\eeq
\eel
In this section we will use the following notation. Introduce the random variable $s$ by
\beq
\xb = \mu_1 + N^{-2/3} s.
\eeq
Note that by definition $ s>0$. We first prove the following.
\bel \label{lem:conv-1}
The random variables $s$ and $s^{-1}$ are tight.
\eel
\proof For a lower bound,
\beq
\beta -1 \geq \theta \frac{g_1^2}{s^2}.
\eeq
Fix $\eps >0$ and consider the event of Lemma \ref{lem:RoughParticleEst} and the corresponding $K = K(\eps)$.  For all $k_0 > K$ we have,
\beq
\frac{\beta -1}{\theta} \leq \frac{1}{s^2} \sum_{j=1}^{k_0} g_j^2 + C \sum_{j > k_0} \frac{g_j^2}{ j^{4/3}}
\eeq
for some $C>0$ independent of $k_0$.  For $k_0$ sufficiently large, there is an event of probability at least $1 - \eps$ on which the second term is less then $\frac{\beta-1}{2 \theta}$, by Markov's inequality.  Hence, there is an event of probability at least $1- 2 \eps$ on which,
\beq
s^2 \leq \frac{2 \theta}{ \beta -1 } \sum_{j=1}^{k_0} g_j^2.
\eeq
This yields the claim. \qed

We now define $s_n$ to be the unique positive solution to,
\beq
\beta - 1 = \frac{ \theta}{N^{4/3}} \sum_{i=1}^n \frac{ g_i^2} { ( s_n N^{-2/3} + \mu_1 - \mu_i )^2 }.
\eeq
We now show that $s_n$ approximates $s$ for $n$ large enough, uniformly in $N$. The proof is deferred to Appendix \ref{a:conv-2}. 
\bel \label{lem:conv-2}
Let $ \eps >0$.  There is an $n_0$, depending on $\eps$ so that for all $N$ large enough, there is an event $\F_{n_0}$ (depending on $N$) on which the following estimate holds for all $n \geq n_0$, 
\beq
|s-s_n | \leq \eps
\eeq
and $\pp [ \F_{n_0} ] \geq 1 - \eps$.
\eel

The next lemma shows that the quantity $\xi$ is well-approximated by a quantity involving only $n$ of the Gaussians $\{ g_j \}_{j=1}^n$ and the quantity $s_n$, uniformly in $N$. The proof is deferred to Appendix \ref{a:conv-3}. 

\bel \label{lem:conv-3} Let $\eps >0$.  There is an $n_1$ so that the following holds for every $n \geq n_1$.  There is an event $\F_{n}$ which holds with probability at least $ 1- \eps$ on which,
\beq
\frac{1}{N^{2/3}} \left| \sum_{j=2}^N \left( \frac{ g_j^2}{ ( \mu_j - \mu_1 - sN^{-2/3} ) } - \frac{1}{ \gamma_j - 2 } \right) - \sum_{j=2}^n \left( \frac{ g_j^2}{ ( \mu_j - \mu_1 - s_n N^{-2/3} ) } - \frac{1}{ \gamma_j - 2 } \right) \right| \leq \eps,
\eeq
and
\beq
|g_1^2 / s - g_1^2 / s_n | \leq \eps.
\eeq
The same statement holds with the $g_j^2$ replaced by $1$ everywhere.
\eel

With the previous lemmas, we can provide the following.

\noindent{\bf Proof of Theorem \ref{thm:conv1}}.  Denote by $y_n$ the random variable
\beq
y_n := N^{2/3} ( \mu_1 + s_n N^{-2/3} - 2 ) + \frac{ \theta}{N^{2/3}} \sum_{i=1}^n \left( \frac{ g_i^2}{\mu_i -mu_1 - s_n N^{-2/3} } - \frac{1}{ \gamma_i - 2 } \right).
\eeq
Note that $y_N$ is the random variable whose convergence we are interested in.  
For any bounded Lipschitz function $F$ and any $\eps >0$ we see from Lemmas \ref{lem:conv-2} and \ref{lem:conv-3} that there is an $n$ so that for all $N$ large enough,
\beq
\left| \ee[ F ( y_N) ] - \ee[F(y_n) ] \right| \leq \eps.
\eeq
Since $g_1^2 >0$ almost surely, and that the GOE eigenvalues are a.s. distinct, we see that $s_n$ and $y_n$ are continuous functions of $\{ \mu_i \}_{i=1}^n$.  Let $\{ \chi_i \}_{i=1}^\infty$ be the Airy$_1$ random point field, normalized so that the joint limit of  $\{ N^{2/3} ( \mu_i - 2 ) \}_{i=1}^n$ is the first $n$ particles of the Airy$_1$ random point field.

If $a_n$ denotes the solution to 
\beq
\beta-1 = \theta \sum_{i=1}^n \frac{g_i^2}{ (\chi_i - \chi_1 - a_n )^2}
\eeq
and $z_n$ the random variable,
\beq
z_n := \chi_1 + a_n + \theta \sum_{i=1}^n \left( \frac{g_i^2}{ \chi_i - \chi_1 - a_n } + \frac{1}{ ( \pi 2 i/3)^{2/3} } \right).
\eeq
We have by the joint convergence of the $n$ largest GOE eigenvalues to the Airy$_1$ random point field,
\beq
\lim_{N \to \infty} \ee[ F(y_n) ] = \ee[ F(z_n ) ].
\eeq
It remains to prove that,
\beq
\lim_{n \to \infty} \ee[ F(z_n ) ] = \ee[ F(z_\infty)].
\eeq
The arguments are almost identical to the arguments involving the analagous statements for $y_n$ and $y_N$ that proved Lemmas \ref{lem:conv-1} and \ref{lem:conv-2}.  Indeed, all the estimates for the GOE eigenvalues $\mu_i$ we used have direct analogues for the Airy$_1$ random point field --  this is proven in Section 6 of \cite{SSK-ls}.

The only substantial difference is that one  truncates the analogue of \eqref{eqn:gj} using the Kolmogorov Three Series lemma, as it is an infinite random sum that is not absolutely summable.   Note that these arguments also show that the random variables $a_\infty$ and $z_\infty$ are well defined.  \qed

The same method as above gives a proof of the following generalization.
\bep
Let $\eps >0$, $k \geq 2$ an integer, and let $s_n$ be as above.  There is an $n_1 >0$ so that for all $n > n_1$ the following holds.  There is an event $\F_n$ with probability at least $1- \eps$ so that
\beq
\frac{1}{N^{2k/3}} \left| \sum_{i=1}^N \frac{g_j^2}{ ( \mu_j - \mu_1 - s N^{-2/3} )^k } - \sum_{i=1}^n \frac{g_j^2}{ ( \mu_j - \mu_1 - s_n N^{-2/3} )^k } \right| \leq \eps.
\eeq
The same statement holds with the $g_j^2$ replaced by $1$.
\eep

Similarly the following is easily proven using the same methods as above. 
\bet \label{thm:int-conv-2}
For any $k_0 \geq 2$ we have the joint convergence of the random variables,
\begin{align}
\left\{ \frac{1}{N^{2 k /3}} \sum_{i=1}^N \frac{g_i^2}{ ( \mu_i - \mu_1 - s N^{-2/3} )^k }, \frac{1}{N^{2 k /3}} \sum_{i=1}^N \frac{1}{ ( \mu_i - \mu_1 - s N^{-2/3} )^k } \right\}_{k=2}^{k_0}
\end{align}
and
\begin{align}
\left\{ s, \frac{1}{N^{2/3}} \sum_{i=1}^N \frac{ g_i^2}{ \mu_i - \mu_1 - sN^{-2/3} } - \frac{1}{ \gamma_i - 2},  \frac{1}{N^{2/3}} \sum_{i=1}^N \frac{ 1}{ \mu_i - \mu_1 - sN^{-2/3} } - \frac{1}{ \gamma_i - 2} \right\}
\end{align}
to the corresponding quantities of the Airy$_1$ random point field.
\eet

\subsection{Overlap with external field} \label{sec:int-ext}

In this short section we prove Theorem \ref{thm:mr-int-ext}.  
We continue to consider the regime of parameters $\beta >1$ and $h^2 \beta = \theta N^{-1/3}$ for a fixed $\theta >0$.  Recall also our definition of the event $\Feps$ in Definition \ref{def:int-F}.   We consider now,
\beq
\log \langle \exp \left[ \beta^{1/2} t N^{-1/6} v^T \sigma \right] \rangle
\eeq
which is then a difference of contour integrals, involving the function $G(z)$ as defined in \eqref{eqn:int-fe-Gz-def} as well as the function
\beq
G_u (z) := G(z) - u N^{-4/3} v^T (M-z)^{-1} v,
\eeq
where
\beq \label{eqn:qa-6}
u = 2 t \theta^{1/2} + t^2.
\eeq
Note that due to our choice of scaling, the function $G_u(z)$ differs slightly from the definition in Section \ref{sec:gauss}.  
We assume that
\beq
|u| \leq N^{- \alpha}
\eeq
for some $\alpha >0$ which is fixed throughout this section.  Let $\gamu, \gamma$ be defined by,
\beq
G_u' ( \gamu ) = 0, \qquad G' ( \gamma) = 0.
\eeq
Elementary arguments similar to the proof of \eqref{eqn:gamexp-a1}   show that on the event $\Feps$ that for sufficiently small $\eps_1 >0$, we have
\beq \label{eqn:jan-a1}
| \gamu - \gamma| \leq N^{\eps_1/20} u N^{-2/3}.
\eeq
As a consequence of Proposition \ref{prop:int-steepest} and representing the Laplace transform as a ratio of contour integrals using Proposition \ref{prop:rep} we obtain the following.
\bel
There is a $c_1 >0$ so that the following holds on $\Feps$ for sufficiently small $\eps_1 >0$. 
For all $|u| \leq N^{-\alpha}$,  
\begin{align}
\log \langle \exp \left[ \beta^{1/2} t N^{-1/6} v^T \sigma \right] \rangle &= \frac{N}{2} \left( G_u ( \gamu ) - G(\gamma) \right) - \frac{1}{2} \log G_u'' ( \gamu ) + \frac{1}{2} \log G'' ( \gamma) + \O (N^{-c_1} ).
\end{align}
\eel

We now turn to estimating the quantities appearing in the above lemma, encapsulated in the following.
\bel  There is a small $c_1 >0$ so that the following holds on $\Feps$ for sufficiently small $\eps_1 >0$.
 For all $|t| \leq C$, we have
\beq
\log \langle \exp \left[ \beta^{1/2} t N^{-1/2} v^T \sigma \right] \rangle = - \frac{1}{2} (t^2  )  N^{-1} v^T (M - \gamma)^{-1} v  - \theta^{1/2} t N^{1/3}  N^{-1}v^T(M- \gamma)^{-1} v + \O (N^{-c_1} ).
\eeq
\eel
\proof   First we have by Taylor expansion that on $\Feps$,
\beq
\gamu - \gamma = \frac{u}{ G'' ( \gamma ) } N^{-4/3} v^T (M- \gamma)^{-2} v + N^{\eps_1/5} \O ( u^2 N^{-2/3} ),
\eeq
where we used the lower bounds $\gamma, \gamu \geq \lambda_1 (M)+N^{-2/3-\eps_1/100}$ which hold on $\Feps$ as well as the delocalization and rigidity estimates.   Furthermore we have on $\Feps$ via Taylor expansion,
\begin{align}
N( G_u ( \gamu ) - G ( \gamma)) =& \frac{1}{2} N G'' ( \gamma)( \gamu - \gamma)^2 \nonumber\\
-& u N^{-1/3} v^T (M- \gamma)^{-1} v - u ( \gamu - \gamma) N^{-1/3} v^T (M-\gamma)^{-2} v \nonumber \\
 +& N^{\eps_1} \O \left( N^{1/3} |u|^3 \right) \nonumber \\
=& - \frac{1}{2} (N^{1/3} u^2 ) \frac{ (N^{-4/3} v^T (M- \gamma)^{-2} v )^2}{ N^{-2/3} G'' ( \gamma) } \nonumber\\
-& u N^{-1/3} v^T (M- \gamma)^{-1} v  + N^{ \eps_1} \O \left( N^{1/3} |u|^3 \right).
\end{align}
Note also that
\beq
|G''_u ( \gamu) - G'' ( \gamma) | \leq N^{2/3+\eps_1/10} |u|,
\eeq
on $\Feps$.

Hence we see that on $\Feps$
\begin{align}
\log \langle \exp \left[ \beta^{1/2} t N^{-1/6} v^T \sigma \right] \rangle &= - \frac{1}{4} (N^{1/3} u^2 ) \frac{ (N^{-4/3} v^T (M- \gamma)^{-2} v )^2}{ N^{-2/3} G'' ( \gamma) }\nonumber \\
&- \frac{1}{2} u N^{-1/3} v^T (M- \gamma)^{-1} v  + N^{\eps_1} \O \left( N^{1/3} |u|^3 + |u| \right).
\end{align}
Now, subbing in \eqref{eqn:qa-6}, we have on  $\Feps$,
\begin{align}
& \frac{1}{4} (N^{1/3} u^2 ) \frac{ (N^{-4/3} v^T (M- \gamma)^{-2} v )^2}{ N^{-2/3} G'' ( \gamma) } - \frac{1}{2} u N^{-1/3}   v^T (M- \gamma)^{-1} v \nonumber \\
=& - \frac{1}{2} (t^2 N^{2/3} ) N^{-1} v^T (M - \gamma)^{-1} v  - \theta^{1/2} t N^{-1/3} v^T(M- \gamma)^{-1} v +\O ( N^{\eps_1 } N^{1/3} t^2 ).
\end{align}
We have therefore proven the  lemma. \qed

For small $\eps_2 >0$ we have, on the event $\Feps$ for $\eps_1 < \eps_2$,
\begin{align}
\frac{1}{N} v^T (M- \gamma)^{-1} v &= \frac{1}{N} v^T (M - \gamma - N^{\eps_2-2/3} )^{-1} v + \O \left( N^{\eps_2+\eps_1-1/3} \right) \nonumber\\
&=\msc( \gamma + N^{\eps_2 - 2/3} ) +  \O \left( N^{\eps_1+\eps_2-1/3} \right) \nonumber \\
&= -1 + \O (N^{-c} ),
\end{align}
for some small $c >0$.  Note that in passing to the second line we used the fact that the isotropic estimates of Theorem \ref{thm:iso} are assumed to hold on $\Feps$.  
This together with the previous lemma proves Theorem \ref{thm:mr-int-ext}.  The equation \eqref{eqn:exp-exp-yy1} of Theorem \ref{thm:mr-int-ext-exp} follows similarly to the proof of Lemma \ref{thm:mr-gauss-over-b1}; the rest of the theorem follows from Theorem \ref{thm:int-conv-2}.

\subsection{Overlap between two replicas}\label{sec:int-over}

In this section we will prove Theorem \ref{thm:mr-int-over}.   
We fix $\alpha >0$ and assume that
\beq
|t| \leq N^{-\alpha}.
\eeq
According to the representation formulas of Proposition \ref{prop:rep},
\beq
\left\langle \exp \left[ \beta N^{-2/3} t \sigma^{(1)} \cdot \sigma^{(2)} \right] \right\rangle = \frac{ \int_{\Gamma^2} \exp \frac{N}{2} G_1 (z, w) \d z \d w }{\left( \int_\Gamma \exp \frac{N}{2} G (z) \d z \right)^2}
\eeq
where,
\beq
G_1 (z, w) = \beta(z+w) - \frac{1}{N} \sum_{i=1}^N \log ((z- \lambda_i)(w - \lambda_i )- N^{-4/3} t^2 ) - \frac{\theta}{N^{4/3}} \sum_i v_i^2 \frac{2\lambda_i - z - w - 2 tN^{-2/3}}{ (\lambda_i - w )( \lambda_i - z ) - N^{-4/3} t^2 }.
\eeq
We define $\gam1$ to be the largest solution to
$
\partial_z G_1 ( \gam1 , \gam1) = 0.
$
We recall $\gamma$ as the solution to $G' ( \gamma ) = 0$.  For $t = o (1)$ we will develop an expansion of $\gam1$ in terms of $\gamma$ and $t$. 
Calculating the derivative, we see that this is the solution to the equation
\beq
\beta + \frac{1}{N} \sum_{i=1}^N \frac{ \lambda_i - \gam1}{ ( ( \lambda_i - \gam1 )^2 - N^{-4/3}t^2 } = \frac{ \theta}{ N^{4/3}} \sum_{i=1}^N v_i^2 \frac{ 1}{ ( \lambda_i - \gam1 + N^{-2/3}t)^2}.
\eeq
The limit of the LHS as $\gam1 \to \lambda_1 +N^{-2/3} |t|$ from above is $-\infty$ whereas the limit on the RHS is either $+ \infty$ or a positive number.  As $\gam1 \to \infty$, the LHS goes to $\beta$ and the RHS goes to $0$, so $\gam1 > \lambda_1 +N^{-2/3} |t|$.  For $\gam1> \lambda_1 +N^{-2/3} |t|$, both sides of the equation are monotonic functions of $\gam1$.  Furthermore, the LHS is less than $\beta$ for such $\gam1$ and the RHS is larger than
\beq
\frac{ \theta}{ N^{4/3}} \sum_{i=1}^N v_i^2 \frac{ 1}{ ( \lambda_i - \gam1 + N^{-2/3}t)^2} \geq \frac{ \theta v_1^2}{(N^{2/3}(\lambda_1 - \gam1 + N^{-2/3} t ))^2}
\eeq
Hence, on the event $\Feps$ we have that, $ N^{2/3} (\gam1  - ( \lambda_1 + N^{-2/3} t ) ) \geq c N^{-\eps_1/200}$.  An upper bound can be proven via a similar argument to Lemma \ref{lem:int-gams-comp}.  Hence, we have the following.
\bel \label{lem:int-gam1-est}
Assume that $|t| \leq N^{-\alpha}$.  For all $\eps_1$ sufficiently small, we have on $\Feps$ that
\beq
N^{-\eps_1/10} \leq N^{2/3} ( \gam1 - \lambda_1 ) \leq N^{\eps_1 /100}.
\eeq
\eel
Using the above lower bound above and the rigidity estimates, we see that on $\Feps$,
\beq
\left| 1 + \frac{1}{N} \sum_{i=1}^N \frac{ \lambda_i - \gam1}{ ( ( \lambda_i - \gam1 )^2 - N^{-4/3} t^2 } \right| \leq N^{-1/3+\eps_1/10}.
\eeq
Therefore, in a similar fashion to the second estimates of Lemma \ref{lem:int-gams-comp}, we see that on $\Feps$,
\beq \label{eqn:int-over-a2}
| \gamma - ( \gam1 - N^{-2/3}t ) | \leq N^{-1+\eps_1}.
\eeq
We next develop some estimates for the Hessian of $G_1$.  The proof is based on an explicit calculation and is so deferred to Appendix \ref{a:G-hess}.
\bel \label{lem:G-hess} 
On the event $\Feps$ we have that the eigenvalues of $\nabla^2 G_1 ( \gam1, \gam1)$ are both positive and between $N^{-2/3 - \eps_1/2}$ and $N^{-2/3+\eps_1/2}$. 
\eel

Before we begin the steepest descent analysis, we require the following preliminary lemma.   
\bel
Let $0.01 > \alpha >0$.  Let $G_1(z, w)$ be as above.  There is an $\eps_1 >0$ depending on $\alpha$ so on the event $\Feps$ we have the following.  Assume $|t| \leq N^{-\alpha}$ and $z = \gam1 + \i u$ and $w = \gam1+ \i v$.  Assume $| u | \geq N^{-2/3-\alpha/100}$.  Then the following estimate holds.
\begin{align}
\Re[ G_1 (z, w) - G_1 ( \gam1, \gam1 ) ]  \leq - N^{-2/3-1/10}- c \1_{ \{ |u| > N^{1/100} \} } \log( |u| / N^{1/100} ).
\end{align}
\eel
\proof First when $|v| \leq N^{-2/3- \alpha/100}$ and $u = N^{-2/3 -\alpha/100}$ we see that for $\eps_1 > 0$ sufficiently small  that on $\Feps$, we have by Taylor expansion and the estimates of Lemma \ref{lem:G-hess},
\beq
\Re[ G_1 (z, w) - G_1 ( \gam1, \gam1 ) ] \leq - N^{2/3-\eps_1}  (u^2 + v^2 ) + \O (N^{4/3+\eps_1} (|u|^3+|v|^3 ) \leq - N^{-2/3-1/10},
\eeq
for $N$ sufficiently large. 
 By direct calculation we see that on $\Feps$, for any $u, v \in \rr$,
\begin{align}
\partial_z (G_1 (z, w) &= \partial_z G(z) + N^{\eps_1} \O \left( \frac{t^2}{N N^{4/3} } \sum_i \frac{1}{ |z- \lambda_i |^2 |w- \lambda_i |} \right) \nonumber \\
&+ N^{\eps_1} \O \left(  \frac{t}{N^{6/3} } \sum_i \frac{1}{ | \lambda_i - z|^2 | \lambda_i -w |} \right). \label{eqn:int-over-a1}
\end{align}
 For $u \geq N^{-2/3-\alpha/2}$ we have $|t N^{-2/3} / ( \lambda_i -z ) | \leq N^{-\alpha/2}$.  So with $\eps_1 >0$ small enough we see that (using the rigidity estimates) all of the big $\O$ the error terms in \eqref{eqn:int-over-a1} are $\O (N^{-\alpha/4})$.  On the other hand, if $N^{-2/3-\alpha/100} \leq u \leq N^{-2/3+ \alpha/100}$,
\beq
\Im [ \partial_z G(z) ] \geq \frac{1}{N^{4/3}} \frac{ v_1^2 u (\gam1 - \lambda_1 ) }{ ( ( \lambda_1 - \gam1)^2 + u^2 )^2} \geq N^{-\alpha/10},
\eeq
on $\Feps$ as long as $\eps_1$ is sufficiently small.  Note we used \eqref{eqn:F-v-lb} and Lemma \ref{lem:int-gam1-est}.  Hence, we see that on $\Feps$ we have that $\Im [ \partial_z G_1 (z, w) ] \geq 0$ for any $v$  
 and $N^{-2/3- \alpha /100} < u < N^{-2/3 + \alpha /100}$ and so we conclude the estimate of the Lemma for this range of $u$ and $v$ satisfying $|v| \leq N^{-2/3 - \alpha/100}$.  For $u > N^{-2/3+\alpha /100}$ we see from the fact that the isotropic estimates 
 are assumed to hold on $\Feps$ as well as the rigidity estimates that,
\beq
\Im [ \partial_z G(z) ] = N^{-1/3} \Im[ \msc' ( z )] + \O ( N^{-\alpha/1000} (N^{2/3} u)^{-1})
\eeq
on $\Feps$ for $\eps_1 >0$ sufficiently small.  A straightforward calculation gives, for $u \leq 1/10$,
\beq
\Im [ \msc' (z) ] \geq \frac{ c}{ u^{1/2} }.
\eeq
as long as $| \gam1 - 2 | \leq N^{-2/3+ \alpha/10000}$.  We see also that the errors in \eqref{eqn:int-over-a1} are $\O (N^{-\alpha/4} (N^{2/3} u)^{-1})$ for any $v$.  Hence, we see that for $N^{-2/3+\alpha/10} \leq u \leq 1/10$ that $\Im [ \partial_z G_1 (z, w) ] \geq 0$ for this range of $u$ as well.  For $u > 1/10$, we see that the errors in \eqref{eqn:int-over-a1} are $\O (N^{-1/10} |u|^{-2} )$ and that $\Im [ \partial_z G(z) ] \geq c/ u$.  So we have proven the required estimates as long as $|v| \leq N^{-2/3-\alpha/100}$.

 Reversing the roles of $z$ and $w$ in the above arguments we see that on $\Feps$,
\beq
\Im [ \partial_w G_1 (z, w) ] \geq 0
\eeq
for any $u \geq N^{-2/3 - \alpha / 100}$ and $v \geq N^{-2/3-\alpha/100}$.  The arguments for $u < N^{-2/3 - \alpha /100}$ and/or $v \leq N^{-2/3 - \alpha/100}$ are almost identical. This completes the proof. \qed

The following is an easy consequence of the previous lemma and the method of steepest descent.
\bep
There is a $c_1>0$ so that the following holds.  Let $\alpha >0$ and assume $|t| \leq N^{-\alpha}$.   Then there is an $\eps_1 >0$ depending on $\alpha$ so that on $\Feps$,
\begin{align}
\log \langle \exp [ \beta N^{-2/3} t \sigma^{(1)} \cdot \sigma^{(2)} ] \rangle &= \frac{N}{2} \left( G_1 ( \gam1, \gam1 ) - 2 G ( \gamma ) \right) \nonumber \\
&+ \log \det \nabla^2 G_1 ( \gam1, \gam1 ) - 2 \log G'' ( \gamma) + \O (N^{-c_1} ).
\end{align}
\eep
We now prove the following.
\bep \label{prop:int-over-1}
There is a $c_1>0$ so that the following holds.  Let $\alpha >0$.   Then on the event $\Feps$ for $\eps_1 >0$ sufficiently small depending on $\alpha >0$ we have,
\begin{align}
&\log \langle \exp [ \beta N^{-2/3} t \sigma^{(1)} \cdot \sigma^{(2)} ] \rangle =\frac{N}{2} \left( 2  \beta tN^{-2/3} - \frac{1}{N} \sum_{i=1}^N \log ( 1 + 2 tN^{-2/3} ( \gamma - \lambda_i )^{-1} ) \right) \nonumber\\
+& \log \left( 1 - \frac{1}{ \theta N^{-1/3} \mv'' ( \gamma ) } \left[\theta N^{-1/3} \mv'' ( \gamma)  -  \frac{ 2 \theta}{N^{4/3}} \sum_i v_i^2 \frac{1}{ ( \lambda_i - \gamma - 2 t N^{-2/3} )( \lambda_i - \gamma)^2 }  \right] \right) + \O(N^{-c_1} )
\end{align}
\eep
\proof 
Now from \eqref{eqn:int-over-a2} we have from a second order Taylor expansion that on $\Feps$ for $\eps_1 >0$ sufficiently small,
\beq
G_1 ( \gam1, \gam1 )  = G_1 ( \gamma+N^{-2/3} t , \gamma+N^{-2/3} t ) + \O \left(N^{-7/6} \right).
\eeq
We have the equality,
\begin{align}
G_1 ( \gamma + tN^{-2/3}, \gamma + tN^{-2/3} ) - 2 G ( \gamma)&= 2 \beta tN^{-2/3} - \frac{1}{N} \sum_{i=1}^N \log ( 1 + 2 tN^{-2/3} ( \gamma - \lambda_i )^{-1} ).
\end{align}
Similarly, we see in a straightforward manner that on $\Feps$,
\begin{align}
\left| \log \det \nabla^2 G_1 ( \gamma+N^{-2/3}t, \gamma+N^{-2/3}t ) - \log \det \nabla^2 G_1 ( \gam1, \gam1 ) \right| \leq C N^{-1/6}.
\end{align}
We begin by calculating the determinant of the Hessian.  We have,
\begin{align}
&\det \nabla^2 G_1 ( \gamma+N^{-2/3}t, \gamma+N^{-2/3}t ) = \left(  \frac{1}{N} \sum_{i=1}^N  \frac{1}{ ( \lambda_i - \gamma)^2 - 2 N^{-2/3}t ( \lambda_i - \gamma ) } - \frac{2 \theta}{N^{4/3}} \sum_i v_i^2 \frac{1}{ ( \lambda_i - \gamma)^3}  \right) \nonumber \\
\times & \left(  \frac{1}{N} \sum_i \frac{ ( \lambda_i - \gamma)^2 + 2 N^{-4/3} t^2 - 2 tN^{-2/3} ( \lambda_i - \gamma) }{ ( ( \lambda_i - \gamma)^2 - 2tN^{-2/3} ( \lambda_i - \gamma) )^2}   - 2 \frac{ \theta}{N^{4/3}} \sum_i v_i^2 \frac{1}{ ( \lambda_i - \gamma - 2 N^{-2/3} t)( \lambda_i - \gamma)^2 }\right).
\end{align}
For the first factor, we have on $\Feps$,
\begin{align}
 \frac{1}{N} \sum_{i=1}^N  \frac{1}{ ( \lambda_i - \gamma)^2 - 2 N^{-2/3}t ( \lambda_i - \gamma ) } - \frac{2 \theta}{N^{4/3}} \sum_i v_i^2 \frac{1}{ ( \lambda_i - \gamma)^3} 
=  G'' ( \gamma )(1  + \O(N^{-1/3+\eps_1/10} ) ).
\end{align}
Recall our notation,
\beq
\mv(z) := \frac{1}{N} v^T \frac{1}{M -z } v
\eeq
so that
\beq
G'' (\gamma) = - \theta N^{-1/3} \mv'' (\gamma) + \O (N^{1/3+\eps_1/10 } ).
\eeq
Then,
\begin{align}
& \frac{1}{N} \sum_i \frac{ ( \lambda_i - \gamma)^2 + 2 N^{-4/3} t^2 - 2 tN^{-2/3} ( \lambda_i - \gamma) }{ ( ( \lambda_i - \gamma)^2 - 2tN^{-2/3} ( \lambda_i - \gamma) )^2}   - 2 \frac{ \theta}{N^{4/3}} \sum_i v_i^2 \frac{1}{ ( \lambda_i - \gamma - 2 N^{-2/3} t)( \lambda_i - \gamma)^2 }\nonumber \\
=~& G'' ( \gamma)+ \left[ \theta N^{-1/3} \mv'' ( \gamma)  -  \frac{ \theta}{N^{4/3}} \sum_i v_i^2 \frac{1}{ ( \lambda_i - \gamma - 2 t N^{-2/3} )( \lambda_i - \gamma)^2 } \right]+ \O (N^{1/3+\eps_1/10}).
\end{align}
This proves the claim. \qed

From Proposition \ref{prop:int-over-1} we see that there is a small $c_1>0$ so that for every $\alpha>0$ and $k>0$ we have for $\eps_1>0$ sufficiently small depending on $\alpha, k$ that on $\Feps$,
\begin{align} \label{eqn:int-over-gen}
 \log \langle \exp \left[ \beta N^{-2/3} t \sigma^{(1)} \cdot \sigma^{(2)} \right] \rangle = & t \left( \beta N^{1/3}  + \frac{1}{N^{2/3}} \sum_{i=1}^N \frac{2}{ \lambda_i - \gamma } - \frac{2 \theta}{3} \frac{ N^{-5/3} \mv''' ( \gamma)}{ \theta N^{-1} \mv'' ( \gamma) } \right) \nonumber \\
+& \sum_{j=2}^k t^j \left( X_j + Y_j \right) + \O \left( |t|^{k+1} N^{2 k \eps_1} +N^{-c_1} \right),
\end{align}
where,
\begin{align}
X_j &:= \frac{1}{N^{2j/3}}\frac{2^j}{j!} N m^{(j-1)}(\gamma) \nonumber \\
Y_j &:= \frac{1}{j!} \frac{d^j}{ds^j} \log \left( 1 - \frac{1}{ \theta N^{-1/3} \mv'' ( \gamma ) } \left[\theta N^{-1/3} \mv'' ( \gamma)  -  \frac{ 2 \theta}{N^{4/3}} \sum_i v_i^2 \frac{1}{ ( \lambda_i - \gamma - 2 s N^{-2/3} )( \lambda_i - \gamma)^2 }  \right] \right) \bigg\vert_{s=0}.
\end{align}
Now, $Y_j$ is a sum of terms which each are products of the form,
\beq
\frac{ N^{1-2(l+1)/3} \mv^{(l)} ( \gamma)}{ N^{-1} \mv'' ( \gamma ) },
\eeq
for $l \geq 2$. 
So for any $\eps >0$ there is a $\eps_2 >0$ so that each term $X_j$ and $Y_j$ is $\O(N^{\eps})$, with probability at least $1- N^{-\eps_2}$.  Note that 
\beq \label{eqn:mvdenom}
N^{-\eps} \leq N^{-1} \mv'' ( \gamma )  \leq N^{\eps}
\eeq
with the same probability with probability at least $1-N^{\eps/100}$.  We recall the definition of $H$, the GOE matrix associated to $M$ via $H = M + V$, and $V$ is a diagonal matrix of iid centered normal random variables with variance $2/N$, as well as the definition of $\hatgam_H$ in \eqref{eqn:int-gam-def}.  From Lemma \ref{lem:int-gams-comp} and the estimates in its proof, we have the following.
\bel
There is a small $c_1 >0$ so that the following holds on $\Feps$ for all sufficiently small $\eps_1 >0$.  First,
\beq
| \gamma - \hatgam_H | \leq N^{-2/3-c_1}.
\eeq 
Second, for any
 $k$ and $\eps >0$ we have, for any $N^{-2/3-\eps} \leq E - \lambda_1 (M) \leq N^{-2/3+\eps}$ 
\begin{align}
\frac{1}{N^{2k/3}} \left( \left| \tr \frac{1}{ (M-E)^k} - \tr \frac{1}{(H-E)^k} \right| + \left| v^T(M-E)^{-k}v - v^T (H-E)^{-k} v \right| \right) \leq N^{-c_1+2k \eps}.
\end{align}
\eel
With this lemma we see that we can replace all the quantities in the coefficients in our Taylor series \eqref{eqn:int-over-gen} by the corresponding quantities involving $H$ and $\hatgam_H$, at a cost of $\O (N^{-c_1/2})$ on the event $\Feps$ where now $\eps_1>0$ must be taken small enough depending on the order of the Taylor expansion.  Note that the quantities in the denominator of the $Y_j$'s are only products of $N^{-1} \mv'' ( \gamma)$ for which we have a lower bound  \eqref{eqn:mvdenom}, so we do not have any small denominator difficulties.

From this observation and Theorem \ref{thm:int-conv-2} we conclude Theorem \ref{thm:mr-int-over}.

\section{$h= \O (N^{-1/2} )$ and $\beta >1 $} \label{sec:micro}

We now consider the regime of low temperature and very small magnetic field.  We assume that,
\beq
C \geq \beta \geq 1 + c, \qquad h^2 \beta = N^{-1} \theta,
\eeq
for a fixed $\theta \in \rr$.  In this section we prove the results of Section \ref{sec:mr-micro}. 

\subsection{Free energy} \label{sec:micro-fe}
In this section we will examine the free energy, and prove Theorem \ref{thm:mr-micro-fe}.  We will also develop estimates which will be used to study the overlap in the next subsection. 
Eventually we will expand $G(z)$ around a saddle $\gamma$, and apply the method of steepest descent.   The saddle will turn out to be distance $\O(N^{-1})$ from $\lambda_1$. 

For our analysis we will again use the event $\Feps$ that was defined in Definition \ref{def:int-F}.  Fix a small $\kapp>0$ satisfying at least $\kapp < \frac{1}{100}$.  The portion of the steepest descent contour that will contribute to the integral lies in the region 
\beq
\Dk :=\{ z:  N^{-1+\kapp} + \lambda_1 > \Re[z] > \lambda_1 - N^{-1+\kapp}, \ \Im[z] \leq N^{-1+2 \kapp} \}.
\eeq
The portion not lying in $\Dk$ will have an exponentially small contribution.  We first analyze the behavior of $G(z)$ inside $\Dk$.  
For $z \in \Dk$ we have the following expansion that holds on $\Feps$, using the level repulsion, rigidity and delocalization estimates:
\beq
G' (z) = (\beta -1 ) + \frac{1}{N ( \lambda_1 - z)} - \frac{ \theta v_1^2}{ N^2 ( \lambda_1 - z)^2} + \O (N^{-1/3+\eps_1/10} ).
\eeq
This prompts us to define $\gamt = \lambda_1 + \cb/N$ as the solution of the equation resulting from setting the first three terms on the RHS to $0$. This leads to a quadratic equation for $\cb$,
\beq
\beta-1 = \frac{1}{ \cb} + \frac{ \theta v_1^2}{ \cb^2},
\eeq
which has a positive solution,
\beq
\cb = \frac{ 1 + \sqrt{1 + 4 ( \beta-1) \theta v_1^2}}{2 ( \beta-1)}.
\eeq
Note the estimate
\beq
\cb \asymp 1 + |v_1|.
\eeq
Note that our notation has changed slightly compared to previous sections: $\gamt$ is not defined as a solution to the equation $G' (z) = 0$.  Instead $\gamt$ solves an approximation to the saddle point equation which, while random, is more tractable than the full equation.  Note that $\cb$ is random through its dependence on the quantity $v_1^2$. We use the same notation $\cb$ as in the paper \cite{SSK-ls}, where $\theta = 0$ (so that in that paper $\cb$ was deterministic).

Recall the notation in \eqref{eqn:m-tilde-def-f} for the quantities $\tilm (z)$ and $\tilmv (z)$ which 
 separate the contribution of the largest eigenvalue from $m$ and $m_v$.  For $|z| \leq N^{-1+\kapp}$ the following expansion holds on $\Feps$ for sufficiently small $\eps_1 >0$,
\begin{align} \label{eqn:aug30-2}
G (z + \gamt) - G ( \gamt ) &=(\beta + \tilm_N ( \gamt ) ) z - \frac{1}{N} \log(1 + Nz / \cb ) + \frac{z^2}{2} \tilm_N' ( \gamt) \nonumber \\
&+ \frac{v_1^2 \theta}{N \cb} \left( \frac{1}{ 1 + N z / \cb} -1 \right) - \frac{z \theta}{N} \tilmv'  ( \gamt) + \O (N^{-2+3 \kapp+\eps_1/10} ).
\end{align}
Making the change of variable $u = Nz / \cb$ and multiplying the above function by $N$ we see that we should examine the function
\beq
f(u) := \cb ( \beta-1) u - \log(1+u) + \frac{\theta v_1^2}{ \cb} \left( \frac{1}{1+u} -1 \right).
\eeq
Define $B$ by,
\beq
B+1 = \cb (\beta-1) = \frac{ \theta v_1^2}{ \cb} + 1
\eeq
so that, (note $B>0$)
\beq
f(u) = (B+1) u - \log(1+u) + B \left( \frac{1}{1+u} - 1 \right).
\eeq
We denote $u = E + \i \eta$ and look for solutions to 
\beq \label{eqn:imp-eta}
\eta \left( B+1 - \frac{B}{ (1+E)^2 + \eta^2 } \right) = \arg \left( (1 + E ) + \i \eta \right)
\eeq
for $E <0$.  Note that this is the steepest descent contour for $f(z)$, i.e. $\Im [f(z) ] = 0$.   We will consider $\eta$ as a function of $E$ and look to find estimates on $\etat(E)$.   We first want to check that $\etat(E)$ is in fact well-defined.  Accordingly, we first provide the following lemma which shows that $\etat(E)$ is well-defined and derives some basic properties. It is based on calculus so is deferred to Appendix \ref{a:micro-1-f}.  
\bel \label{lem:micro-1-f}
For $E <0$ we define the function $\etat(E)$ as follows.  First, $\etat(0) = 0$ and,
\beq
\etat (-1) = \frac{ \frac{\pi}{2} + \sqrt{\frac{ \pi^2}{4} + 4B(B+1)}}{2 (B+1)}.
\eeq
For $-1 < E < 0$ there is a unique solution to the equation 
\beq \label{eqn:imp-eta2}
\eta \left( B+1 - \frac{B}{ (1+E)^2 + \eta^2 } \right) = \arg \left( (1 + E ) + \i \eta \right)
\eeq
satisfying $0  < \eta < \etat(-1)$.  We define $\etat(E)$ to be this solution.  Then $\etat(E)$ is strictly decreasing on $-1 < E < 0$ and is smooth.  For $E <-1$ there is a unique solution of \eqref{eqn:imp-eta2} on the domain $\eta > \eta_* (E)$ where 
\beq
\eta_* (E) := \inf\left\{ \eta>0 : (1+E)^2 + \eta^2 > B(1+B)^{-1} \right\}
\eeq
The function $\etat(E)$ is smooth for $E < -1$.
\eel

We require some estimates on the steepest descent contour which are summarized in the following. The proof is deferred to Appendix \ref{a:micro-2-f}. 
\bel \label{lem:micro-2-f}
There is a $c>0$ so that if $-c < E < 0$ we have,
\beq
 \etat (E) =  \sqrt{ -E \frac{3+6B}{1+3B} } \left( 1+ \O(|E| ) \right) .
\eeq
For any $c>0$ there is a $c'>0$ so that $\etat (E) \geq c'$ if $-1 < E < -c$.  There is a $C>0$ so that $\etat (E) \leq C$ for all $E<0$.  Finally for $E<-1$ we have,
\beq
\etat (E) \geq \frac{ \pi}{2 (B+1)}.
\eeq
Note that the constants above do not depend on $B$. 
\eel

In order to control the function $G$ on the steepest descent contour we establish the following estimates on $\Re[f']$. The proof is deferred to Appendix \ref{a:micro-f-est-1}. 
\bel \label{lem:micro-f-est-1} 
Let $c>0$.  Then there is a $c'>0$ and $C'>0$ so that for $-1 < E < -c$ we have,
\beq
\Re[f'(E+ \i \eta) ] > c'(1+B)
\eeq
and 
\beq
\Re[ f' (E + \i \eta ) ] > c'
\eeq
for $E < -1$.  We have $|\Im [f' (E+ \i \eta ) ] \leq C$ as long as $E < -c$.
\eel

Now we define the contour $\Gamk$ by,
\beq
\Gamk := \{ z = (E+ \i \etat (E) ) \cb N^{-1} : 0 \leq -E \leq N^{\kapp} \}.
\eeq
Finally, with our estimates on $f$ established, we can obtain the following estimates on $G$ along the steepest descent contour.
\bel \label{lem:micro-G}  Let $\kapp >0$.  For sufficiently small $\eps_1 >0$ the following holds on $\Feps$ for $N$ large enough.  
For any $c>0$ there is a $c'>0$ we have for $z \in \Gamk$,
\beq
N\Re[ G (z + \gamt) - G ( \gamt ) ] \leq - N \cb^{-1} c' (\Re[z] - c \cb N^{-1} )  \1_{ \{ \Re[z] \leq - c \cb N^{-1} \}} + N^{-1/3+\eps_1+ \kapp}.
\eeq
Moreover, 
there are $C>0$ and $c''>0$ so that on the line $z = -N^{\kapp} + \i \eta$ we have the estimate,
\beq
N( \Re[ G (z + \gamt) - G ( \gamt ) ])  \leq - c'' N^{\kapp} \cb^{-1} - \1_{ \{ \eta > C+1 \} } N c'' \log ( \eta).
\eeq
For any $|E| \leq 5$ and all $\eta$ sufficiently large we have,
\begin{align}
N( \Re[ G (z + \gamt) - G ( \gamt ) ])  \leq   - N c'' \log ( \eta).
\end{align}
\eel
\proof
Starting from \eqref{eqn:aug30-2} and applying rigidity and level repulsion, we see that
 for $|u| \leq N^{ \kapp}$,
\beq
G ( u\cb /N + \gamt) - G ( \gamt) = \frac{1}{N} f(u) + \O (N^{-1/3-1+\kapp+ \eps_1/10} ),
\eeq
on $\Feps$. 
Recall the equations,
\beq
\frac{\d \eta}{ \d E} = - \frac{ \Im [ f']}{ \Re[f']}, \qquad \frac{\d}{\d E} \Re[f (E, \etat (E) )] = \Re[f'] + \frac{ \Im[f']^2}{ \Re[f']}.
\eeq
Hence, by integrating the estimates of Lemma \ref{lem:micro-f-est-1} we find that for $-N^{ \kapp} < E < 0$ and any $c>0$ there is a $c'>0$ such that on the event $\Feps$,
\beq
N\Re[ G ( (E+\i \etat(E) ) \cb /N + \gamt) - G ( \gamt)] \leq - c' (|E|-c) \1_{ \{E \leq -c\} } + N^{-1/3+\kapp+\eps_1}.
\eeq
This proves the first part of the lemma.  We now control $G$ on the vertical line $E = -N^{\kapp}$ and $\eta \geq \etat ( - N^{\kapp})$.  We calculate,
\beq
\frac{\d}{ \d \eta} \Re[f ( E+ \i \eta) ] = - \Im [ f'] = - \frac{ \eta}{ (1+E)^2 + \eta^2 }\left(1 + \frac{B(1+E)}{ (1+E)^2 + \eta^2} \right)
\eeq
As long as $|B| \leq N^{\kapp/2}$, which we can guarantee by taking $\eps_1 >0$ sufficiently small, then the fraction in the brackets above is $o(1)$ and so in this case $-\Im[f'] <0$.  Therefore, for $N^{\kapp} > \eta > \etat (-N^{\kapp})$ we have,
\beq
N \Re[ G ( (-N^{\kapp} + \i \eta ) \cb/N + \gamma ) - G ( \gamma ) ] \leq - c N^{\kapp} \label{eqn:vertG}
\eeq
for some $c>0$.  Note that we also clearly have $\Im[f'] \leq C$ when $|B| \leq N^{\kapp/2}$ and since $\etat (-N^{\kapp} ) \leq C$, the estimate \eqref{eqn:vertG} holds also for $0 \leq \eta \leq \etat(-N^{\kapp})$.

Now, turn to $\eta \geq N^{\kapp}$.  We calculate,
\begin{align}
\frac{\d}{ \d \eta} \Re[ G ( E + \i \eta + \gamt )] &= - \Im [ m (z) ] - \theta N^{-2} \Im [ v^T (M-z)^{-2} v ] \nonumber  \\
&\leq - \frac{1}{N} \sum_{i=1}^N \frac{ \eta}{ | \lambda_i - z|^2} - \frac{ |v_i|^2 N^{-1} \theta}{ | \lambda_i - z|^2}.
\end{align}
Hence as long as $\eta \geq N^{-1+\kapp}$ for any $\kapp>0$ this will be negative as long as $|v_i|^2 \leq N^{\kapp/2}$ which holds on $\Feps$ by our choice of $\eps_1 >0$.  Hence the estimate \eqref{eqn:vertG} holds for any $\eta \geq \etat(-N^{\kapp})$.

Now for any $|E| \leq 5$ and $\eta \geq 10$ we have,
\begin{align}
\Im[ G' (E + \i \eta ) ] \geq \Im [m (E + \i \eta ) ]  - \frac{C}{ \eta^2} \geq \frac{c}{ \eta} - \frac{C} {\eta^2}
\end{align}
where in the first inequality we used the fact that $||v||_2^2 = N$.   This completes the proof.  \qed

Note that during the proof we established the following.
\bel \label{lem:micro-f-est-2}  For any $c>0$ there is a $c' >0$ so that 
for $z \in \Gamk$,
\beq
\Re[ f( z N / \cb )] \leq  - N \cb^{-1} c' ( \Re[z] - c \cb N^{-1} ) \1_{ \{ \Re[z] \leq - c \cb N^{-1} \} }.
\eeq
\eel

We denote by $\tilGamk$ the contour
\beq
\tilGamk = \Gamk \cup \{ z = \cb N^{-1} (- N^{\kapp} + \i \eta), \eta < \etat (-N^{\kapp} ) \}.
\eeq
We can now calculate the free energy.  Using our representation formula we have,
\begin{align}
\frac{1}{N} \log Z_{N, \beta, h} &= \frac{1}{N} \left( (1-N/2) \log (\beta) +\frac{N}{2} \log(2 \pi) + \frac{1}{2} \log(N) \right) \nonumber \\
&+ \log \frac{1}{ 2 \pi \i} \int_{ \gamt - \i \infty}^{ \gamt+\i \infty} \exp \left[ \frac{N}{2} G(z) \right] \d z,
\end{align}
where $\gamt$ is as above.  In the following we apply the method of steepest descent to calculate the above contour integral.
\bep \label{prop:micro-int}  Fix $\kapp>0$ sufficiently small.  Then for all sufficiently small $\eps_1 >0$ it holds on $\Feps$ that,
\begin{align}
 & \int_{ \gamt - \i \infty}^{ \gamt+\i \infty} \exp \left[ \frac{N}{2} (G(z)-G( \gamt) ) \right] \d z  \nonumber\\
= & \int_{\tilGam}  \exp \left[ \frac{\beta+\tilm_N ( \gamt)}{2} \cb u - \frac{1}{2} \log(u) + \frac{v_1^2 \theta}{ 2 \cb} \frac{1}{u} \right] \times\nonumber \\
& \exp \left[ - \cb ( \beta + \tilm_N ( \gamt ) )/2 - \frac{v_1^2 \theta}{ 2 \cb} \right] \left(1 + \frac{\cb^2}{ 4 N} (u-1)^2 \tilm_N' ( \gamt)  - \theta \cb \frac{u-1}{ 2 N} \tilmv' ( \gamt  ) \right) \frac{ \d u \cb}{N} + \O(N^{-2+4 \kapp+\eps_1} ) \nonumber\\
=& \int_{\tilGam}  \exp \left[ \frac{\beta+\tilm_N ( \gamt)}{2} \cb u - \frac{1}{2} \log(u) + \frac{v_1^2 \theta}{ 2 \cb} \frac{1}{u}  - \cb ( \beta + \tilm_N ( \gamt ) )/2 - \frac{v_1^2 \theta}{ 2 \cb} \right]  \frac{ \d u \cb}{N} + \O(N^{-5/3+4 \kapp+\eps_1} ) 
\end{align}
where $\tilGam$ is a keyhole contour  circling $u=0$ and continuing above and below the real axis to $- \infty$.
\eep
\proof Using Lemma \ref{lem:micro-G} we can first replace the contour consisting of a vertical line in the complex plane by $\tilGamk$ at the cost of an exponential error,
\begin{align}
\int_{ \gamt - \i \infty}^{ \gamt+\i \infty} \exp \left[ \frac{N}{2} (G(z)-G( \gamt) ) \right] \d z = \int_{\tilGamk} \exp \left[ \frac{N}{2} ( G( z + \gamt) - G ( \gamt) ) \right] \d z + \O (\e^{-N^{c} }),
\end{align}
on the event $\Feps$ for some $c>0$.
   Via Taylor expansions we have on $\Feps$, (see \eqref{eqn:aug30-2})
\begin{align}
 &\int_{\tilGamk} \exp \left[ \frac{N}{2} ( G( z + \gamt) - G ( \gamt) ) \right] \d z \nonumber \\
 =& \int_{\tilGamk} \exp \left[ \frac{N}{2}( \beta + \tilm_N ( \gamt ) )z - \frac{1}{2} \log(1 + Nz / \cb ) + \frac{N z^2}{4} \tilm_N' ( \gamt) + \frac{v_1^2 \theta}{ 2 \cb} \left( \frac{1}{ 1 +N z / \cb } - 1 \right) - \frac{z \theta}{2} \tilmv' ( \gamt ) \right] \d z \nonumber \\
+ &\O ( N^{-2 + 4 \kapp+\eps_1/5} ) \nonumber \\
=& \int_{ \tilGamk} \exp \left[ \frac{N}{2} ( \beta + \tilm_N ( \gamt ) ) z  - \frac{1}{2} \log\left( 1 + \frac{Nz}{ \cb} \right) + \frac{v_1^2 \theta}{2 \cb } \left( \frac{1}{ 1 + \frac{Nz}{\cb} } - 1 \right) \right] (1 + \frac{Nz^2}{4} \tilm_N' ( \gamt) - \frac{z \theta}{2} \tilmv' ( \gamt ) ) \d z \nonumber \\
+& \O (N^{-2 +4 \kapp+\eps_1/5} ) \label{eqn:micro-int-1}
\end{align}
Now for $E \leq - N^{\kapp}/N$ we see that
\beq
\Re \left\{ \frac{N}{2} ( \beta + \tilm_N ( \gamt ) ) z  - \frac{1}{2} \log\left( 1 + \frac{Nz}{ \cb} \right) + \frac{v_1^2 \theta}{2 \cb } \left( \frac{1}{ 1 + \frac{Nz}{\cb} } - 1 \right) \right\}  \leq - c N|E|
\eeq
and so if we define $\Gamma$ as a keyhole contour around $z = -\cb/N$, we see we can replace $\tilGamk$ by $\Gamma$ at an error of $\O (\e^{-N^{\kapp/2}})$.  Making the change of variable $u= 1 + Nz / \cb$ (so that $\tilGam$ is the image of $\Gamma$ under this change of variable) we see that,
\begin{align}
 & \int_{ \gamt - \i \infty}^{ \gamt+\i \infty} \exp \left[ \frac{N}{2} (G(z)-G( \gamt) ) \right] \d z\nonumber \\
= & \int_{\tilGam}  \exp \left[ \frac{\beta+\tilm_N ( \gamt)}{2} \cb u - \frac{1}{2} \log(u) + \frac{v_1^2 \theta}{ 2 \cb} \frac{1}{u} \right] \times\nonumber \\
& \exp \left[ - \cb ( \beta + \tilm_N ( \gamt ) )/2 - \frac{v_1^2 \theta}{ 2 \cb} \right] \left(1 + \frac{\cb^2}{ 4 N} (u-1)^2 \tilm_N' ( \gamt)  - \theta \cb \frac{u-1}{ 2 N} \tilmv' ( \gamt  ) \right) \frac{ \d u \cb}{N} + \O(N^{-2+4 \kapp+\eps_1} ).
\end{align}
This yields the first claim. The second is obtained in a similar fashion, but by dropping the polynomial terms in \eqref{eqn:micro-int-1}. \qed

We need the following representation for Bessel functions. The elementary proof is deferred to Appendix \ref{a:bessel}.
\bel \label{lem:bessel}
Let $\tilGam$ be a keyhole contour encircling $z = 0$ and continuing to $-\infty$.  Then for $a>0, b>0$ and $\alpha \in \rr$,
\begin{align}
\frac{1}{ 2 \pi \i} \int_{\tilGam}\exp\left[ a z + b z^{-1} -  \alpha \log(z) \right] \d z \left( \frac{a}{b}\right)^{(1-\alpha)/2} &= \frac{\sin ( \alpha \pi )}{ \pi} \int_0^\infty \exp \left[ - \lambda \cosh(u) + (1- \alpha) u  \right]  \d u \nonumber \\
&+\frac{1}{ \pi} \int_0^\pi \exp \left[ \lambda \cos ( \theta ) \right] \cos ((1- \alpha) \theta ) \d \theta \nonumber \\
&= I_{\alpha-1} ( \lambda)
\end{align}
where  $\lambda = 2 \sqrt{ ab}$ and $I_\beta (x)$ denotes the modified Bessel function of the first kind.

Moreover,
\beq \label{eqn:aug31-1}
I_{1/2} (x) = \frac{\sqrt{2}}{ \sqrt{  \pi}} \frac{ \sinh(x)}{ \sqrt{x}}, \qquad I_{-1/2} (x) = \frac{\sqrt{2}}{ \sqrt{  \pi }} \frac{ \cosh(x)}{ \sqrt{x}}.
\eeq
\eel

Note that once we have the formula for one Bessel function, we can find the others an integer value of $\beta$ apart due to the recursions,
\beq
I_{\beta-1} (x) = \frac{1}{ x^\beta} \frac{\d}{ \d x} x^\beta I_{\beta} ( x), \qquad I_{\beta+1} (x) = x^\beta \frac{\d}{ \d x} x^{-\beta} I_\beta (x).
\eeq

From the explicit formula for the Bessel function $I_{-1/2}$ we have,
\bel
There are $c>0$ and $C>0$ such that for any $a>0$ and $b>0$, with $\lambda = 2 \sqrt{ab}$,
\beq
 c \frac{ \e^{ \lambda}}{ \sqrt{ \lambda} } \leq I_{-1/2} ( \lambda) \leq C \frac{ \e^{ \lambda}}{ \sqrt{ \lambda} }
\eeq
\eel

We will use all of the above with
\beq
a := ( \beta + \tilm_N ( \gamt) ) \cb /2, \qquad b:= v_1^2 \theta / (2 \cb ).
\eeq
With this notation, the leading order contribution in the second estimate of Proposition \ref{prop:micro-int} is,
\beq
\frac{1}{ 2\pi \i} \int_{\tilGam} \exp [ a u - \frac{1}{2} \log(u) + b / u  - a - b ] \cb N^{-1} \d u = b^{1/4} a^{-1/4} I_{-1/2} ( 2 \sqrt{ab} ) \exp[-a-b] \cb N^{-1}
\eeq
In the asymptotics for $I_{-1/2} ( 2 \sqrt{ ab})$  we have an exponential term $\e^{ 2 \sqrt{ab}}$ which could cause problems.  However, it is balanced by, the $\e^{-a-b}$ term as follows. First,
\beq
2\sqrt{ab} - a - b = - ( \sqrt{a} - \sqrt{b} )^2 = - \frac{ (a - b)^2}{ (\sqrt{a} + \sqrt{b} )^2} .
\eeq
Second, on the event $\Feps$,
\beq
2(a-b) = \cb (\beta -1 ) - v_1^2 \theta \cb^{-1} + \O (N^{-1/3\eps_1/10} ) = 1 + \O (N^{-1/3+\eps_1/10} ),
\eeq 
where we used the definition of $\cb$ in the second equality.  Since $a \geq c$ we therefore conclude that,  
\beq \label{eqn:micro-Z-asymp}
b^{1/4} a^{-1/4} I_{-1/2} ( 2 \sqrt{ab} ) \exp[-a-b] \asymp \frac{1}{ \sqrt{a}}
\eeq
and $N^{-\eps_1/10} \leq a^{-1/2} \leq C$.  In summary, we have on $\Feps$ that
\begin{align}
\frac{1}{N} \log Z_{N, \beta, h} &= \frac{1}{2} G( \gamt) +  \frac{1}{N} \left( (1-N/2) \log (\beta) +\frac{N}{2} \log(2 \pi) + \frac{1}{2} \log(N) \right) \nonumber \\
& + \frac{1}{N} \log\left( \cb N^{-1} (b/a)^{1/4} I_{-1/2} ( 2 \sqrt{ab} )\e^{-a-b}  \right) + \O(N^{-1-2/3+4\kapp+\eps_1} ).
\end{align}
Furthermore, since on $\Feps$,
\beq
\tilm ( \gamt) +1 = \O (N^{-1/3+\eps_1/10} ),
\eeq
we find the estimate,
\beq
\left| \log\left(  (b/a)^{1/4} I_{-1/2} ( 2 \sqrt{ab} )\e^{-a-b}  \right)  - \log\left(  (b/\hata)^{1/4} I_{-1/2} ( 2 \sqrt{\hata b} ) \e^{-\hata-b}  \right)  \right| \leq N^{-1/3+\eps_1}
\eeq
where
\beq
\hata := \frac{ \cb ( \beta -1)}{2}.
\eeq
 For $G(\gamt)$, we have on $\Feps$
\beq
G( \gamt) = (\beta-1)( \lambda_1 - 2) + 2 \beta - \int \log(2 - x) \rhosc (x) \d x + \O (N^{-1+\eps_1} ).
\eeq

These estimates are summarized in the following proposition. 
\bep For sufficiently small $\eps_1 >0$ we have on $\Feps$ that,
\begin{align}
\frac{1}{N} \log Z_{N, \beta, h} &= \frac{1}{2} G( \gamt) +  \frac{1}{N} \left( (1-N/2) \log (\beta) +\frac{N}{2} \log(2 \pi) + \frac{1}{2} \log(N) \right)  \\
& + \frac{1}{N} \log\left( \cb N^{-1} (b/\hata)^{1/4} I_{-1/2} ( 2 \sqrt{ \hata b} )\e^{-\hata-b}  \right) + \O(N^{-1-1/3+\eps_1} ),
\end{align}
as well as
\begin{align}
\frac{1}{N} \log Z_{N, \beta, h} &= \frac{\beta-1}{2} ( \lambda_1 -2) +2 \beta - \int \log(2-x) \rhosc (x) \d x  \\
& \frac{1}{N} \left( (1-N/2) \log (\beta) +\frac{N}{2} \log(2 \pi) + \frac{1}{2} \log(N) \right) + \O(N^{-1+\eps_1} ).
\end{align}
\eep
From the above, we conclude Theorem \ref{thm:mr-micro-fe}.

\subsection{Overlap between two replicas} \label{sec:micro-over}
In this section we will calculate a few moments of the overlap between two replicas $N^{-1} \sigma^{(1)} \cdot \sigma^{(2)}$.  We will build on the estimates established in the previous subsection and prove the remaining results of Section \ref{sec:mr-micro}. 

 We start with the following representation from Proposition \ref{prop:moment-rep}. 
\begin{align} \label{eqn:micro-over-rep-f}
\frac{1}{N} \langle \sigma^{(1)} \cdot \sigma^{(2)} \rangle = \frac{ \int_{\Gamma^2} \frac{ \theta}{N^2 \beta} v^T \frac{1}{ (M-w)(M-z) } v \exp\left[ \frac{N}{2} (G(z) + G(w) ) \right] \d z \d w }{ \int_{\Gamma^2} \exp \left[ \frac{N}{2} (G(z) + G(w) ) \right]  \d z \d w},
\end{align}
where $\Gamma$ denotes the vertical line in the complex plane $\gamt  + \i t$. The function $G$ is as in the previous subsection.  
We will again use the approximate steepest descent contour $\Gamk$ as defined in the previous subsection.  
It will be convenient to introduce
\beq
\tilf (z) := \frac{N}{2} ( \beta + \tilm_N ( \gamt) ) z - \frac{1}{2} \log \left( 1 + \frac{Nz}{\cb} \right) + \frac{v_1^2 \theta}{ 2 \cb} \left( \frac{1}{ 1 + \frac{Nz}{ \cb} } - 1  \right).
\eeq
For $|z| \leq N^{-1+\kappa}$ we have that on $\Feps$,
\beq
\tilf (z) = f (Nz \cb^{-1} ) + \O (N^{-1/3+\kapp+\eps_1/10 } ),
\eeq
for $\eps_1$ sufficiently small.  
From this observation as well as Lemma \ref{lem:micro-f-est-2} we have the first estimate of the following lemma.  The second estimate easily follows from the fact that the first term in the definition of $\tilf$ dominates the others for large $z$.   
\bel \label{lem:micro-tilf} Let $\kapp >0$, and $\eps_1 >0$ sufficiently small.   Then for any $c>0$ there is a $c'>0$ such that for all $ z \in \Gamk$ we have that on $\Feps$,
\beq
\Re[ \tilf(z) ] \leq - N \cb^{-1} c' (\Re[z] - c \cb N^{-1} )  \1_{ \{ \Re[z] \leq - c \cb N^{-1} \}} + N^{-1/3+\kapp+ \eps_1/10}.
\eeq
Also, 
\beq
\Re[ \tilf(z) ] \leq - c N |E|
\eeq
if $E \leq -N^{\kapp}/N$.
\eel
The following proposition contains our first estimates for the quantities in the numerator and denominator.  Note that the proof is similar to that of Proposition \ref{prop:micro-int}, and is therefore deferred to Appendix \ref{a:micro-steepest-1}. 
\bep \label{prop:micro-steepest-1}  Let $\kapp >0$ be sufficiently small.  Then for all
sufficiently small $\eps_1 >0$, the following hold on $\Feps$; we have
\begin{align}
\int_{\gamt-\i\infty}^{\gamt+\i \infty} \exp\left[ \frac{N}{2} \left( G(z) - G ( \gamt) \right) \right] \d z &= \frac{ \cb}{N} \int_{\tilGam} \exp\left[ au+bu^{-1} - \frac{1}{2} \log(u) -a-b\right] \d u \nonumber \\
& + \O(N^{-1+4\kapp-2/3+\eps_1} ),
\end{align}
where
\beq
a = \cb \frac{\beta+ \tilm ( \gamt ) }{ 2}, \qquad b = \frac{v_1^2 \theta }{ 2 \cb}.
\eeq
Moreover,
\begin{align}
 & \int_{\Gamma^2} \frac{ \theta}{N^2 \beta} v^T \frac{1}{ (M-w)(M-z) } v \exp\left[ \frac{N}{2} (G(z) + G(w)-2 G ( \gamt) ) \right] \d z \d w \nonumber \\
= & \frac{\cb^2}{N^2} \int_{\tilGam^2} \exp\left[ a (u+w) + b(u^{-1} + w^{-1} ) - \frac{1}{2}( \log(u) + \log(w))-2a-2b \right] \frac{2b}{\cb \beta} \frac{1}{uw} \d u \d w\nonumber\\
 +& \O(N^{-2+4 \kapp +\eps_1-2/3} ).
\end{align}
where $\tilGam$ is a keyhole contour encircling $u=0$ and continuing along above and below the negative real axis.
\eep

\noindent{\bf Proofs of Theorems \ref{thm:mr-micro-over-2} and \ref{thm:mr-micro-over-3}}.  Now from \eqref{eqn:micro-Z-asymp} and Lemma \ref{lem:bessel} we have that,
\beq
N^{-1-\eps_1/10 } \leq \left|  \frac{ \cb}{N} \int_{\tilGam} \exp\left[ au+bu^{-1} - \frac{1}{2} \log(u) -a-b\right] \d u \right| \leq N^{\eps_1/10-1}
\eeq
for any sufficiently small $\eps_1 >0$ on $\Feps$.  Hence, on the event $\Feps$, we have by applying Lemma \ref{lem:bessel} and the previous proposition (recall also \eqref{eqn:micro-over-rep-f}),
\beq
\frac{1}{N} \langle \sigma^{(1)} \cdot \sigma^{(2)} \rangle = \frac{ I_{1/2} ( \lambda)^2  \lambda (\cb \beta)^{-1} }{ (b/a)^{1/2} I_{-1/2} ( \lambda)^2}  + \O(N^{-2/3+2\eps_1 + 10 \kapp} )
\eeq
by taking $\kapp>0$ and $\eps_1>0$ sufficiently small.  By direct calculation using \eqref{eqn:aug31-1} we have,
\beq
\frac{ I_{1/2} ( \lambda)^2  \lambda (\cb \beta)^{-1} }{ (b/a)^{1/2} I_{-1/2} ( \lambda)^2} =\frac{ \beta + \tilm ( \gamt)}{\beta} \left( \tanh(\sqrt{v_1^2 \theta(\beta+\tilm ( \gamt))}\right)^2.
\eeq
From this and the estimate
\beq
| \tilm( \gamt ) - \tilm ( \lambda_1 (M) ) | \leq N^{-2/3+\eps_1/10}
\eeq
which holds on $\Feps$, 
 we conclude Theorem \ref{thm:mr-micro-over-2}.  Theorem \ref{thm:mr-micro-over-3} follows from a Taylor expansion developing $\tilm( \lambda_1 (M))$ around $-1$ the the calculation,
\begin{align}
&\frac{\mathrm{d}}{\mathrm{d} m} \frac{ \beta+m}{\beta} \left( \tanh(\sqrt{v_1^2 \theta(\beta+m)}\right)^2 \bigg\vert_{m=-1} \\
= &\frac{1}{\beta} \tanh^2 ( \tillam) + \frac{\sqrt{ \tillam} }{ \beta}\frac{\tanh( \tillam)}{ \cosh^2 ( \tillam) } = \tanh ( \tillam) \frac{\sinh ( \tillam)\cosh ( \tillam)+\sqrt{ \tillam}}{ \beta \cosh^2 ( \tillam) }
\end{align}
as well as Lemma \ref{lem:cond}, where we defined $\tillam := \sqrt{ 2 v_1^2 \theta (\beta-1 ) }$. \qed

\noindent{\bf Proof of Theorem \ref{thm:mr-micro-over-4}}. 
We  turn to calculating the second moment of the overlap between two replicas.  From Proposition \ref{prop:moment-rep} we have
\begin{align}
 & \frac{1}{N^2} \langle ( \sigma^{(1)} \cdot \sigma^{(2)})^2 \rangle \nonumber \\
=& \int_{\Gamma^2} \bigg\{  \frac{1}{N^2 \beta^2} \sum_i \frac{1}{ ( \lambda_i - z )( \lambda_i -w)} - \frac{ \theta}{N^3 \beta^2} v^T[  \frac{1}{(M-z)^2(M-w)} + \frac{1}{(M-w)^2(M-z)} ]v \nonumber \\
+& \frac{\theta^2}{N^4 \beta^2} \left( v^T \frac{1}{ (M-z)M-w)} v \right)^2 \bigg\} \exp\left[ \frac{N}{2} (G(z) + G(w) -2 G(\gamt)) \right] \d w \d z \nonumber \\
\times& \left( \int_{\Gamma} \exp\left[ \frac{N}{2} (G(z)-G(\gamt)) \right] \d z \right)^{-2} := A B^{-2}
\end{align}
For the denominator we use the same analysis as above.  Fixing $\kapp >0$ and $\eps_1 >0$ small enough, we proceed by arguing along similar lines to Proposition \ref{prop:micro-steepest-1}.  We obtain the following chain of estimates that hold on $\Feps$,
\begin{align}
A &= \frac{1}{ \beta^2} \int_{\Gamk^2} \bigg\{ \frac{1}{ \cb^2(1+ \frac{Nz}{ \cb} )(1+ \frac{Nw}{ \cb} )} + 2 \frac{ \theta v_1^2}{ \cb^3} \frac{1}{ (1+ \frac{Nz}{ \cb} )^2 (1+ \frac{Nw}{ \cb} ) } + \theta^2 \frac{ v_1^4}{ \cb^4 (1+ \frac{Nz}{ \cb} )^2(1+ \frac{Nw}{ \cb} )^2}  \bigg\} \nonumber\\
&\times \exp \left[ \tilf (z) + \tilf ( w) \right] \d z \d w + \O(N^{-2+4\kapp-2/3+\eps_1} ) \nonumber \\
&= \frac{\cb^2}{N^2 \beta^2} \int_{ \tilGam^2} \bigg\{ \frac{1}{ \cb^2 uw} + \frac{2 \theta v_1^2}{ \cb^3} \frac{1}{u^2 w} + \frac{ \theta^2 v_1^4}{ \cb^4 u^2 w^2 } \bigg\} \nonumber\\
&\times \exp \left[ a(u+w) + b (u^{-1} + w^{-1} ) - \frac{1}{2} ( \log(u) + \log(w)) - 2 a - 2 b \right] \d u \d w + \O(N^{-2+4\kapp-2/3+ \eps_1} ) .
\end{align}
Hence, using Lemma \ref{lem:bessel} we see that on $\Feps$,
\begin{align}
\frac{1}{N^2} \langle ( \sigma^{(1)} \cdot \sigma^{(2)} )^2 \rangle &= \frac{1}{   \cb^2 \beta^2} \frac{a}{ b I_{-1/2}^2} \left( I_{1/2}( \lambda) ^2 + 2 \lambda I_{3/2} (\lambda) I_{1/2} ( \lambda) + \lambda^2 I_{3/2} ( \lambda) ^2 \right) + \O(N^{-2/3+10\kapp + 2\eps_1}).
\end{align}
From \eqref{eqn:bessel-recursion}, we have $I_{3/2} ( \lambda) = I_{-1/2} ( \lambda) - I_{1/2} ( \lambda) / \lambda$ and so,
\beq
 I_{1/2}( \lambda) ^2 + 2 \lambda I_{3/2} (\lambda) I_{1/2} ( \lambda) + \lambda^2 I_{3/2}^2 = \lambda^2 I_{-1/2} ( \lambda)^2.
\eeq
Therefore,
\begin{align}
\frac{1}{N^2} \langle ( \sigma^{(1)} \cdot \sigma^{(2)} )^2 \rangle &=  \frac{1}{   \cb^2 \beta^2} \frac{a}{ b I_{-1/2}^2} \lambda^2 I_{-1/2}^2 ( \lambda)+\O(N^{-2/3+10\kapp + 2\eps_1}) \nonumber\\
&= \frac{1}{ \beta^2} ( \beta + \tilm_N ( \gamt ) )^2 + \O(N^{-2/3+10\kapp + 2\eps_1}).
\end{align}
We conclude Theorem \ref{thm:mr-micro-over-4} from this and our previous estimate of $N^{-1} \langle \sigma^{(1)} \cdot \sigma^{(2)} \rangle$ from Theorem \ref{thm:mr-micro-over-2}. \qed

Finally, we turn to the proof of Theorem \ref{thm:mr-micro-over-1}.  It will suffice to show that the squared overlap $N^{-2} ( \sigma^{(1)} \cdot \sigma^{(2)} )^2$ concentrates with respect to the Gibbs measure, with high probability w.r.t. the disorder variables.  For this, we will show that its variance is $o(1)$ with high probability.  We must therefore calculate the fourth moment of the overlap.  For this we prove the following in Appendix \ref{a:prop-aug-31}. 
\bep \label{prop:aug-31}
For all sufficiently small $\eps >0$ there is an $\eps_2 >0$ so that we have with probability at least $1-N^{-\eps_2}$ that,
\begin{align}
\frac{1}{N^4} \left\langle \left( \sigma^{(1)} \cdot \sigma^{(2) } \right)^4 \right\rangle &=\left( \frac{ \beta + \tilm ( \gamt ) }{ \beta} \right)^4 + \O(N^{-2/3+\eps} ). 
\end{align}
\eep

\noindent{\bf Proof of Theorem \ref{thm:mr-micro-over-1}}.   Recall our notation $R_{12} = N^{-1} \sigma^{(1)} \cdot \sigma^{(2)}$.  From Proposition \ref{prop:aug-31} we see that there is a $c >0$ so that for all small $\eps >0$ there is an event with probability at least $1 - N^{-c\eps}$ on which the following holds,
\beq
\langle \1_{ \{ | R_{12}^2 - (1- \beta^{-1} )^2 ) | \geq t \} } \rangle \leq N^{-2/3+\eps} t^{-2}.
\eeq
Let $N^{-\delta} \geq t \geq N^{-1/3+\delta}$, and let $p$ be
\beq
p_\pm := \langle \1_{ \{ | R_{12} \mp ( 1- \beta^{ -1} ) | < t \} } \rangle.
\eeq
Then, 
\beq
1 = p_+ + p_- + \langle \1_{ \{ | R_{12}^2 - (1-\beta^{-1} )^2 | \geq t (|R_{12} | + (1-\beta^{-1} )) \} } \rangle
\eeq
For the latter term,
\beq
\langle \1_{ \{ | R_{12}^2 - (1-\beta^{-1} )^2 | \geq t (|R_{12} | + (1-\beta^{-1} )) \} } \rangle \leq \langle \1_{ \{ | R_{12}^2 - (1-\beta^{-1} )^2 | \geq t  (1-\beta^{-1} ) \} } \rangle \leq CN^{-2/3+\eps} t^{-2}
\eeq
and so $1 = p_+ + p_- + \O (N^{-2/3+\eps} t^{-2})$.  Applying Theorem \ref{thm:mr-micro-over-2} we have,
\begin{align}
(1-\beta^{-1} ) \tanh^2 ( \tillam ) &= \langle R_{12} \rangle + \O(N^{-1/3+\eps } ) \nonumber \\
&= (1-\beta^{-1} )(p_+ - p_- ) +\langle |R_{12} | \1_{ \{| |R_{12}| - (1-\beta)^{-1} |>t }\} \rangle + \O(t + N^{-1/3+\eps} ),
\end{align}
where we denoted $\tillam := \sqrt{ v_1^2 \theta ( \beta + \tilm ( \lambda_1 (M) ) ) }$.  
Using $|R_{12} | \leq 1$ we have,
\beq
\langle |R_{12} | \1_{ \{| |R_{12}| - (1-\beta)^{-1} |>t }\} \rangle \leq \langle  \1_{ \{ | R_{12}^2 - (1-\beta^{-1} )^2 | > t(1-\beta^{-1} )  \} } \rangle \leq C  N^{-2/3+\eps} t^{-2} ,
\eeq
we conclude
\beq
p_{\pm} = \frac{1 \pm \tanh ( \tillam )}{2} + \O ( N^{-1/3+\eps} + t +  N^{-2/3+\eps} t^{-2} ).
\eeq
This yields the claim. \qed

\appendix

\section{Proofs of representation formulas} \label{a:rep}

\subsection{Proof of Proposition \ref{prop:rep}}

Fix $z \in \cc$ with $\Re [z] > \lambda_1 (M)$ and consider the Gaussian integral over $\rr^N$,
\beq
\int_{\rr^N} \exp \left[ \frac{\beta}{2} x^T (M - z) x + \lambda v^T x \right] \d x = \left( \frac{ 2 \pi}{\beta} \right)^{ \frac{N}{2} } \exp \left[ - \frac{\lambda^2}{2 \beta} v^T (M-z)^{-1} v -\frac{1}{2} \sum_i \log( z - \lambda_i) \right].
\eeq
On the other hand,
\beq
\int_{\rr^N} \exp \left[ \frac{\beta}{2} x^T (M - z) x + \lambda v^T x \right] \d x = \int_{0}^\infty \exp \left[  - \frac{ \beta z}{2} r^2 \right] r^{N-1} \J (r) \d r ,
\eeq
where 
\beq
\J (r) := \int_{ \| \sigma \|_2 =1 } \exp \left[ \frac{\beta}{2} r^2 \sigma^T M \sigma + r \lambda v^T \sigma \right] \d \sigma ,
\eeq
and  $\d \sigma$ is uniform measure over the $N-1$ sphere of radius $1$.  
Note that the LHS of \eqref{eqn:representation1} equals $N^{\frac{N-1}{2}} \J ( \sqrt{N})$.   Making a change of variable $u = \frac{\beta}{2} r^2$, we see that
\beq
 \int_{0}^\infty \exp \left[  - \frac{ \beta z}{2} r^2 \right] r^{N-1} \J (r) \d r = \int_{0}^\infty \exp\left[ - z u \right]  \frac{1}{\beta} \left( \frac{2 u}{ \beta} \right)^{\frac{N-2}{2}} \J ( (2 u / \beta)^{1/2} ) \d u.
\eeq
The RHS is the Laplace transform of a function that equals the LHS of \eqref{eqn:representation1} at $u= \frac{ \beta}{2} N$ times the constant $N^{-1/2} \beta^{-1}$.  The result \eqref{eqn:representation1} follows from the Laplace inversion formula.  

For \eqref{eqn:representation2}, we begin by considering for fixed $z, w \in \cc$ and $\Re[z], \Re[w] > \lambda_1 (M)$, the  Gaussian integral,
\begin{align}
&\int_{\rr^N \times \rr^N} \exp \left[ \frac{\beta}{2} x^T (M- z ) x + \frac{\beta}{2} y^T(M- w ) y + t x^T y + \lambda v^T( x + y ) \right] \d x \d y \nonumber \\
= &( 2 \pi)^N \exp \left[ - \frac{ \lambda^2}{2} \sum_{i=1}^N ( v^T u_i )^2 \frac{ \beta ( 2 \lambda_i - w - z ) - 2 t}{\beta^2 (  \lambda_i-z )(\lambda_i - w ) - t^2 }  - \frac{1}{2} \sum_{i=1}^N \log ( \beta^2 (z- \lambda_i )( w - \lambda_i ) - t^2 ) )\right].
\end{align}
On the other hand, considering this integral in polar coordinates and making the same change of variables as above we have,
\begin{align}
& \int_{\rr^N \times \rr^N} \exp \left[ \frac{\beta}{2} x^T (M- z ) x + \frac{\beta}{2} y^T(M- w ) y + t x^T y + \lambda v^T( x + y ) \right] \d x \d y \nonumber \\
= & \int_{0}^\infty \int_0^\infty \exp \left[ - \frac{ \beta z}{2} r_1^2 - \frac{ \beta w }{2} r_2^2 \right] (r_1 r_2)^{N-1} \J (r_1, r_2 ) \d r_1 \d r_2 \nonumber \\
= & \int_{0}^\infty \int_0^\infty \exp [ -zu_1 - w u_2 ] \frac{1}{ \beta^2} \left(  \frac{ 4 u_1 u_2}{ \beta^2} \right)^{ \frac{N-2}{2} } \J ( (2 u_1 / \beta)^{1/2}, (2 u_2 / \beta)^{1/2} ) \d u_1 \d u_2 \label{eqn:rep-proof1}
\end{align}
where
\beq
\J (r_1, r_2 ) := \int_{ || \sigma_1||_2 , ||\sigma_2||_2 =1} \exp \left[ \frac{\beta}{2} r_1^2 \sigma_1^T M \sigma_1 + \frac{ \beta}{2} r_2^2 \sigma_2^T M \sigma_2 + t r_1 r_2 \sigma_1^T \sigma_2 + \lambda v^T ( r_1 \sigma_1 + r_2 \sigma_2 ) \right] \d \sigma_1 \d \sigma_2 .
\eeq
We recognize the last line of \eqref{eqn:rep-proof1} as the multidimensional Laplace transform of a function that equals the LHS of \eqref{eqn:representation2} at $u_1 = u_2 = \frac{\beta}{2} N$ times the constant $N^{-1} \beta^{-2}$.  The result follows. \qed

\subsection{Proof of Proposition \ref{prop:moment-rep}}

We first consider the Gaussian integral
\begin{align}
&\int_{\rr^{2N}} x \cdot y \exp \left[ \frac{\beta}{2} \left( x^T(M-z)x + y^T (M-w) y \right) + \lambda v^T (x+y) \right] \nonumber \\
= & \frac{ \lambda^2}{ \beta^2} v^T \frac{1}{ (M-z)(M-w)} v \left( \frac{2 \pi}{ \beta} \right)^N \exp \bigg\{ - \frac{ \lambda^2}{ 2 \beta} \left( v^T (M-z)^{-1} v + v^T (M-w)^{-1} v \right)\nonumber\\
-& \frac{1}{2} \sum_{i} \log(z- \lambda_i ) + \log(w - \lambda_i )\bigg\}.
\end{align}
On the other hand,
\begin{align}
&\int_{\rr^{2N}} x \cdot y \exp \left[ \frac{\beta}{2} \left( x^T(M-z)x + y^T (M-w) y \right) + \lambda v^T (x+y) \right] \nonumber \\
= & \int_0^\infty \int_0^\infty \exp\left[ \frac{-\beta z}{2} r^2 - \frac{\beta w}{2} s^2 \right] (rs)^{N-1} \J_1 (r, s) \d r \d s
\end{align}
where,
\beq
\J_1(r, s) = \int_{ || \sigma||_2 = ||\omega||_2 =1 } (rs \omega \cdot \sigma ) \exp\left[ \frac{\beta}{2} \left( r^2 \sigma M \sigma + s^2 \omega M \omega \right) + \lambda v^T (r\sigma + s \omega ) \right] \d \sigma \d \omega.
\eeq
The quantity on the line \eqref{eqn:2repa} that we want to calculate equals $N^{N-1} \J_1 ( \sqrt{N}, \sqrt{N} )$.   Proceeding as in the proof of Proposition \ref{prop:rep} yields the claimed representation. 

For the second representation we begin with the Gaussian integral,
\begin{align}
&\int_{\rr^{2N}} (x \cdot y)^2 \exp \left[ \frac{\beta}{2} \left( x^T(M-z)x + y^T (M-w) y \right) + \lambda v^T (x+y) \right] \nonumber \\
=& \bigg\{ \frac{1}{ \beta^2} \sum_i \frac{1}{ ( \lambda_i - z ) ( \lambda_i  - w ) }  - \frac{ \lambda^2}{ \beta^3} v^T \left( \frac{1}{(M-z)^2(M-w) } + \frac{1}{ (M-w)^2 (M-z)} \right) v \nonumber \\
+& \frac{ \lambda^4}{ \beta^4} \left( v^T \frac{1}{ (M-w)(M-z) } v \right)^2 \bigg\}  \nonumber \\
\times &  \left( \frac{ 2 \pi}{ \beta} \right)^N \exp \left[ - \frac{ \lambda^2}{ 2 \beta} \left( v^T (M-z)^{-1} v + v^T (M-w)^{-1} v \right)- \frac{1}{2} \sum_{i} \log(z- \lambda_i ) + \log(w - \lambda_i )  \right]
\end{align}
We have also that,
\begin{align}
&\int_{\rr^{2N}} (x \cdot y)^2 \exp \left[ \frac{\beta}{2} \left( x^T(M-z)x + y^T (M-w) y \right) + \lambda v^T (x+y) \right] \nonumber \\
=&\int_0^\infty \int_0^\infty \exp\left[ - \frac{\beta z}{2} r^2 - \frac{\beta w}{2} s^2 \right] (rs)^{N-1} \J_2 (r, s) \d r \d s
\end{align}
where,
\beq
\J_2 (r, s) = \int_{ || \sigma||_2 = ||\omega||_2 =1 } (rs \omega \cdot \sigma )^2 \exp\left[ \frac{\beta}{2} \left( r^2 \sigma M \sigma + s^2 \omega M \omega \right) + \lambda v^T (r\sigma + s \omega ) \right] \d \sigma \d \omega.
\eeq
Note that $N^{N-1} \J_2 ( \sqrt{N}, \sqrt{N} )$ is the quantity we want to calculate.  We proceed as before. \qed

\section{Proofs of auxiliary results in Gaussian regime} \label{a:gaussian}

\subsection{Proof of Lemma \ref{thm:mr-gauss-over-b1}} \label{a:gauss-1}

  We start with the first part of the lemma, computing $A_N$:
\begin{align}
A_N \beta & = - \frac{ \theta ( \msc' ( \hatgam) )^2}{ g'' ( \hatgam) } - \frac{1}{2} \msc ( \hatgam) \nonumber \\
&= \frac{\msc' ( \hatgam)}{ 2 g'' ( \hatgam)}\left( - \msc ( \hatgam) + \theta \msc ( \hatgam) \msc'' ( \hatgam) / \msc' ( \hatgam) - 2 \theta \msc' ( \hatgam) \right) \nonumber \\
&= \frac{\msc' ( \hatgam)}{ 2 g'' ( \hatgam)}\left(  - \msc ( \hatgam) +2 \theta \msc' ( \hatgam) / (1 - \msc^2 ( \hatgam ))- 2 \theta \msc' ( \hatgam) \right) \nonumber \\
&=  \frac{\msc' ( \hatgam)}{ 2 g'' ( \hatgam)}\left( - \msc ( \hatgam) + 2 \theta \msc' (\hatgam)  \frac{\msc^2 ( \hatgam)}{ 1 - \msc^2 ( \hatgam)} \right) \nonumber \\
&= \frac{\msc' ( \hatgam)}{ 2 g'' ( \hatgam)}\left( - \msc ( \hatgam) + 2 \theta ( \msc' ( \hatgam))^2 \right).
\end{align}
In the third line we used the identity $ \msc(z) \msc'' (z) = 2 ( \msc' (z))^2 / (1 - \msc^2 (z) )$, and in the last line the identity $\msc' (z) = \msc^2 (z) / (1 - \msc^2 (z))$.  For the last line, the prefactor is of order $(1 + \theta / \kappa)^{-1}$ and the term in the brackets is of order $(1 + \theta / \kappa)$. 

For the second part of the lemma, let $X := N^{-1/2} v \cdot \sigma-B_N$,   with  $B_N$ as in Theorem \ref{thm:gauss-ext}. We have, on the event of that theorem,
\beq
\frac{1}{t} \left( \langle \e^{t X} \rangle -1 \right) = \frac{1}{t} \left( \e^{t^2 A_N  + \O(N^{-c_1} ) } -1 \right) =  \O ( t^{-1} N^{-c_1} + t),
\eeq
for $|t| \leq 1$. 
On the other hand,
\beq
\frac{1}{t} \left( \langle \e^{t X} \rangle -1 \right) = \langle X \rangle + t \O( \langle X^2 \e^{ t|X| }  \rangle).
\eeq
We choose $t = N^{-c_1/2}$. With this choice,
\beq
 \langle X^2 \e^{ t|X| }  \rangle \leq C \langle \e^{X} + \e^{-X} \rangle \leq C'.
\eeq
We conclude the proof.   \qed

\subsection{Proof of Lemma \ref{lem:gauss-1-f}} \label{a:gauss-1-f}

First we note that by taking $\delta >0$ sufficiently small we can assume that the LHS of \eqref{eqn:zw-ass} is smaller than $N^{-\eps} \kappa$ for some small $\eps >0$.  This, combined with our assumption that $| \gamh - \gamma| \leq \kappa / \log(N)$ implies that with overwhelming probability that $\Re[z]$ and $\Re[w]$ will both be larger than $2+N^{-2/3+\tau/2}$.  This allows for the use of the rigidity and local estimates of Theorems \ref{thm:rigi} and \ref{thm:ll-derivs}, as well as the Taylor expansion of quantities appearing in denominators below.

 We start by rewriting $G_1(z, w)$ as follows.
\begin{align}
G_1 (z, w) &= \beta(z + w) - \frac{1}{N} \sum_{i=1}^N \log  ( ( z - \lambda_i )( w - \lambda_i ) ) \\
&- \frac{1}{N} \sum_{i=1}^N \log ( 1 - t^2 / (  ( z - \lambda_i )( w - \lambda_i ) )) \label{eqn:G1aa} \\
&- \frac{ \theta}{N} \sum_{i=1}^N v_i^2 \frac{ ( 2 \lambda_i - z - w ) }{ ( \lambda_i - w)( \lambda_i - z ) - t^2 } \label{eqn:G1bb} \\
&+ 2 t \frac{ \theta \beta}{N} \sum_{i=1}^N v_i^2 \frac{1}{ \beta^2 ( \lambda_i  - w)( \lambda_i - z) - t^2 }. \label{eqn:G1cc}
\end{align}
First, via Taylor expansion and the rigidity estimates of Theorem \ref{thm:rigi} we have,
\begin{align}
- \frac{1}{N} \sum_{i=1}^N \log ( 1 - t^2 / (( z - \lambda_i )( w - \lambda_i ) )) &= t^2 \frac{1}{N } \sum_i \frac{1}{( z - \lambda_i )(w- \lambda_i )} + \O ( t^4 \kappa^{-5/2} ) \nonumber \\
&= t^2  m' ( \gamh ) + \O ( t^4 \kappa^{-5/2} + t^2\kappa^{-5/4} N^{\delta-1/2} ).
\end{align}
We next turn to \eqref{eqn:G1bb}.  We have
 the expansion,
\begin{align}
\frac{ \lambda_i - z}{ ( \lambda_i - w)( \lambda_i - z) - t^2 } &= \frac{ \lambda_i - \gamh}{ ( \lambda_i - \gamh )^2 - t^2 } + ( \gamh - w) \left( - \frac{ ( \lambda_i - \gamh )^2}{ ( ( \lambda_i - \gamh )^2 - t^2 )^2 } \right)\nonumber \\
&+ ( \gamh - z) \left( \frac{1}{ ( \lambda_i - \gamh )^2 -t^2 }  - \frac{ ( \lambda_i - \gamh )^2}{ ( ( \lambda_i - \gamh )^2 - t^2 )^2 } \right) \nonumber \\
&+ ( \gamh - z)^2 \left( \frac{ ( \lambda_i - \gamh)^3}{ ( ( \lambda_i - \gamh)^2 -t^2 )^3 } - \frac{ ( \lambda_i - \gamh )}{ ( ( \lambda_i - \gamh)^2 - t^2 )^2 } \right) \nonumber \\
&+ ( \gamh - w)^2 \left( \frac{ ( \lambda_i - \gamh )^3 }{ ( ( \lambda_i - \gamh )^2 -t^2 )^3} \right) \nonumber \\
&+2 ( \gamh-z)( \gamh-w) \left( \frac{ ( \lambda_i - \gamh )^3}{ ( ( \lambda_i - \gamh)^2 - t^2 )^3} - \frac{ \lambda_i - \gamh}{ ( ( \lambda_i - \gamh )^2 - t^2 )^2} \right)\nonumber \\
& + \O ( (|z - \gamh|^3 + |w  - \gamh|^3)| \lambda_i - \gamh|^{-4} ),
\end{align}
as well as
\begin{align}
\frac{1}{ \lambda_i - z } &= \frac{1}{  \lambda_i - \gamh} + ( z - \gamh)\frac{1}{ ( \lambda_i - \gamh)^2} + (z- \gamh)^2 \frac{1}{ ( \lambda_i - \gamh)^3 } + \O (|z- \gamh |^3 | \lambda_i - \gamh|^{-4} ).
\end{align}
Using the above two expansions as well as Theorems \ref{thm:rigi} and \ref{thm:iso} (the latter to estimate $|v_i|^2 \leq N^{\eps}$ with overwhelming probability) we find that with overwhelming probability,
\begin{align}
&- \frac{ \theta}{N} \sum_{i=1}^N v_i^2 \frac{ ( 2 \lambda_i - z - w ) }{  ( \lambda_i - w)( \lambda_i - z ) - t^2 } + \frac{ \theta}{N} \sum_{i=1}^N v_i^2 \left( \frac{1}{ \lambda_i - z } + \frac{1}{ \lambda_i - w } \right) \nonumber \\
=&2 \frac{ \theta}{N} \sum_{i=1}^N v_i^2 \left( \frac{1}{  \lambda_i - \gamh }  - \frac{ \lambda_i - \gamh }{ ( \lambda_i - \gamh)^2 - t^2 } \right) \nonumber \\
+& [ ( \gamh - z) + ( \gamh - w) ] \left( \frac{ \theta}{N} \sum_{i=1}^N v_i^2 \left( \frac{ 2 ( \lambda_i - \gamh )^2}{ ( ( \lambda_i - \gamh)^2 - t^2 )^2} - \frac{1}{ ( \lambda_i - \gamh )^2 - t^2 } - \frac{1}{( \lambda_i - \gamh )^2 } \right) \right) \nonumber \\
+& N^{\eps}\O \left( \theta \kappa^{-2}N^{-1}N^{2 \delta}(1+ \theta / \kappa)^{-1} (t^2 \kappa^{-1} + N^{\delta} \kappa^{1/4} N^{-1/2} )  \right) .
\end{align}
Note that the terms quadratic in $(z-\gamh )$ and $(w- \gamh)$ were estimated using the cancellation between the above two expansions which gives the error terms that have the $t^2$ terms above.  

We next turn to \eqref{eqn:G1cc}. 
We use,
\begin{align}
\frac{1}{ ( \lambda_i - w )( \lambda_i - z) - t^2 } &= \frac{1}{ ( \lambda_i - \gamh)^2 - t^2}- ( ( \gamh - z) + ( \gamh - w) ) \frac{ ( \lambda_i - \gamh)}{ ( ( \lambda_i - \gamh )^2 - t^2 )^2} \nonumber \\
+& N^{\eps}  \O \left( ( | z- \gamh|^2 + |w - \gamh |^2 ) | \lambda_i - \gamh |^{-4} \right)
\end{align}
to find, with overwhelming probability,
\begin{align}
2 t \frac{ \theta }{N} \sum_{i=1}^N v_i^2 \frac{1}{ ( \lambda_i  - w)( \lambda_i - z) - t^2 } &= \frac{ 2 t \theta}{N} \sum_{i=1}^N v_i^2 \frac{1}{ ( \lambda_i - \gamh)^2 - t^2 } \nonumber \\
&+ [ ( \gamh - z )+ ( \gamh - w ) ] \left( \frac{ - 2 \theta t}{N}  \sum_{i=1}^N v_i^2 \frac{ \lambda_i - \gamh }{ ( ( \lambda_i - \gamh )^2 - t^2 )^2 } \right) \nonumber \\
&+ \O ( N^{2 \delta} |t| N^{-1} \theta ( 1+ \theta/ \kappa)^{-1} \kappa^{-2} ).
\end{align}
Note,
\beq
2\frac{  ( \lambda_i - \gamh)^2-  t( \lambda_i - \gamh)}{ ( ( \lambda_i - \gamh )^2 - t^2 )^2} - \frac{1}{ ( \lambda_i - \gamh)^2 - t^2 } = \frac{1}{ ( \lambda_i - \gamh + t)^2}.
\eeq
We have so far arrived at the following expansion which holds with overwhelming probability,
\begin{align}
G_1 (z, w) &= G(z) + G(w)\nonumber  \\
&+ [ ( \gamh - z) + ( \gamh - w)] \left( - \theta N^{-1} v^T(M-\gamh)^{-2} v   + \theta N^{-1} v^T (M- \gamh + t)^{-2} v \right)\nonumber \\
&+ t^2 m_N' ( \gamh) + 2 \theta N^{-1} v^T (M- \gamh)^{-1} v - 2 \theta N^{-1} v^T (M - \gamh + t )^{-1} v  \nonumber \\
&+N^{\eps}\O ( t^4 \kappa^{-5/2} + t^2\kappa^{-5/4} N^{\delta-1/2} ) \nonumber\\
&+ N^{\eps} \O ( \theta \kappa^{-2}N^{-1}N^{2 \delta}(1+ \theta / \kappa)^{-1} (|t| + t^2 \kappa^{-1} + N^{\delta} \kappa^{1/4} N^{-1/2} ).
\end{align}
Now, using
\begin{align}
&\theta N^{-1} \bigg\{ -v^T (M - z)^{-1} v + ( \gamh - z)( v^T (M- \gamh + t)^{-2} v - v^T (M- \gamh)^{-2} v ) \nonumber \\
+& v^T (M- \gamh)^{-1} - v^T(M- \gamh+t)^{-1} v ) \bigg\}\nonumber \\
= & - \theta N^{-1} v^T (M +t - z)^{-1} v+  \O ( N^{2 \delta} |t| N^{-1} \theta ( 1+ \theta/ \kappa)^{-1} \kappa^{-2} ).
\end{align}
We conclude the claim. \qed

\section{Proofs of auxiliary results in intermediate regime} \label{a:intermediate}

\subsection{Proof of Lemma \ref{lem:int-gams-comp}} \label{a:int-gams-comp}

We have,
\beq
\beta - 1 \geq \frac{ \theta (v^T u_1 (M) )^2 }{\hats_M^2},
\eeq
and so on the event $\Feps$, $\hats_M \geq c N^{-\eps_1/100}$.  The same estimates also clearly hold for $s$ and $\hats_H$ (for $\hats_H$ we also use that the estimates  \eqref{eqn:compare1} hold on $\Feps$).  For an upper bound we use the fact that the delocalization estimates hold for $v^T u_i (M)$ to see that,
\beq
\frac{\theta}{N^{4/3}} v^T \frac{1}{(M- \hatgam_M)^2} v \leq \theta N^{\eps} N^{ \eps_1/1000} (\hats_M)^{-2} + \sum_{i > N^{\eps}}  \frac{N^{\eps_1/10^6}}{ N^{4/3} ( \lambda_i (M) - \lambda_1 (M) )^2 }
\eeq
for any $\eps >0$.  Taking $\eps = \eps_1 /1000$, we see that by the rigidity estimates, 
 the second term is $\O (N^{-\eps/3+\eps_1/10^6} )$.   We therefore get the inequality,
\beq
\frac{\beta-1}{2} \leq \theta N^{\eps_1 / 500 } ( \hats_M )^{-2}
\eeq
and so  we conclude the estimate for $\hats_M$.  A similar conclusion holds for $s$ and $\hats_H$.
On $\Feps$ we may now assume that $N^{-\eps_1/100 } \leq s \leq N^{\eps_1/100}$ and the same for $\hats_M$ and $\hats_H$.  In the region $N^{-2/3 - \eps_1/50} \leq E \leq N^{-2/3+\eps_1/50}$, we see that
\beq
N^{- \eps_1/10 } \leq \frac{1}{N^{2/3}} \frac{\d }{ \d E} \frac{ \theta}{N^{4/3}} v^T \frac{1}{( M - ( \lambda_1(M) - E) )^2 } v \leq N^{ \eps_1/10}.
\eeq
From the lower bound and the mean value theorem we see that 
\beq \label{eqn:t-1}
N^{- \eps_1/10} |s - \hats_M | \leq \left| \frac{1}{N^{4/3}} v^T \frac{1}{(M - \gamma)^2 } v - \frac{1}{N^{4/3}} v^T \frac{1}{ (M - \hatgam_M)^2 } v \right| 
\eeq
From the fact that the rigidity and level repulsion estimates are assumed to hold on $\Feps$, we
have that the estimate
\beq
\left| m_N ( \gamma) + 1 \right| \leq N^{ \eps_1/50 - 1/3 }
\eeq
holds on $\Feps$.  This, together with \eqref{eqn:t-1} and the definitions of $\gamma, \hatgam_M$ imply,
\beq
|s - \hats_M | \leq N^{ \eps_1/5 - 1/3}.
\eeq

We now turn to proving the estimate on the difference $\hats_H - \hats_M$.  
Consider now $E = \lambda_1 (M) + u N^{-2/3}$ for $N^{-2 \eps_2 } \leq u \leq N^{2\eps_2}$.  Let $\eps_2 = 10^{-6}$.  Let $\eta = N^{-2/3 - 1/1000}$.  
We begin by rewriting,
\begin{align}
& N^{-4/3} v^T (H- E)^{-2}  v  - N^{-4/3} v^T (M - E)^{-2} v  \label{eqn:aa} \\
=& N^{-4/3} v^T(H - E)^{-2} v - N^{-4/3} v^T (H - (E+ \i \eta ) )^{-2} v \label{eqn:a1} \\
 +& N^{-4/3} v^T (M - (E + \i \eta ) )^{-2} v - N^{-4/3} v^T (M - E)^{-2} v \label{eqn:a2} \\
+ & N^{-4/3} v^T (H- (E+ \i \eta) )^{-2}  v  - N^{-4/3} v^T (M - (E + \i \eta ) )^{-2} v  \label{eqn:a3}
\end{align}
From the delocalization estimates that hold on $\Feps$ we have that $|v^T u_i(H)|+ |v^T u_i (M) | \leq N^{2\eps_1/10^6}$ for all $i$, and we can assume $\eps_1 < \eps_2/100$.  Using the rigidity estimates that hold on $\Feps$ we see that the term  \eqref{eqn:a1} can estimated as,
\begin{align}
\left| N^{-4/3} v^T(H - E)^{-2} v - N^{-4/3} v^T (H - (E+ \i \eta ) )^{-2} v \right| \leq & C \eta N^{\eps_1/10} \sum_{i=1}^N \frac{1}{ N^{4/3} | \lambda_i (H) - E |^3 } \notag \\
\leq & C \eta N^{ 10 \eps_2} N^{2/3}, \label{eqn:a4}
\end{align}
for all $E$ specified above.  
The term \eqref{eqn:a2} is similar.  For \eqref{eqn:a3} we have the estimates \eqref{eqn:isocompare1} assumed to hold on $\Feps$ which gives an estimate of (using the Cauchy integral formula),
\begin{align}
& \left| N^{-4/3} v^T (H- (E+ \i \eta) )^{-2}  v  - N^{-4/3} v^T (M - (E + \i \eta ) )^{-2} v  \right| \notag \\
\leq & \frac{N^{10\eps_2}}{\eta N^{1/3}} \left( \frac{1}{N^{1/2} } + \frac{1}{ (N \eta )^2} + \frac{1}{N \sqrt{ \eta }} + \frac{N^{-1/3}}{N \eta} \right) \label{eqn:a5}
\end{align}
Note we used $\Im[\msc] \leq C ( |E-2|^{1/2} + \eta^{1/2})$.  
Under our assumptions on $\eps_2$ and $\eta$ we see that both \eqref{eqn:a4} and \eqref{eqn:a5} are $\O (N^{- \alpha } )$ for some small $\alpha >0$, on the event $\Feps$; note that $\alpha >0$ is independent of $\eps_1$.  We can then write,
\beq \label{eqn:inter1}
\left| N^{-4/3}v^T (M - \hatgam_H)^{-2} v - N^{-4/3} v^T (H - \hatgam_H)^{-2} v \right| = \left| N^{-4/3} v^T (M - \hatgam_H )^{-2} v - N^{-4/3} v^T (M - \hatgam_M )^{-2} v \right|.
\eeq
On $\Feps$ we can assume that $| \lambda_1 (M) - \lambda_1 (H)| \leq N^{-5/6}$ and that $N^{-\eps} \leq \hats_M, \hats_H \leq N^{\eps}$, with $\eps < \eps_2 /2$, with $\eps_2$ as above, and $\eps < \alpha / 100$, where $\alpha$ is as above.  As we have seen, the LHS of \eqref{eqn:inter1} is  $\O (N^{-\alpha})$ on $\Feps$.  
From the inequality which holds for $a, b > \lambda_1 (M)$,
\begin{align}
\left| v^T (M-a)^{-2} v - v^T (M-b)^{-2} v\right| =& |b-a|  \left| v^T [(M-a)^{-2}(M-b)^{-1} + (M-a)^{-1} (M-b)^{-2} ]v \right| \nonumber \\
\geq & |b-a| (v^T u_1 (M))^2| (\lambda_1 - a)^{-2} (\lambda_1 - b)^{-1}|,
\end{align}
we find
\beq
\left| N^{-4/3} v^T (M - \hatgam_H )^{-2} v - N^{-4/3} v^T (M - \hatgam_M )^{-2} v \right| \geq N^{- \eps_1} | \hats_M - \hats_H |.
\eeq
  On the other hand, we assumed that the LHS of \eqref{eqn:inter1} is less than $N^{-\alpha}$, and $\eps_1 < \alpha /100$.  Hence, we have that there is some $c_1 >0$ so that $|\hats_M - \hats_H | \leq N^{-c_1}$ on the event $\Feps$. \qed

\subsection{Proof of Lemma \ref{lem:conv-2}} \label{a:conv-2}

 First note that $s_n \leq s$.  We have the estimates,
\begin{align}
\sum^N_{j=n+1} \frac{ g_j^2}{ N^{4/3} ( \mu_j - \mu_1 )^2 } &\geq \sum^N_{j = n+1 } \frac{ g_j^2}{ N^{4/3} ( \mu_j - \mu_1 - sN^{-2/3} )^2 } \nonumber \\
&=  \sum_{j=1}^n \frac{g_j^2}{ N^{4/3} ( \mu_j - \mu_1 - s_n N^{-2/3} )^2 } - \sum_{j=1}^n \frac{ g_j^2}{ N^{4/3}( \mu_j - \mu_1  - sN^{-2/3} )^2 } \nonumber\\
&= ( s - s_n ) \sum_{j=1}^n g_j^2 \frac{ N^{2/3} ( [ \mu_1 + s N^{-2/3} - \mu_j ] + [ \mu_1 + s_n N^{2/3} - \mu_j ] )}{ N^{8/3} ( \mu_j - \mu_1 - sN^{-2/3} )^2 ( \mu_j - \mu_1 - s_n N^{-2/3} )^2 } \nonumber\\
& \geq (s- s_n ) \frac{g_1^2}{s^3}.
\end{align}
Due to the previous lemma, for any $\eps >0$ there is a $\delta >0$ on which $g_1^2 /s^3 \geq \delta$ with probability $1 - \eps$.  Hence, with probability at least $1 - \eps$, we have that for all $n$,
\beq
(s -s_n) \leq \frac{1}{ \delta} \sum_{j=n+1} \frac{g_j^2}{N^{4/3} ( \mu_j - \mu_1 )^2 }.
\eeq
Let $K$ be so that the event of Lemma \ref{lem:RoughParticleEst} holds with probability at least $1- \eps$.  As long as $n \geq K$, we have that there is an event with probability at least $1- 2\eps$ on which,
\beq
(s -s_n) \leq \frac{C}{ \delta} \sum_{j=K}^N \frac{g_j^2}{j^{4/3}}.
\eeq
Note that this inequality holds for every $n > K$ and $C, \delta$ are independent of $K$.
Now if $K$ is sufficiently large then by Markov's inequality with probability at least $1- \eps$ the term on the RHS is less than $\eps$.  This yields the claim. \qed

\subsection{Proof of Lemma \ref{lem:conv-3}} \label{a:conv-3}

\proof The second estimate is an easy consequence of Lemmas \ref{lem:conv-1} and \ref{lem:conv-2}.  For the first estimate, we first let $n_0$, $C_0 >0$, $\delta, c_0$ and $\F$ be the event with $\pp [ \F] \geq 1- \eps$ on which,
\beq
c_0 j^{2/3} \leq N^{2/3} | \mu_j - \mu_1 |  \leq c_0^{-1} j^{2/3}, \qquad j \geq n_0
\eeq
and also
\beq
| \mu_1 - 2 | N^{2/3} \leq C_0
\eeq
and $\delta \leq s \leq \delta^{-1}$.  Note that $c_0$ does not depend on $\eps$ but the other parameters do.  For $n\geq n_0$ (possibly increasing $n_0$), we estimate 
\begin{align}
\frac{1}{N^{2/3}} & \left| \sum_{j =n+1}^N \left(  \frac{ g_j^2}{ ( \mu_j - \mu_1 - sN^{-2/3} ) } - \frac{1}{ \gamma_j - 2 } \right) \right| \1_\F \nonumber\\
\leq &   \sum^N_{j=n+1}C  |g_j|^2 \frac{ N^{2/3}| \mu_j - \gamma_j | \1_{ \{ (\mu_j - 2) N^{2/3} \leq - c_0/2 j^{2/3}  \} }+ \delta^{-1} + C_0  } {j^{4/3}}\nonumber \\
+ & \frac{1}{N^{2/3}} \left| \sum^N_{j=n+1} \frac{ g_j^2 -1 }{( \gamma_j- \gamma_1 ) } \right|. \label{eqn:gj}
\end{align}
The second moment of the last term is bounded by $C n^{-1/3}$.   From Section 6 of \cite{SSK-ls} we have the estimate (possibly increasing $n_0$ if necessary), for all $j > n_0$,
\beq
\ee[ N^{2/3}| \mu_j - \gamma_j | \1_{ \{ (\mu_j - 2) N^{2/3} \leq - c_0/2 j^{2/3}  \} }] \leq \frac{C \log(j)^2}{ j^{1/3}}.
\eeq
Hence, the first term on the RHS of \eqref{eqn:gj} has expectation bounded by $C n^{-1/5}$.  So if $n > n_1$, with $n_1$ sufficiently large, than these terms are less than $\eps$ with probability at least $1- \eps$. 

Now, fix $\eps_1 >0$.  By taking $n_1$ sufficiently large, we may assume that for all $n > n_1$ that $|s-s_n| < \eps_1$ with probability at least $1- \eps$.  Assume $\eps_1 < \delta/2$.  Call this event $\F_{1}$.  Let $C_2 $ be a constant so that $\sum_{j=1}^{n_0} g_j^2 \leq C_2$ with probability at least $1- \eps$ and call this event $\F_2$.  Note that the constant $C_2$ depends only on $n_0$ and not on $\eps_1$ or $n_1$.  On $\F \cap \F_1 \cap \F_2$,
\begin{align}
\frac{1}{N^{2/3}}& \left| \sum_{j=2}^n g_j^2 \left( \frac{1}{ \mu_j - \mu_1 - s N^{-2/3} } - \frac{1}{ \mu_j - \mu_1 - s_n N^{-2/3} } \right) \right| \nonumber\\
\leq & |s - s_n | 10 \delta^{-2} C_2 + \eps_1 \sum_{j=n_0+1}^n \frac{g_j^2}{ c_0^2 j^{4/3} }\nonumber\\
 \leq & \eps_1 \left( 10 \delta^{-2} C_2 + \sum_{j=n_0+1}^n \frac{g_j^2}{ c_0^2 j^{4/3} } \right)
\end{align}
By Markov's inequality, there is an $M>0$ so that the second term in the brackets on the last line is less than $M$ with probability at least $1-\eps$ (this $M$ does not depend on $\eps_1$ or $n$).  Hence, by taking $\eps_1$ small enough depending on $C_2$, $\delta$ and $M$, we get that this is less than $\eps$. \qed

\subsection{Proof of Lemma \ref{lem:G-hess}} \label{a:G-hess}

The diagonal elements are equal to,
\begin{align}
\partial_z^2 G_1 ( \gam1, \gam1 ) &= \frac{1}{N} \sum_{i=1}^N \frac{ ( \lambda_i - \gam1 )^2}{ ( ( \lambda_i - \gam1)^2 - N^{-4/3} t^2 )^2 } \nonumber\\
-& 2 \frac{ \theta}{N^{4/3}} \sum_{i=1}^N v_i^2 \frac{ (\lambda_i - \gam1 )}{ ( ( \lambda_i - \gam1 )^2 - N^{-4/3} t^2 )(\lambda_i - \gam1 +N^{-2/3} t )^2}
\end{align}
and off-diagonal elements equal to,
\begin{align}
\partial_z \partial_w G_1 ( \gam1, \gam1 ) &= \frac{1}{N} \sum_{i=1}^N \frac{ N^{-4/3} t^2}{ ( ( \lambda_i - \gam1 )^2 -N^{-4/3} t^2 )^2} \nonumber\\
+& 2 \frac{ \theta}{N^{4/3}} \sum_{i=1}^N v_i^2 \frac{N^{-2/3} t}{ ( ( \lambda_i - \gam1 )^2 - N^{-4/3} t^2 ) ( \lambda_i - \gam1 + N^{-2/3} t)^2}
\end{align}
By \eqref{eqn:int-over-a2} and the fact that $\gamma - \lambda_1 \geq N^{-2/3-\eps_1/100}$ on $\Feps$ (see Lemma \ref{lem:int-gams-comp}) we see that all the terms contributing to $\del_z^2 G( \gam1, \gam1 )$ and $\del_z \del_w G_1 ( \gam1, \gam1)$ are positive. We then see that for sufficiently small $\eps_1>0$, that
\beq
N^{-\eps_1/10} \leq \tr  N^{-2/3} \nabla^2 G_1 ( \gam1 , \gam1 ) \leq N^{\eps_1/10}.
\eeq
We calculate,
\begin{align}
\partial_z^2 G_1 ( \gam1, \gam1 ) - \partial_z \partial_w G_1 ( \gam1, \gam1 ) 
=  \frac{1}{N} \sum_{i=1}^N  \frac{1}{ ( \lambda_i - \gam1 )^2 - N^{-4/3} t^2 } - 2 \frac{ \theta}{N^{4/3}} \sum_{i=1}^N v_i^2 \frac{1}{  ( \lambda_i - \gam1 + N^{-2/3} t )^3}
\end{align}
and conclude from the formula $\det \nabla^2 G_1 ( \gam1 , \gam1 ) = ( \del_z^2 G_1  - \del_z \del_w G_1 )( \del_z^2 G_1 + \del_z \del_w G_1 ) ( \gam1, \gam1)$ that
\beq
N^{-\eps_1/5} \leq N^{-4/3} \det \nabla^2 G_1 ( \gam1, \gam1 ) \leq N^{\eps_1/5},
\eeq
on $\Feps$.  The claim now follows from our estimates on the trace and determinant of the $2 \times 2$ matrix $\nabla^2 G_1 ( \gam1, \gam1)$. \qed 

\section{Proofs of auxiliary results in microscopic regime} \label{a:microscopic}

\subsection{Proof of Lemma \ref{lem:micro-1-f}} \label{a:micro-1-f}

 At $E=-1$ we obtain a quadratic equation for $\eta$ (the RHS simplifies to $\frac{\pi}{2}$ for $\eta >0$) and see that the only positive solution to \eqref{eqn:imp-eta} is given by
\beq
\etat(-1) = \frac{ \frac{\pi}{2} + \sqrt{\frac{ \pi^2}{4} + 4B(B+1)}}{2 (B+1)}.
\eeq
At $E = 0$ we see that $\etat(0) = 0$ is a solution to \eqref{eqn:imp-eta}.  Now, 
\beq
\frac{\mathrm{d}}{\mathrm{d} \eta} \left( \eta(B+1) - \frac{ \eta B}{ 1 + \eta^2 } - \arg (1 + \i \eta ) \right) = \frac{ (B+1) \eta^2}{  1 + \eta^2} + \frac{2 \eta^2 B}{ ( 1 + \eta^2)^2} > 0,
\eeq
and so this is the unique solution to \eqref{eqn:imp-eta}.  To simplify discussion, define
\beq
F_E ( \eta) := \eta \left( B+1 - \frac{B}{ (1+E)^2 + \eta^2 } \right)  - \arg \left( (1 + E ) + \i \eta \right)
\eeq
with the convention that $F_{-1} (0) = - \frac{ \pi}{2}$.  Note also $F'_{-1} ( \eta ) >0$.   From all of this, we see that on the vertical half-line in the complex plane given by
\beq
\{ z = -1 + \i \eta : \eta >0 \}
\eeq 
 that $F_{-1} (\eta ) <0$ for $\eta < \etat(-1)$ and $F_{-1} ( \eta ) >0$ for $\eta > \etat(-1)$.  Moreover, $F_0 ( \eta) >0$ for $\eta >0$.  Consider now the horizontal line segments $-1 < E < 0$ at fixed height $0 < \eta < \etat (-1)$.  We have,
\beq
\frac{\mathrm{d}}{\mathrm{d}E} F_E ( \eta ) = \frac{ 2 B \eta (1+E)}{ ( (1+E)^2 + \eta^2 )^2} + \frac{ \eta}{ (1+E)^2 + \eta^2 } > 0.
\eeq
So for fixed $0 < \eta < \etat (-1)$ we have that the function $E \mapsto F_E (\eta)$ is strictly increasing, and satisfies $F_{-1} ( \eta) <0$ and $F_0  (\eta) >0$, so for each $\eta$ in between $0 < \eta < \etat(0)$ there is some $-1 < E< 0$ s.t. $F_E ( \eta )=0$. We will see later that this function $\eta \mapsto E$ is invertible and so we can define $E \mapsto \etat (E)$ as its inverse. 

We now consider $E< -1$.  Define first $\eta_* (E)$ by,
\beq
\eta_* (E) := \inf\left\{ \eta>0 : (1+E)^2 + \eta^2 > B(1+B)^{-1} \right\}
\eeq
Note that $F_E ( \eta_* (E) ) <0$ whether or not $\eta_* (E)$ is $0$ or strictly positive.  Next,
\beq
F'_E ( \eta ) = (B+1) - \frac{B}{ (1+E)^2 + \eta^2 } + \frac{2 \eta^2 B}{ ( (1+E)^2 + \eta^2 )^2 } - \frac{1+E}{ (1+E)^2 + \eta^2 } 
\eeq
For $\eta > \eta_* (E)$, the sum of the first two terms is a positive quantity and so $F'_E ( \eta )>0$ for $\eta> \eta_* (E)$.  Moreover, $F_E ( \infty) = \infty$ so we find a unique solution in the domain $\eta > \eta_* (E)$ for each $E$.  Here, we see also by the implicit function theorem that $E \to \etat(E)$ is a smooth function.

We return now to the region $-1 < E < 0$, and consider
\beq
f' (u)  = B+1 - \frac{1}{1+u} - \frac{B}{(1+u)^2}.
\eeq
We claim that if $u$ lies on the contour that we have described, then $\Re[f'(u)] >0$.  In the region $-1 < E < 0$ this shows a strictly monotonic relation between the $E$ and $\etat(E)$ constructed above (due to the Cauchy-Riemann equations $\del_\eta \Im[f] = \Re[f']$), fulfilling our above promise to show that the relation $\eta \to E$ is invertible.  Moreover, this shows that $E \to \etat (E)$ is smooth.  Let $1+u = x + \i y$.  We have,
\beq
\Re[f' (u) ] = B+ 1 - \frac{x}{x^2+y^2} - \frac{B}{ x^2 + y^2} + 2 B \frac{ y^2}{ (x^2 + y^2 )^2}.
\eeq
Now, on the contour we have that $B+1-B/ (x^2+y^2) = y^{-1} \arctan(y/x)$.  Hence, on the contour we have
\begin{align} \label{eqn:aug30-1}
\Re[ f' (u) ] &= y^{-1} \arctan (y/x) - \frac{x}{x^2 + y^2 } + 2 B \frac{y^2}{ (x^2 + y^2 )^2} \nonumber \\
=& \frac{1}{y} \left( \arctan(y/x) - \frac{y/x}{1+(y/x)^2} \right) + 2 B \frac{y^2}{ (x^2 + y^2 )^2}.
\end{align}
We claim that for $w>0$, that the function $w \to \arctan(w) - w/(1+w^2)$ is strictly positive.  At $w=0$ the function is $0$. Its derivative is $2w^2 (1+w^2)^{-2} >0$, which proves this claim.  We therefore conclude that $\Re[f'(u)]>0$ on the contour. \qed

\subsection{Proof of Lemma \ref{lem:micro-2-f}} \label{a:micro-2-f}

We first examine the behavior of the contour near $E=0$.    For $|u| < \frac{1}{2}$, we can write
\beq
f(u) = \frac{1}{2}u^2(1+2B)- \frac{u^3}{6}(2+6B)  + f_1 (u)
\eeq
where $f_1(u)$ is a function satisfying $| \Im [f_1 (u)] | \leq C(1+B) \eta (|E|^3+|E|\eta^2 +  \eta^4)$.  Hence, we find the equation,
\beq
\frac{1}{3} \eta^2 (1 + 3 B) = - E (1+2B) + (1+3B) E^2 + (1+B) \O ( |E|^2 + \eta^4 ).
\eeq
Hence,
\beq
\etat(E) = \sqrt{ -E \frac{3+6B}{1+3B} } \left( 1+ \O(|E| ) \right)
\eeq
for $|E| \leq c$, some $c>0$.  

Due to the monotonicity of $E \to \etat(E)$ between $-1 <E < 0$ we see from the above that for any $c>0$, there is a $c' >0$ such that $\etat(E) >c'$ if $-1 < E < -c$.   Note also that in this range, $ \etat(E) \leq \eta (-1) \leq C$.  Now for $E <-1$ we have,
\beq
\eta \geq \frac{\pi}{2 (B+1) }.
\eeq
For an upper bound, we first have
\beq
\pi \geq \eta (1+B - \frac{B}{ \eta^2} )
\eeq
If $\eta^2 >10$, then we see that $\eta \leq C \pi / (1+B)$.  So, we have $\etat(E) \leq C$ for $E<-1$.  Hence, we have the lemma. \qed
\subsection{Proof of Lemma \ref{lem:micro-f-est-1}} \label{a:micro-f-est-1}

 Letting $1+u = x+ \i y$ we have as in \eqref{eqn:aug30-1},
\beq
\Re[f'(u)] = \frac{1}{y} \left( \arctan(y/x) - \frac{y/x}{1+(y/x)^2} \right) + \frac{2B y^2}{(x^2 + y^2)^2} .
\eeq
Now suppose that $-1 < E < -c$.  Then as above $\etat(E) > c'$ for some $c'>0$,  and so $y/x > c'$.  From the argument immediately below \eqref{eqn:aug30-1} we see that $y/x > c'$ implies that
\beq
\left( \arctan(y/x) - \frac{y/x}{1+(y/x)^2} \right) \geq c''
\eeq
for some $c''>0$.  Therefore, $\Re[f'(u)] \geq c''' (1+B)$ for some $c'''$.  For $E < -1$ we have,
\beq
\Re[f'(u) ] > \frac{\pi}{2y} > c_1
\eeq
for some $c_1 >0$.  Now,
\beq
\Im[f' (u) ] = \frac{y}{ x^2 + y^2} + \frac{B xy}{ (x^2+y^2)^2}
\eeq
and we see $|\Im [f'(u)] | \leq C(B+1)$ for $E < -c$.  Hence we conclude the proof of the lemma. \qed

\subsection{Proof of Lemma \ref{lem:bessel}} \label{a:bessel}

 We first make the substituion $z = (b/a)^{1/2} u$ to find,
\beq
\frac{1}{ 2 \pi \i} \int_{\tilGam} \exp\left[ a z + b z^{-1} - \alpha \log(z) \right] \d z  = \frac{1}{ 2 \pi \i }  \int_{\tilGam} \exp\left[ \frac{ \lambda}{2}(u + u^{-1} ) - \alpha \log(u) \right] \d u (b/a)^{(1-\alpha)/2}
\eeq
where $\lambda =2 (ab)^{1/2}$.  Now we break up $\tilGam$ into a circle $|z|=1$ and two horizontal segments at $\pm \i 0$ and $\Re[z] \leq -1$.  The line from $- \infty$ to $-1$ sitting in the lower half-plane at $-\i 0$ can be parameterized by $u= - \e^{t}$, $t>0$.    This portion of the integral contributes (dropping the $b/a$ factor)
\beq
\frac{1}{2 \pi \i } \int_0^\infty \exp\left[ - \lambda \cosh (t) +(1-\alpha)t + \alpha \i \pi \right] \d t 
\eeq
The contribution from the other line segment is identical except that the term $\alpha \i \pi$ becomes $- \alpha \i \pi$, due to the branch cut of the logarithm.  For the circle we parameterize it as $z = \e^{ \i t}$ for $- \pi \leq t \leq \pi$ and see that it contributes,
\beq
\frac{1}{2 \pi} \int_{- \pi}^{ \pi} \exp \left[ \lambda \cos ( \theta) + \i (1-\alpha) \theta \right] \d \theta= \frac{1}{ \pi} \int_0^\pi \exp[ \lambda \cos \theta] \cos((1-\alpha) \theta) \d \theta.
\eeq
We find first equality in the lemma after adding up all of the contributions.  
The second equality in the lemma, i.e., the integral representation for the Bessel function, is on page 181 of \cite{watson}.

We now turn to the proof of \eqref{eqn:aug31-1}. One can check that $I_{\beta} (x)$ satisfy,
\beq
x^2 I_\beta''(x) + x I_\beta' (x) -(x^2 + \beta^2) I_\beta (x) = 0.
\eeq
Making the substitution $\varphi (x) := \sqrt{x} I_{1/2} (x)$ we see that $\varphi''(x) = \varphi (x)$ and so $\varphi(x)$ is a linear combination of $\sinh(x)$ and $\cosh(x)$.  Since $I_{1/2} (0)$ is finite, we see that $\varphi(x)$ is a linear multiple of $\sinh(x)$.   From the equations, (see page 79 of \cite{watson}),
\begin{align}
I'_\beta ( \lambda) &= \frac{\beta}{ \lambda} I_\beta ( \lambda) +  I_{\beta+1} ( \lambda)\nonumber\\
I'_\beta ( \lambda) &= - \frac{ \beta}{ \lambda} I_\beta ( \lambda) +I_{\beta-1} ( \lambda) , \label{eqn:bessel-recursion}
\end{align}
we see that $\varphi' (x) = \sqrt{x} I_{-1/2} (x)$.  Then,
\begin{align}
\lim_{x \to 0} \sqrt{x} I_{-1/2} (x) &= \lim_{x \to 0} \frac{1}{\pi} \sqrt{x} \int_0^\infty \exp[ - x \cosh (t) + t/2] \d t \nonumber \\
&= \lim_{x \to 0} \frac{1}{\pi} \sqrt{x} \int_0^\infty \exp[ - x \cosh (t) ] 2 \sinh(t/2) \d t \nonumber \\
&= \lim_{x \to 0} \frac{1}{ \pi} \sqrt{x}  \int_0^\infty \exp[- 2 x \cosh(t/2)^2 ] 2 \sinh(t/2) \d t \nonumber\\
&= \lim_{ x \to 0 } \frac{1}{ \pi} \sqrt{x} \int_1^\infty 4 \exp[- 2 x t^2] \d t = \frac{\sqrt{2}}{ \sqrt{  \pi}}.
\end{align}
Therefore, we conclude \eqref{eqn:aug31-1}. \qed

\subsection{Proof of Proposition \ref{prop:micro-steepest-1}} \label{a:micro-steepest-1}

 The first estimate is from Proposition \ref{prop:micro-int}. We now prove the second. On the event $\Feps$, we have the estimate for $z \in \Gamk$: 
\beq
\left| \frac{1}{N^2} \sum_i \frac{1}{ | \lambda_i - z|^2}  \right| \leq N^{\eps_1/10}.
\eeq
 This estimate also holds  if $\eta \geq N^{-1}$. Hence, from Lemma \ref{lem:micro-G} we can replace the contour in the integral  by $\Gamk$ with only an error of size $\O(\e^{-N^c})$ some $c>0$.  For $z \in \Gamk$ we have the estimate, which holds on the event $\Feps$,
\begin{align}
\frac{N}{2} \left( G(z + \gamt) - G (\gamt)\right) &= \tilf (z) + \O(N^{2\kapp+\eps_1/10-2/3} ),
\end{align}
as well as,
\beq
\frac{1}{N^2} v^T \frac{1}{(M-\gamt - z)(M-\gamt -w) } v = \frac{v_1^2}{\cb^2 (1 + \frac{Nz}{ \cb} )( 1 + \frac{Nw}{\cb} )} + \O(N^{-2/3+\eps_1/10} ),
\eeq
and so,
\begin{align}
&\int_{\Gamk^2} \frac{}{N^2 } v^T \frac{1}{ (M-w)(M-z) } v \exp\left[ \frac{N}{2} (G(z) + G(w) - 2 G ( \gamt)) \right] \d z \d w \nonumber\\
= & \int_{\Gamk^2} \exp\left[ \tilf (z) + \tilf (w) \right] \frac{v_1^2}{\cb^2 (1 + \frac{Nz}{ \cb} )( 1 + \frac{Nw}{\cb} )} \d z \d w + \O ( N^{-2-2/3+4\kapp+\eps_1}).
\end{align}
Due to the second estimate of Lemma \ref{lem:micro-tilf} we may then turn $\Gamk$ into a keyhole contour by adding in the portion above and below negative real axis at only an exponentially small error.  The claim then follows from the substitution $Nz / \cb +1 = u$. \qed

\subsection{Proof of Proposition \ref{prop:aug-31}} \label{a:prop-aug-31}

Before embarking on the proof we state the following representation for the fourth moment of the overlap. It is proven by the same methodology as the proof of Proposition \ref{prop:moment-rep} in addition to some tedious calculations which we omit for brevity.
\bep \label{prop:rep-4m}
Denote the functions,
\beq
F_1 (z, w) = \frac{1}{ \beta^2} \sum_i \frac{1}{ ( \lambda_i - z ) ( \lambda_i  - w)} , \qquad F_2 (z, w) = \frac{ \lambda^2}{\beta^2} \sum_i \frac{ v_i^2}{ ( \lambda_i - z)( \lambda_i - w) }.
\eeq
The following holds.
\begin{align}
& \int ( \sigma^{(1)} \cdot \sigma^{(2)} )^4  \exp \left[ \frac{\beta}{2} \left( (\sigma^{(1)})^T M \sigma^{(1)} +\sigma^{(2)})^T M \sigma^{(2)}\right) + \lambda v^T ( \sigma^{(1)} + \sigma^{(2)}) \right] \d \sigma^{(1)} \d \sigma^{(2)} \nonumber \\
&= \frac{ \beta^2 N}{ (2 \pi \i )^2} \left( \frac{ 2 \pi}{ \beta} \right)^N \int \d z \d w \exp \left[ \frac{N}{2} ( G_o(z, v, \lambda, \beta) + G_o (w, v, \lambda, \beta) ) \right] \times\bigg\{ 6 \beta^{-2} F_{1zw} + 3 (F_1)^2 \nonumber \\
&+3\beta^{-2}( F_{2z}^2 + F_{2w}^2) +F_2^4 + 4! \beta^{-2} F_2 F_{2zw} + 6\beta^{-2} F_{2z} F_{2w} + 6 F_2^2 (- \beta^{-1} F_{2z} + -\beta^{-1} F_{2w} + F_1 )\nonumber \\
-& 6\beta^{-3} (F_{2zzw} + F_{2wwz} ) - 6\beta^{-1} F_1( F_{2z} + F_{2w} ) \bigg\}
\end{align}
Above, we have omited the arguments of the $F_i$ for brevity (they are always $F_i (z, w)$).  Furthermore, the notation $F_{1zw}$, etc., denote partial derivatives wrt $z$, $w$. 
\eep

We apply Proposition \ref{prop:rep-4m} to find the representation,
\beq
\frac{1}{N^4} \left\langle \left( \sigma^{(1)} \cdot \sigma^{(2) } \right) \right\rangle = D B^{-2}
\eeq
where,
\begin{align}
D &:= \frac{1}{ \beta^4} \int_{\Gamma^2} \d z \d w \exp\left[ \frac{N}{2} (G(z) + G(w)-2 G( \gamt) ) \right]  \bigg\{ 6 N^{-2} F_{1zw} + 3 (F_1)^2 \nonumber \\
&+3 N^{-2}( F_{2z}^2 + F_{2w}^2) +F_2^4 + 4! N^{-2} F_2 F_{2zw} + 6N^{-2} F_{2z} F_{2w} + 6 F_2^2 (- N^{-1} F_{2z} -N^{-1} F_{2w} + F_1 ) \nonumber \\
-& 6N^{-3} (F_{2zzw} + F_{2wwz} ) - 6N^{-1} F_1( F_{2z} + F_{2w} ) \bigg\}
\end{align}
where,
\beq
F_1 (z, w) = \frac{1}{N^2} \sum_i \frac{1}{ ( \lambda_i - z)(\lambda_i - w ) }, \qquad F_2 (z, w) = \frac{ \theta}{N^2} \sum_i \frac{ v_i^2}{ ( \lambda_i - z ) ( \lambda_i - w) }.
\eeq
In the definition of $D$ above we have suppressed the arguments of $F_1$ and $F_2$ as they are all just the integration variables $(z, w)$.   It is no problem to argue as in our calculations of the first and second moments of the overlap to move the contour $\Gamma$ to $\Gamk$, and then expand the function $G(z) - G ( \gamt)$ that appears in the exponential around $\tilf(z)$ along $\Gamk$ on the event $\Feps$.  Similarly, for $z, w \in \Gamk$ we use the following estimates which hold on the event $\Feps$
\begin{align}
F_1 (\gamt+z, \gamt + w) &= \frac{1}{\cb^2} \frac{1}{ (1+ \frac{Nz}{ \cb} )(1+ \frac{Nw}{\cb} ) } + \O (N^{-2/3+\eps_1/10} ) \nonumber \\
N^{-2} F_{1zw} ( \gamt+z, \gamt+w) &= \frac{1}{ \cb^4} \frac{1}{ (1+ \frac{Nz}{ \cb} )^2(1+ \frac{Nw}{\cb} )^2 } + \O (N^{-4/3+\eps_1/10} ) \nonumber \\
F_2( \gamt+z, \gamt+w)  &= \frac{ \theta v_1^2}{ \cb^2} \frac{1}{ (1+ \frac{Nz}{ \cb} )(1+ \frac{Nw}{\cb} ) } +\O (N^{-2/3+\eps_1/10} ) \nonumber \\
N^{-1} F_{2z} ( \gamt+z, \gamt+w) &=-  \frac{ \theta v_1^2}{ \cb^3} \frac{1}{ (1+ \frac{Nz}{ \cb} )^2(1+ \frac{Nw}{\cb} ) } + \O (N^{-1+\eps_1/10} ) \nonumber\\
N^{-2} F_{2zw}( \gamt+z, \gamt+w)  &=  \frac{ \theta v_1^2}{ \cb^4} \frac{1}{ (1+ \frac{Nz}{ \cb} )^2(1+ \frac{Nw}{\cb} )^2 } + \O (N^{-4/3+\eps_1/10} ) \nonumber\\
N^{-3} F_{2zzw} ( \gamt+z, \gamt+w) &= -2 \frac{ \theta v_1^2}{ \cb^5} \frac{1}{ (1+ \frac{Nz}{ \cb} )^3(1+ \frac{Nw}{\cb} )^2 } + \O (N^{-5/3+\eps_1/10} ) .
\end{align}
Changing the contour from $\Gamk$ to the keyhole $\tilGam$ after making the same change of variables $1+\frac{Nz}{\cb} = u$ incurs the same exponential error as before.  Calculating all of the contributions from $F_1$, $F_2$ and their derivatives yields the following, where we drop the error $\O(N^{-2/3-2+4\kapp+\eps_1})$ for brevity:
\begin{align}
&\frac{1}{ ( 2 \pi \i )^2} D = \frac{ \e^{-2a-2b}}{\beta^4 N^2 \cb^2} \bigg\{ 9 I_{3/2} ( \lambda)^2  (a/b)^{3/2} + (2b)^4 I_{7/2} ( \lambda)^2 (a/b)^{7/2} + (4!+6+6) (2b)^2 I_{5/2} ( \lambda)^2 (a/b)^{5/2} \nonumber \\
+&12(2b)^3 I_{5/2} (\lambda) I_{7/2} (\lambda)(a/b)^3+12(2b)(2+1) I_{3/2} ( \lambda) I_{5/2} ( \lambda) (a/b)^2 + 6 I_{7/2} ( \lambda) I_{3/2}  ( \lambda) (2b)^2 (a/b)^{5/2} \bigg\} \nonumber \\
=& \frac{ a^{3/2} \e^{-2a-2b}}{\beta^4 N^2 \cb^2 b^{3/2} }
\times  \bigg\{ 9 I_{3/2} ( \lambda)^2  +  I_{7/2} ( \lambda)^2 \lambda^4 +36 I_{5/2} ( \lambda)^2 \lambda^2 \nonumber\\
+& 12 I_{5/2} ( \lambda) I_{7/2} ( \lambda) \lambda^3 + 36 I_{3/2} ( \lambda) I_{5/2} ( \lambda) \lambda + 6 I_{7/2} ( \lambda) I_{3/2}  ( \lambda) \lambda^2 \bigg\}. \label{eqn:micro-D}
\end{align}
Using
\begin{align}
I_{3/2} ( \lambda) &= I_{-1/2} ( \lambda ) - \frac{ I_{1/2} ( \lambda)}{ \lambda} \nonumber\\
 I_{5/2} ( \lambda) &= I_{1/2} ( \lambda) (1 + \frac{3}{ \lambda^2} ) -\frac{3}{\lambda} I_{-1/2} ( \lambda) \nonumber\\
I_{7/2} ( \lambda) &= I_{-1/2} ( \lambda) (1 + \frac{15}{ \lambda^2} ) - I_{1/2} ( \lambda) ( \frac{6}{ \lambda} + \frac{15}{ \lambda^3} )
\end{align}
we find for the term in braces in the last line of \eqref{eqn:micro-D}, 
\begin{align}
& 9 I_{3/2} ( \lambda)^2  +  I_{7/2} ( \lambda)^2 \lambda^4 +36 I_{5/2} ( \lambda)^2 \lambda^2 + 12 I_{5/2} ( \lambda) I_{7/2} ( \lambda) \lambda^3 + 36 I_{3/2} ( \lambda) I_{5/2} ( \lambda) \lambda + 6 I_{7/2} ( \lambda) I_{3/2}  ( \lambda) \lambda^2 \nonumber \\
= & \lambda^4 I_{-1/2} ( \lambda)^2.
\end{align}
Hence, on $\Feps$ we have for sufficiently small $\kapp$ and $\eps_1 >0$,
\beq
\frac{1}{N^4} \langle ( \sigma^{(1)} \cdot \sigma^{(2)} )^4 \rangle = \left( \frac{ \beta + \tilm ( \gamt ) }{ \beta} \right)^4 + \O(N^{-2/3+10\kapp+ 2\eps_1} ) .
\eeq
This yields the claim. \qed

\section{Proofs of Theorem \ref{thm:compare} and Lemma \ref{lem:HViso}} \label{a:compare}

We begin by proving Lemma \ref{lem:HViso}.  We start with the resolvent identity,
\beq \label{eqn:resolv}
\frac{1}{M-z} - \frac{1}{H-z} = \frac{1}{M-z} \sum_{k=1}^m (V (M -z)^{-1} )^k   + \frac{1}{H-z} (V (M -z ))^{-(m+1)}.
\eeq
Denote,
\beq
A_k := \frac{1}{ M -z } \left( V \frac{1}{ M - z } \right)^k, \qquad R(z) := \frac{1}{M-z}.
\eeq
We first prove the following lemma.
\bel \label{lem:resolv1}
Let $C>0$ be a constant.  On the event,
\beq
\max_{i, j} |R_{ij} | \leq C,
\eeq
we have for $k \geq 2$ and even $p$,
\beq
\ee_V \left| v^T A_k v \right|^p \leq C(k, p) \left( \frac{1}{N^p} + \sup_{i} | (R v)_i - \msc(z) v_i |^{2p} \right)
\eeq
for any unit vector $v$, where $\ee_V$ denotes expecation over $V$. 
\eel
\proof We write,
\beq
v^T A_k v = \sum_{j_1, \dots j_k } (R v)_{j_1} V_{j_1 j_1} R_{j_1 j_2} V_{j_2 j_2} \dots R_{j_{k-1} j_k }  V_{j_k j_k }(R v)_{j_k}.
\eeq
For even $p$, 
\beq
\ee_V \left| v^T A_k v \right|^p = \sum_{ \ulj} (R^\# v)_{j_1} (R^\# v)_{j_k} (R^\# v)_{j_{k+1} } \dots (R^\# v)_{j_{kp }}  M ( \ulj ) \ee_V [ V_{j_1} V_{j_2} \dots V_{j_{kp}} ],
\eeq
where $R^\#$ denotes $R$ or $R^*$ where necessary - the distinction will make no difference for us. The sum is over $kp$-tuples of indices $\ulj$.   The term $M ( \ulj)$ is a monomial in the matrix elements $R_{ij}$ or $\bar{R}_{ij}$.  Its specific form is unimportant.  Note that importantly since $k \geq 2$, the index $j_1$ is different from $j_k$.

 We now rewrite the summation over $\ulj$ as a sum over partitions induced by the coincidences of $\ulj$, that is
\beq
\sum_{\ulj} = \sum_{ \P} \sum_{ \ulj \in \P}.
\eeq
The first sum is over partitions $\P$ on $kp$ elements, and the second sum means over all multi-indices $\ulj$ so that $j_a = j_b$ whenever $a$ and $b$ are in the same block of $\P$ and $j_a \neq j_b$ if $a$ and $b$ are in different blocks of $\P$.  Now, note that the expectation $\ee_V [ V_{j_1} \dots V_{j_{kp}}]$ vanishes unless the partition induced by $\ulj$ has every block at least size $2$.  We denote these partitions by $\P_2$.  Using the assumption that the resolvent entries are all bounded, we have
\beq
\ee_V \left| v^T A_k v \right|^p \leq C(k, p)  \sum_{\P \in \P_2} \sum_{ \ulj \in \P} | (Rv)_{j_1} (R v)_{j_k} \dots (R v)_{j_{kp } } | \frac{1}{N^{kp/2} }.
\eeq
Fix a partition $\P \in \P_2$.  Assume that it has $\ell$ blocks, each of size $s_k \geq 2$.  Let $a_k$ be the number of powers of $(R v)$ that get assigned to the $k$th block, so that
\beq
\sum_{ \ulj \in \P} | (Rv)_{j_1} (R v)_{j_k} \dots (R v)_{j_{kp } } | \frac{1}{N^{kp/2} } \leq \prod_{j =1 }^\ell \left( \frac{1}{N^{s_j / 2} } \sum_{i} | (R v)_i |^{a_j} \right).
\eeq
Consider the $j$th term in the product on the RHS.  If $a_j = 0$, then it is bounded by $1$ because $s_j \geq 2$.  If $a_j =1$, then it is bounded by 
\beq
\frac{1}{N^{s_j/2} } \sum_i | (Rv)_i | \leq \frac{1}{N} \sum_i | (Rv)_i | \leq \frac{C}{N} \|v\|_1 + \max_i | (Rv)_i - \msc v_i | \leq \frac{C}{N^{1/2} } + \max_i | (Rv)_i - \msc v_i |
\eeq
We now turn to the case $a_j \geq 2$. Since the $k$ in $A_k$ is $k \geq 2$  we have $s_j \geq a_j$ (i.e., each index in $\ulj$ carries at most one power of $(Rv)_{j}$ and so the size of a partition block must be at least the number of powers of $(Rv)_j$ assigned to it). Therefore, for $a_j \geq 2$, we have, 
\beq
\frac{1}{N^{s_j/2}} \sum_i | (Rv)_i |^{a_j} \leq \frac{1}{ N^{a_j /2 } } | \msc v_i + (Rv)_i - \msc v_i |^{a_j} \leq C \frac{ \|v\|_2^2}{ N^{a_j/2}} + C \max_i | (Rv)_i  - \msc v_i |^{a_j}.
\eeq
Therefore, the inequality
\beq
\frac{1}{N^{s_j/2}} \sum_i | (Rv)_i |^{a_j}  \leq \frac{C}{ N^{a_j/2}} + C \max_i | (Rv)_i  - \msc v_i |^{a_j}.
\eeq
holds no matter the value of $a_j \geq 0$. The claim then follows  after noting that $\sum_{j} a_j = 2 p$. \qed

\vspace{5 pt}

\noindent{\bf Proof of Lemma \ref{lem:HViso}}.  We use the resolvent identity \eqref{eqn:resolv},
\beq
v^T \frac{1}{M-z} v - v^T \frac{1}{H-z} v = v^T A_1 v + \sum_{k=2}^m v^T A_k v + v^T \frac{1}{H-z} ( V (M-z)^{-1} )^{m+1} v.
\eeq
The term on the RHS can be taken to be less than $N^{-100}$ with overwhelming probability using Lemma A.4 of \cite{SSK-ls} by taking $m$ sufficiently large.  The terms involving $A_{k}$ for $2 \leq k \leq m$ are handled using Lemma \ref{lem:resolv1} and the estimates of Theorem \ref{thm:iso}.  Conditionally on $M$, the term $v^T A_1 v$ is a Gaussian with variance less than
\beq
\frac{1}{N} \sum_{i} | (Rv)_i |^4 \leq \frac{C}{N} + \sup_i | (Rv)_i - \msc v_i |^4.
\eeq
The claim follows, again using Theorem \ref{thm:iso}.  \qed

\vspace{5 pt}

\noindent{\bf Proof of Theorem \ref{thm:compare}}.  Equation \eqref{eqn:compare1} was proven in \cite{SSK-ls}, so it remains to prove \eqref{eqn:compare2}.  Fix a small $0 < \delta_{lr} \leq 0.1$ and small $\eps >0$, with $\eps < 0.1$.  We may assume that $| \lambda_i (M) - \lambda_i (H) | \leq N^{-1+\eps}$ for $i=1, 2$.  Assume that the level repulsion events of Theorem \ref{thm:lr} hold with $s = N^{-\delta_{lr}}$.  Let $\delta_1 >0$ with $\delta_1 < 0.1$, and set $\eta_1 = N^{-2/3-\delta_1}$ and let $\eta_2 = N^{-2/3-\delta_{lr}}/2$.  We denote by $\Gamma$ the contour in $\cc$ that is a rectangle with sides parallel to the real and imaginary axes, symmetric across the real axis, centered at the point $\lambda_1 (M)$ and horizontal side length $2 \eta_2$ and vertical side length $\eta_1$.  

Due to our assumptions, we have that only $\lambda_1 (M)$ and $\lambda_1(H)$ are inside the contour $\Gamma$, and that $\lambda_2 (M)$ and $\lambda_2 (H)$ are at distance $N^{-2/3-\delta_{lr}}/2$ at least from $\Gamma$.  It follows that, on the event described above that,
\beq
(v^T u_1(M) )^2 - (v^T u_1 (H) )^2 = \frac{1}{ 2 \pi \i } \int_\Gamma v^T (M - z)^{-1} v - v^T (H-z)^{-1} v \d z.
\eeq
Fix a small $\delta_c >0$ with $\delta_c  \leq 0.1$.  We estimate in the above integral the contribution from $| \Im [z] | \leq N^{\delta_c-1}$.  Due to the orientation of the integral, we may estimate the contribution as
\beq
\int_{ |y| \leq N^{\delta_c-1} } | v^T ( A - (\lambda_1 (M) + \eta_2 + \i y ) )^{-1} v - v^T (A - ( \lambda_1 (M) - \eta_2 + \i y ) )^{-1} v | \d y
\eeq
for $A = M, H$.  We bound first the contribution coming from $\lambda_1 (A)$; recalling that $| \lambda_1 (H) - \lambda_1 (M) | \ll \eta_2$, we have with overwhelming probability,
\beq
\int_{ |y| \leq N^{\delta_c -1 } } \frac{ ( v^T u_1 (A) )^2}{ | \lambda_1 (A) - ( \lambda_1 (M) \pm \eta_2 + \i y ) |} \d y \leq \frac{N^{\eps}}{N} \int_{ |y| \leq N^{\delta_c -1} } \frac{C}{ \eta_2} \d y \leq \frac{N^{ 2 \eps} N^{\delta_c}}{ N^2 \eta_2 }. \label{eqn:compareest1}
\eeq
We denote $\eta_{lr} = N^{-2/3-\delta_{lr}}$.  For the contribution from the remaining eigenvalues, we have, on the level repulsion event with overwhelming probability,
\begin{align}
 & \int_{ |y| \leq N^{ \delta_c -1 }} \sum_{i=2}^N (v^T u_i (A) )^2 \left| \frac{1}{  \lambda_i (A) - ( \lambda_i (M) - \eta_2 + \i y ) } - \frac{1}{ \lambda_i (A) - ( \lambda_i (M) + \eta_2 + \i y ) } \right| \d y \notag \\
\leq & \frac{N^{ 2 \eps } \eta_2 }{N} \int_{ |y| \leq N^{\delta_c-1}} \sum_{i=2}^N \frac{1}{ ( \lambda_i (A) - \lambda_1 (A) )^2 + ( \eta_{lr} )^2 } \d y \notag \\
\leq & N^{2 \eps} \eta_2  N^{\delta_c-1} \sup_{ |E-2| \leq N^{\eps-2/3} } (\eta_{lr} )^{-1} \Im [ N^{-1} \tr (A -(E + \i \eta_{lr} ) )^{-1} ] \notag  \\
\leq & \frac{N^{3 \eps } N^{\delta_c} \eta_2 }{N^2 \eta_{lr}^2} + N^{-1/3+3 \eps} \frac{ N^{\delta_c} \eta_2}{N \eta_{lr} } \label{eqn:compareest2}
\end{align}
where we applied Theorem \ref{thm:locallaw} in the final inequality and bounded $\Im[ \msc (z) ] \leq N^{-1/3+\eps}$ for $z$ in the indicated region. 
For the remaining portion of the vertical segments of $\Gamma$, we instead use \eqref{eqn:isocompare1}.  The error we get is 
\begin{align}
& \int_{\eta_1 > |y| > N^{ \delta_c -1 } } N^{\eps} \left( \frac{1}{ \sqrt{N} } + \frac{1}{ N^2 y^2} + \frac{y^{1/2} + N^{-1/3+\eps}}{N y } + \frac{1}{N \sqrt{ y + N^{-2/3+\eps} } } \right) \d y \notag \\
\leq & N^{4 \eps} \left(  \frac{\eta_1}{N^{1/2}} + \frac{1}{N^{1 + \delta_c }} + \frac{\eta_1^{1/2}}{N} + N^{-1-1/3} \right) \label{eqn:compareest3}
\end{align}
The contribution of the horizontal segments is bounded by,
\begin{align}
& \int_{ |x - \lambda_1 (M) | \leq \eta_2 } N^{\eps} \left(  \frac{1}{ \sqrt{N} } + \frac{1}{ N^2 \eta_1^2} + \frac{\eta_1^{1/2} + N^{-1/3+\eps}}{N \eta_1 } + \frac{1}{N \sqrt{ \eta_1 + N^{-2/3+\eps} } } \right) \notag \\
\leq & N^{2 \eps } \left( \frac{\eta_2}{ \sqrt{N}}   + \frac{\eta_2}{N^2 \eta_1^2 } + \frac{ \eta_2}{ N \sqrt{ \eta_1}} + \frac{N^{-1/3} \eta_2}{N \eta_1} \right) \label{eqn:compareest4}
\end{align}
If we choose, e.g., $\eta_{lr} = N^{-2/3-0.1}$,  $\eta_1 = \eta_2 = \eta_{lr}/2$, $\eps = 10^{-10}$, and say $\delta_c = 0.01$, then all of the errors \eqref{eqn:compareest1}, \eqref{eqn:compareest2}, \eqref{eqn:compareest3} and \eqref{eqn:compareest4} are seen to be $\O (N^{-1-c})$ for some $c >0$. \qed

\section{Isotropic CLT for zero-diagonal GOE} \label{a:iso}

In this section we recall that $M$ is a matrix from the GOE with a zero diagonal.  Let $v$ be a unit vector, and $\gamma > 2$.  Define,
\beq
\kappa := \gamma -2.
\eeq
We prove the following theorem. Let $R(z) := (M - z)^{-1}$. 
\bet \label{thm:iso-clt}
Fix $\eps >0$, and let $C \geq \kappa \geq N^{\eps-2/3}$.  Then,
\beq
V_N^{-1/2} N^{1/2} \kappa^{1/4} \left( v^T R( \gamma ) v -  \msc ( \gamma ) \right) 
\eeq
converges to a standard normal random variable.  Here, 
\beq
V_N := \frac{ \gamma + \sqrt{ \gamma^2 - 4 } }{ \sqrt{ \gamma+2} } \msc^4 \left( \msc^2 + (1 - \|v\|_4^4 ) ( 1 - \msc^2 ) \right)
\eeq
satisfies $c \leq V_N \leq C$.  If $\kappa \to 0$ then $V_N \to 1$.
\eet
\remark In the following proof we work with expectations of matrix elements of the resolvent on the real line, $(M - E)^{-1}$, with $|E| > 2+N^{-2/3+\eps}$, some $\eps >0$.   Due to integrability concerns, one should instead work with the regularization $\Re[ ( M - (E + \i \eta ) )^{-1}]$ for, e.g., $\eta = N^{-100}$.  With overwhelming probability, the difference for any matrix element between these two quantities is $\O ( N^{-90})$.  For notational convenience we omit this regularization in the proof below, but it is elementary to restore it and check that the proof goes through. \qed

\vspace{5 pt}

\proof Define the characteristic function,
\beq
\psi ( \lambda ) = \ee[ e ( \lambda ) ], \qquad e ( \lambda) := \exp \left[ \i \lambda N^{1/2} \kappa^{1/4}  v^T R^\circ ( \gamma ) v \right],
\eeq
where we introduced the notation $X^\circ := X - \ee[X]$ for any random variable.  We simplify notation and write $R = R(\gamma)$. 
We apply Stein's method and calculate,
\beq
\psi' ( \lambda ) = \i N^{1/2} \kappa^{1/4} \sum_{i, j} v_i v_j \ee[ e ( \lambda ) R_{ij}^\circ  ].
\eeq
From the matrix identity $R(M- \gamma ) = \1$ and Gaussian integration by parts,
\begin{align}
\gamma \ee[ e ( \lambda ) R_{ij}^\circ ] &= \sum_{a \neq j } \ee[ e ( \lambda ) ( R_{ia} M_{aj} - \ee[ R_{ia} M_{aj} ] ) ] \notag \\
&= \frac{1}{N} \sum_{a \neq j} \ee[ ( \partial_{aj} e ( \lambda ) ) R_{ia} ] \notag \\
& - \frac{1}{N} \sum_{a \neq j } \ee[ e ( \lambda ) ( R_{ia} R_{ja} + R_{ij} R_{aa} )^\circ ].
\end{align}
We begin with,
\begin{align}
& \sum_{i, j} v_i v_j N^{1/2} \kappa^{1/4} \frac{1}{N} \sum_{a \neq j } \ee[ e ( \lambda ) ( R_{ia} R_{ja} )^\circ ] \notag\\
=& \kappa^{1/4} N^{-1/2} \ee[ e ( \lambda) (v^T R^2 v)^\circ ] - \kappa^{1/4} N^{-1/2} \sum_{i, j} \ee[ e ( \lambda) v_i ( R_{ij} R_{jj} )^\circ v_j ].
\end{align}
From Theorem \ref{thm:iso} and the Cauchy integral formula, the first term on the RHS is $N^{\eps} \O ( N^{-1/2} \kappa^{-1/4} )$ for any $\eps >0$.  For the second term,
\begin{align}
 \kappa^{1/4} N^{-1/2} \sum_{i, j} \ee[ e ( \lambda) v_i ( R_{ij} R_{jj} )^\circ v_j ] &= \kappa^{1/4} N^{-1/2} \msc \ee[ e ( \lambda) (v^T R v )^\circ ] \notag \\
+& \kappa^{1/4} N^{-1/2} \sum_{i, j} \ee[ e ( \lambda) v_i ( R_{ij} (R_{jj}-\msc) )^\circ v_j ] \notag \\
&= \O (N^{\eps-1} ) + \O ( N^{\eps} N^{-1/2} \kappa^{-1/4} )
\end{align}
where for the terms on the second line we used $|R_{ij} (R_{jj} - \msc ) | \leq N^{\eps} (N \sqrt{ \kappa } )^{-1} + \delta_{ij} N^{\eps}N^{-1/2} \kappa^{-1/4}$ with overwhelming probability,  as well as the fact that $\|v\|_1 \leq N^{1/2}$ for the $i \neq j$ terms and  $\|v\|_2 \leq 1$ for the $i=j$ terms.  The next term we handle is,
\begin{align}
\kappa^{1/4} N^{-1/2} \sum_{i, j} v_i v_j \sum_{ a \neq j } \ee[ e ( \lambda ) ( R_{ij} R_{aa} )^\circ ] &= \msc \kappa^{1/4} N^{-1/2} \sum_{i, j} v_i v_j \sum_{a \neq j } \ee[ e ( \lambda ) (R_{ij} )^\circ ] \notag \\
&+ \kappa^{1/4} N^{-1/2} \sum_a \ee[ e ( \lambda ) ( v^T R v  ( R_{aa} - \msc ) )^\circ ] \notag \\
&- \kappa^{1/4} N^{-1/2} \sum_{i, j} v_i v_j \ee[ e ( \lambda ) ( R_{ij} (R_{jj} - \msc ) )^\circ ].
\end{align}
The term on the last line appeared above and was shown to be $\O (N^{\eps} N^{-1/2} \kappa^{-1/4} )$.  The term on the second line equals
\beq
\kappa^{1/4} N^{-1/2} \sum_a \ee[ e ( \lambda ) ( v^T R v  ( R_{aa} - \msc ) )^\circ ] = \kappa^{1/4} N^{1/2} \ee[ e ( \lambda ) ( v^T R v (m_N - \msc ) )^\circ ] = \O ( N^{\eps} N^{-1/2} \kappa^{-3/4} ) ,
\eeq
using Theorem \ref{thm:locallaw} and Theorem \ref{thm:iso}. 
The term on the first line equals,
\beq
 \msc N^{-1/2} \kappa^{1/4} \sum_{i, j} v_i v_j \sum_{a \neq j } \ee[ e ( \lambda ) (R_{ij} )^\circ ] = \msc \kappa^{1/4} N^{1/2} \ee[ e ( \lambda ) (v^T R v)^\circ ] + \O( N^{-1/2}). 
\eeq
We observe,
\beq
\partial_{aj} e ( \lambda ) =- 2 \kappa^{1/4} N^{1/2} \i \lambda e ( \lambda ) (R v)_a (R v)_j,
\eeq
and so 
\begin{align}
-N^{-1/2} \kappa^{1/4} \sum_{i, j} v_i v_j \sum_{a \neq j } \ee[ ( \partial_{aj} e ( \lambda ) R_{ia} ] &= 2 \kappa^{1/2} \i \lambda \sum_{a \neq j } \ee[ e ( \lambda ) (R v )_a^2 v_j (R v)_j ]  \notag \\
&= 2 \i \lambda \kappa^{1/2} \ee[ e ( \lambda) v^T R^2 v v^T R v ] - 2 \i \lambda \kappa^{1/2} \sum_{j} \ee[ e ( \lambda ) v_j (R v)_j^3].
\end{align}
The first term is,
\beq
2 \i \lambda \kappa^{1/2} \ee[ e ( \lambda) v^T R^2 v v^T R v ] = 2 \i \lambda \kappa^{1/2} \msc' \msc \psi ( \lambda ) + \O( | \lambda| N^{\eps} N^{-1/2} \kappa^{-3/4} ).
\eeq
For the other term we note that with overwhelming probability,
\begin{align}
v_j (R v)_j^3 &= v_j ( (R v)_j - v_j \msc  + v_j \msc )^3 \notag\\
&= \msc^3 v_j^4 + \msc^2 v_j^3 \O (N^\eps N^{-1/2} \kappa^{-1/4} ) + \O ( v_j^2 N^{-1} \kappa^{-1/2} + v_j N^{-3/2} \kappa^{-3/4} ),
\end{align}
and so 
\beq
2 \i \lambda \kappa^{1/2} \sum_{j} \ee[ e ( \lambda ) v_j (R v)_j^3] = 2 \i \lambda \kappa^{1/2} \msc^3 \|v\|_4^4 \psi ( \lambda ) + \O( |\lambda| N^{-1/2} ).
\eeq
From the above calculations we see that,
\beq \label{eqn:stein-1}
\gamma \psi' ( \lambda ) = - \msc \psi' ( \lambda ) + (2 \kappa^{1/2} \msc' \msc - 2 \|v\|_4^4 \msc^3 \kappa^{1/2} ) \lambda  \psi ( \lambda ) + \O ( (1 + | \lambda | )N^\eps N^{-1/2} \kappa^{-3/4} ).
\eeq
We calculate the quantity which will be seen to be the variance,
\beq
\frac{- 2 \kappa^{1/2} \msc}{ \gamma + \msc  } \left( \msc' - \|v\|_4^4 \msc^2 \right) =\left( \frac{ \kappa^{1/2}}{ 1 - \msc^2} \right) 2 \msc^4 \left(  \msc^2 + (1 - \|v\|_4^4 ) ( 1 - \msc^2 ) \right)
\eeq
where we used \eqref{eqn:quadratic-msc}. 
A calculation using $2\msc (z)= -z + \sqrt{4 -\gamma^2 }$ gives
\beq
\left( \frac{ \kappa^{1/2}}{ 1 - \msc^2} \right) = \frac{1}{2} \frac{ \gamma + \sqrt{ \gamma^2-4} }{ \sqrt{ \gamma + 2 } } \asymp 1.
\eeq
Noting that $1 - \msc^2 \geq 0$ we see that
\beq
V_N = \frac{ \gamma + \sqrt{ \gamma^2 - 4 } }{ \sqrt{ \gamma+2} } \msc^4 \left( \msc^2 + (1 - \|v\|_4^4 ) ( 1 - \msc^2 ) \right)
\eeq
is bounded below and above, $c \leq V_N \leq C$, uniformly in $N$. Therefore, re-arranging and integrating \eqref{eqn:stein-1} we see that
\beq
\psi ( \lambda ) = \exp \left( - \frac{ \lambda^2}{2} V_N \right) + \O \left( (1+ | \lambda | ) N^{\eps} N^{-1/2} \kappa^{-3/4} \right)
\eeq
  This proves that the centered random variable $v^T R^\circ v$ converges to a standard normal random variable after the stated normalization.

Next, we calculate the expectation.  We have 
\begin{align}
\gamma \kappa^{1/4} N^{1/2} \ee[ v^T R v ] = - \kappa^{1/4} N^{1/2} + \sum_{i, j} v_i v_j N^{1/2} \kappa^{1/4} \sum_{a \neq j } \ee[ R_{ia} M_{aj } ].
\end{align}
The second term equals, by Gaussian integration by parts,
\begin{align}
\sum_{i, j} v_i v_j N^{1/2} \kappa^{1/4} \sum_{a \neq j } \ee[ R_{ia} M_{aj } ] &= - \sum_{i, j} v_i v_j  N^{-1/2} \kappa^{1/4} \sum_{a \neq j } \ee[ R_{ia} R_{ja} + R_{ij} R_{aa} ] \notag \\
&=- N^{-1/2} \kappa^{1/4} \ee[ (v^T R^2 v)] - N^{1/2} \kappa^{1/4} \ee[ (v^T R v) m_N ] \notag \\
&+ 2\sum_{i, j} v_i v_j N^{-1/2} \kappa^{1/4} \ee[ R_{ij} R_{jj} ]
\end{align}
The last line is,
\beq
\sum_{i, j} v_i v_j N^{-1/2} \kappa^{1/4} \ee[ R_{ij} R_{jj} ] = N^{-1/2} \kappa^{1/4} \msc \ee[ v^T R v ] + \sum_{i, j} v_i v_j N^{-1/2} \kappa^{1/4} \ee[ R_{ij} (R_{jj}-\msc ) ]
\eeq
These terms are $\O (N^{\eps} N^{-1/2} \kappa^{-1/4} )$.  Next,
\beq
N^{1/2} \kappa^{1/4} \ee[ v^T R v m_N ] = \msc N^{1/2} \kappa^{1/4} \ee[ v^T R v ] + \O ( N^{-1/2} \kappa^{-3/4} ).
\eeq
Finally,
\beq
N^{-1/2} \kappa^{1/4} \msc \ee[ v^T R^2 v ] = \O ( N^{-1/2} \kappa^{-1/4} ).
\eeq
Hence,
\beq
\kappa^{1/4} N^{1/2} \ee[ v^T R v ] = \kappa^{1/4} N^{1/2}  \msc + \O (N^{-1/2} \kappa^{-3/4} ).
\eeq
This yields the claim. \qed

\section{Conditional probability statement}

We require the following elementary lemma. 
\bel \label{lem:cond}
Let $X$ and $Z$ be random variables and $\G$ a sigma-algebra.  Let $F$ be a bounded Lipschitz function with Lipschitz constant $\|F\|_{\mathcal{L}}$ and let $\eps_1, \eps_2 >0$ be constants so that,
\beq
\P[ |X-Z| > \eps_1 ] \leq \eps_2.
\eeq
Then,
\beq
\ee\left[ \left| \ee[ F(X) \mid \G ] - \ee[ F(Z) \mid \G ] \right| \right] \leq 2 \|F\|_{\mathcal{L}} ( \eps_1 + \eps_2),
\eeq
and
\beq
\pp \left[ \left| \ee[ F(X) \mid \G ] - \ee[ F(Z) \mid \G ] \right|  > \eta \right] \leq \frac{2 \|F\|_{\mathcal{L}} ( \eps_1 + \eps_2)}{ \eta}.
\eeq
\eel
\proof We have,
\begin{align}
\ee[ F(X) \mid \G ] - \ee[ F(Z) \mid \G ] = \ee[ (F (X) - F(Z) ) \1_{ \{ |X-Z| > \eps_1 \} } \mid \G] + \ee[ (F (X) - F(Z) ) \1_{ \{ |X-Z| < \eps_1 \} }\mid  \G].
\end{align}
Then by conditional Jensen's inequality, 
\begin{align}
&\ee \left|  \ee[ (F (X) - F(Z) ) \1_{ \{ |X-Z| > \eps_1 \} }\mid \G] + \ee[ (F (X) - F(Z) ) \1_{ \{ |X-Z| < \eps_1 \} }\mid  \G] \right| \nonumber \\ 
\leq& 2 \|F\|_\infty \pp[ |X- Z| > \eps_1 ] + \|F\|_{\mathcal{L}} \eps_1,
\end{align}
which yields the first claim.  The second is of course Markov's inequality.
\qed

\bibliography{ssk_bib}
\bibliographystyle{abbrv}

\end{document}